\documentclass[11pt,twoside]{article}
\usepackage{latexsym}
\usepackage{amssymb,amsbsy,amsmath,amsfonts,amssymb,amscd}
\usepackage{mathrsfs}
\usepackage{epsfig, graphicx}
\usepackage{graphics,epsfig}
\usepackage{color}
\usepackage{enumerate}
\usepackage{epstopdf}
\usepackage{titlesec}
\usepackage[titletoc,toc,title]{appendix}
\usepackage{subfigure}
\usepackage[utf8]{inputenc}
        \usepackage[T1]{fontenc}
        \usepackage{lmodern}
        \usepackage{graphicx}
        \usepackage{caption}
         \usepackage{floatrow}%
\usepackage{tikz}
\usetikzlibrary{arrows,positioning, arrows}

\usepackage{amsthm}
\usepackage{amsfonts}
\usepackage{graphicx}
\usepackage{caption}
\usepackage{floatrow}
\usepackage{amsmath}
\usepackage{amssymb}
\usepackage{amsmath}
\usepackage{bm}
\usepackage{subeqnarray}
\usepackage[linesnumbered,ruled]{algorithm2e}
\usepackage[noend]{algpseudocode}
\usepackage{cases}
\theoremstyle{plain}
    \newtheorem{theorem}{Theorem}
    \newtheorem*{Thm*}{Main Theorem}
    \newtheorem{lemma}{Lemma}
    
    \newtheorem{proposition}{Proposition}
    \newtheorem{corollary}{Corollary}

\theoremstyle{definition} 
    
    \newtheorem{remark}{Remark}

\textheight 9in \textwidth 6.5in \numberwithin{equation}{section}

\newcommand{\eps}{\varepsilon}
\newcommand{\norm}[1]{ \left\| #1 \right\|}
\newcommand{\normM}[1]{ \| #1 \|_{\mathcal{M}^{-1}}}

\newcommand{\normMinf}[1]{ \| #1 \|_{L^\infty_{{\mathcal{M}^{-1}}}}}
\newcommand{\normv}[1]{ \| #1 \|_{L_v^2}}
\newcommand{\normxv}[1]{ \| #1 \|_{L_{x,v}^2}}

\newcommand{\half}{\frac{1}{2}}
\newcommand{\halfd}{\frac{d}{2}}
\newcommand{\RR}{\mathbb{R}}
\newcommand{\rd}{\mathrm{d}}
\newcommand{\equilibrium}{\mathcal{M}}
\newcommand{\Lop}{\mathcal{L}}
\newcommand{\Hop}{\mathcal{H}}

\newcommand{\vpran}[1]{\left(#1\right)}
\newcommand{\sgn}{\textrm{sgn}^+}

\newcommand{\nn}{\nonumber}
\newcommand{\abs}[1]{\left\lvert#1\right\rvert}
\newcommand{\Ni}{\noindent}
\newcommand{\CalM}{\mathcal{M}}
\newcommand{\CalL}{\mathcal{L}}
\newcommand{\dv}{ \, {\rm d}v}

\newcommand{\ftilde}{\tilde{f}}

\newcommand{\egn}{e_1^n}
\newcommand{\een}{e_2^n}
\newcommand{\epsdt}{\frac{\varepsilon^{2s}}{\Delta t}}
\newcommand{\tnp}{t^{n+1}}

\newcommand{\dt}{\Delta t}

\newcommand{\tgn}{\tilde g^n}
\newcommand{\ten}{\tilde \eta^n}
\newcommand{\tgnp}{\tilde g^{n+1}}
\newcommand{\tenp}{\tilde \eta^{n+1}}
\newcommand{\gstar}{g^{n+\frac{1}{2}}}
\newcommand{\np}{{n+1}}
\newcommand{\nh}{{n+\frac{1}{2}}}
\newcommand{\expo}{e^{-\eps^{-2s}t}}
\newcommand{\epst}{\eps^{-2s}t}
\newcommand{\finv}{\mathcal{F}^{-1}}
\newcommand{\id}{\mathcal{I}}

\newcommand{\average}[1]{\left\langle#1\right\rangle}

\graphicspath{{../}}

\title{Uniform error estimate of an asymptotic preserving scheme for the L\'{e}vy-Fokker-Planck equation
\footnote{W. S. is partially supported by NSERC Discovery Grant R611626. L.W. is partially supported by NSF CAREER grant DMS-1846854. We would like to thank the Isaac Newton Institute for Mathematical Sciences, Cambridge, for support and hospitality during the programme, {\it Frontiers in kinetic theory: connecting microscopic to macroscopic scales}, where part of the work on this paper was undertaken.}}
\date{}
\author{Weiran Sun \thanks{Department of Mathematics, Simon Fraser University, Burnaby BC V5A 4X9, CA. (weirans@sfu.edu)} 
and Li Wang \thanks{School of Mathematics, University of Minnesota, Twin cities, MN 55455. (liwang@umn.edu) }
}

\begin{document}
\maketitle
\begin{abstract}
We establish a uniform-in-scaling error estimate for the asymptotic preserving scheme proposed in \cite{XW21} for the L\'evy-Fokker-Planck (LFP) equation. The main difficulties stem not only from the interplay between the scaling and numerical parameters but also the slow decay of the tail of the equilibrium state. 
We tackle these problems by separating the parameter domain according to the relative size of the scaling parameter $\eps$: in the regime where $\eps$ is large, we design a weighted norm to mitigate the issue caused by the fat tail, while in the regime where $\eps$ is small, we prove a strong convergence of LFP towards its fractional diffusion limit with an explicit convergence rate. This method extends the traditional AP estimates to cases where uniform bounds are unavailable. 
Our result applies to any dimension and to the whole span of the fractional power. 
\end{abstract}

\section{Introduction}
Consider the L\'{e}vy-Fokker-Planck (LFP) equation 
\begin{equation} \label{eqn:LFP}
\begin{cases}{}
\partial_t f + v \cdot \nabla_x f = \nabla_v \cdot (v f) - (-\Delta_v)^s f :=\Lop^s (f), \qquad s\in (0,1),
\\ f(0,x,v) = f_{in}(x,v),
\end{cases}
\end{equation}
where $f(t,x,v): (0, \infty) \times \RR^d \times \RR^d  \mapsto  \RR^+$ is the distribution function of a large group of particles, which undergoes a free transport dynamics along with an  interaction with the background, described by the L\'evy-Fokker-Planck operator.  In contrast to the conventional Fokker-Planck operator, here we have $(-\Delta_v)^s$  in place of $\Delta_v$, which models the L\'evy processes at the microscopic level instead of Brownian motion. One way to understand  $(-\Delta_v)^s$ is via the Fourier transform. Namely, for any $\phi(v) \in L^1(\RR^d)$, 
\begin{align*}
-(-\Delta_v)^s \phi := \mathcal{F}^{-1} (|k|^{2s} \mathcal{F} \phi)\,,
\end{align*}
where 
\begin{align*}
{\hat \phi}(k):=\mathcal F[ \phi](k) = \int_{\mathbb R^d} \phi(v) e^{-i  v \cdot k} \rd v, 
\qquad 
\mathcal F^{-1}[ \phi](v) 
= \frac{1}{(2\pi)^d} \int_{\mathbb R^d} \hat \phi(k) e^{i  v \cdot k} \rd k 
\end{align*}
are the Fourier and inverse Fourier transforms. Formally, from this definition $- (-\Delta_v)^s$ reduces to $\Delta_v$ if $s = 1$. Another way of defining this operator is through the integral representation: 
\begin{equation}
    (-\Delta_v)^s f(v) : = C_{s,d}~ \text{P.V.} \int_{\mathbb{R}^d} \frac{f(v)-f(w)}{|v-w|^{d + 2s}} \rd w,
\label{def_fl}
\end{equation}
where $P.V.$ denotes the Cauchy principal value and $C_{s,d} = \frac{4^s \Gamma(d/2+s)}{\pi^{d/2}|\Gamma(-s)|}$. In principle, the fractional Laplacian allows particles to make long jumps at the microscopic scale, which leads to the nonlocal effect as written in \eqref{def_fl} at the mesoscopic scale. 

In the small mean free path regime with a long time, equation \eqref{eqn:LFP} can be rescaled as
\begin{equation} \label{eqn:111}
\begin{cases}{}
\eps^{2s}\partial_t f^\eps + \eps v \cdot \nabla_x f^\eps = \mathcal L^s (f^\eps) ,
\\ f^\eps(0,x,v) = f_{in}(x,v).
\end{cases}
\end{equation}
Formally, $f^\eps$ converges to $\rho(t,x) \equilibrium (v)$ as $\eps \to 0$, where $\equilibrium(v)$ is the unique normalized equilibrium of~\eqref{eqn:111} (see \cite{gentil2008levy}) with the properties
\begin{equation} \label{eqn:equilibrium}
\Lop^s(\equilibrium) = 0, \quad \int_{\RR^d} \equilibrium (v) \rd v = 1, \quad \equilibrium(v) \sim \frac{C}{|v|^{d+2s}} ~\text{as}~ |v| \rightarrow \infty.
\end{equation}
Meanwhile, the limiting density $\rho(t,x)$ solves 
\begin{equation}\label{eqn: limit_system}
\begin{cases}{}
\partial_t \rho + (-\Delta_x)^s \rho= 0,
\\ \rho(0,x) = \rho_{in}(x) := \int_{\RR^d} f_{in}(x,v)\rd v.
\end{cases}
\end{equation}
A weak convergence from $f^\eps$ to $\rho \equilibrium$ has been established in  \cite{cesbron2012anomalous}, and extended to the case with electrical field \cite{aceves2019fractional}. In Section~\ref{sec:APlimit}, we will strengthen this result in two aspects: first, we will show that, this convergence is strong in $L^2_{x, v}$. Second, we will provide an explicit convergence rate in terms of $\eps$. See Theorem~\ref{APlimit} for details. 

A notable difference between the anomalous scaling in \eqref{eqn:111} and classical diffusive scaling (i.e., $s=1$ in \eqref{eqn:111}) is that the classical diffusion limit takes the form 
\begin{align*}
\partial_t \rho + \nabla_x \cdot ( D \nabla_x \rho)  = 0, 
\end{align*}
where $D$ is the diffusion matrix 
\begin{align*}
D = \int v \otimes  v \equilibrium \rd v.  
\end{align*}
It is clear that  the fat tail equilibrium \eqref{eqn:equilibrium} renders $D$ unbounded and therefore necessitates the anomalous scaling. Similar scaling has also been investigated in the framework of linear Boltzmann equation
\begin{align}\label{Bol}
\eps^{2s}\partial_t f + \eps v \cdot \nabla_x f = \int \phi(v,v') (\equilibrium(v) f(v') - \equilibrium(v') f(v)) \rd v'\,,
\end{align}
where $\equilibrium(v)$ is a given heavy tailed equilibrium satisfies \eqref{eqn:equilibrium}, and the diffusion limit is again \eqref{eqn: limit_system} (see \cite{abdallah2011fractional, mellet2011fractional}).

From a computational perspective, it is desirable to design a method that performs uniformly in $\eps$. An asymptotic preserving (AP) scheme would suffice this purpose \cite{jin2022asymptotic}. However, due to the algebraic decay in the tail of the equilibrium, existing methods that work for exponential decay equilibrium (i.e., Maxwellian) cannot be applied. The reason is that these methods always reside in a truncated velocity domain as the tail information is negligible. By contrast, when the anomalous scaling is considered, the tail carries the most important information that cannot be ignored. For this reason, special treatment is needed to preserve the tail information along the dynamics. This is the main idea behind the works \cite{crouseilles2016numerical, CHL16, wang2016asymptotic, wang2019asymptotic}, all of which have dealt with the linear Boltzmann type equation in \eqref{Bol}. Even though the LFP equation shares a similar equilibrium state and diffusion limit with the linear Boltzmann equation~\eqref{Bol}, its structural differences from~\eqref{Bol} prohibit a direct application of methods developed for~\eqref{Bol}. There are two main obstacles. The first one is due to the differential operators involved in LFP. More specifically, the method used for~\eqref{Bol} that directly compensates the heavy tail using mass conservation does not apply to LFP since it will violate the smoothness requirement for LFP's differential operators. Second, the Hilbert expansion, which is key to the design of many AP schemes \cite{jin2022asymptotic, wang2016asymptotic, wang2019asymptotic}, applies well to the linear Boltzmann equation. Unfortunately we do not know yet how to adapt it for LFP.

In \cite{XW21}, a new asymptotic preserving scheme was developed to alleviate the aforementioned difficulties. The main idea is based on a novel type of micro-macro decomposition, with a unique macro part that is inspired by the special choice of the test function in proving the weak convergence in \cite{cesbron2012anomalous}. 
%
More specifically, decompose $f$ as
\begin{equation}\label{eqn: decomp}
    f(t,x,v) = \eta(t,x, v) \equilibrium(v) +g(t,x,v),
\end{equation}
where the regularity of $\eta$ and $g$ are specified in Remark~\ref{rmk:1}.
Moreover, $\eta(t,x,v) $ takes the form
\begin{equation} \label{eqn:tr}
    \eta(t,x,v) = h(t,x+\eps v)
\end{equation} 
for some function $h(t,x)$ and $\equilibrium$ is the equilibrium state satisfying
\begin{equation} \label{equilibrium000}
\partial_v (v \equilibrium) - (-\Delta_v)^s \equilibrium =0, \qquad \int_{\mathbb R} \equilibrium(v) \rd v = 1. 
\end{equation}
Note that although $\eta$ depends on $(x, v)$, intrinsically it lives on a lower dimensional manifold than $f$ does. Direct computation using \eqref{eqn:tr} gives
\begin{equation} \label{eqn032}
\eps \partial_x \eta =  \partial_v \eta, \qquad  (-\Delta_v)^s \eta =\eps^{2s} (-\Delta_x)^s \eta .  
\end{equation}
Inserting $\eqref{eqn: decomp}$ into $\eqref{eqn:111}$, we get
\begin{equation*} 
   \eps^{2s}\partial_t (\eta \equilibrium +g) 
   + \eps v  \partial_x (\eta \equilibrium +g) 
= \partial_v  (v (\eta \equilibrium +g)) 
    - (-\Delta_v)^s (\eta \equilibrium +g).
\end{equation*}
By \eqref{equilibrium000} and \eqref{eqn032}, it simplifies to
\begin{equation}\label{eqn: decomp_sym}
    \eps^{2s}\partial_t (\eta \equilibrium +g) + \eps v  \partial_x g = \Lop^{s} (g)  - \eps^{2s} (-\Delta_x)^s \eta \equilibrium - I(\eta,\equilibrium) ,
\end{equation}
where 
\begin{align} \label{I}
I(p, q)  & = (-\Delta_v)^s (p \, q ) - q \, (-\Delta_v)^s p  - p \, (-\Delta_v)^s q   \nonumber
\\ & =C_{d,s} \int_{\RR^d} \frac{(p(v)-p(w))(q(w)-q(v))}{|v-w|^{d+2s}} \rd w.
\end{align}
Splitting \eqref{eqn: decomp_sym}, we introduce the following macro-micro model 
\begin{equation} \label{split0}
    \begin{cases}{}
     \partial_t \eta = -(-\Delta_x)^s \eta,
 \\  \eps^{2s}\partial_t g + \eps v \partial_x g =\Lop^s(g) - I(\eta, \equilibrium)  ,
    \end{cases}
\end{equation}
with the initial condition given as
\begin{equation} \label{IC}
  \eta_{in}(x,v) = \rho_{in}(x+\eps v), 
\qquad 
  g_{in}(x,v) = f_{in}(x,v) - \eta_{in}(x,v) \equilibrium (v).
\end{equation}
One can recover $f$ in \eqref{eqn:111} by solving \eqref{split0} and using~\eqref{eqn: decomp}.

To numerically solving \eqref{split0}, we propose a semi-discrete scheme based on an operator splitting: 
\begin{subequations}\label{semi_scheme}
\begin{numcases}{}
 \frac{1}{\Delta t} (\eta^{n+1}-\eta^n)  =  - (-\Delta_x)^s \eta^{n+1} , \label{semis3}\\
\frac{\eps^{2s}}{\Delta t} (g^{\nh}-g^n)  = \Lop^s(g^\nh) - \gamma g^\nh - I(\eta^{n+1},\equilibrium)  , \label{semis1} \\
 \frac{\eps^{2s}}{\Delta t}(g^{n+1}-g^\nh)  + \eps v  \partial_x g^{n+1} = \gamma  g^{n+1}  , \label{semis2} 
\end{numcases}
\end{subequations}
where $\gamma$ is a positive constant chosen according to 
\begin{align*}
& \text{When}~ \eps^{2s} \geq \Delta t^{2s\beta}:  ~~ 2 < \gamma \leq \min\{\sqrt{\frac{\lambda_2-4}{\lambda_2}} \frac{\eps^{2s}}{\Delta t} \,,~ 4\}  \quad \text{for some constant} ~ \lambda_2 >4\,;
\\ & \text{When}~\eps^{2s} <  \Delta t^{2s\beta}: \gamma = \left \{ \begin{array}{cc} \sqrt{3} & \text{if}~ \eps^{2s} < \Delta t \\ \sqrt{3} \Delta t^{2s\beta -1} & \text{if}~ \Delta t \leq \eps^{2s} \leq \Delta t^{2s\beta}\,,
\end{array} \right.
\end{align*}
where $0<\beta< 1/{6s}$. See \eqref{cond:gamma1}, \eqref{cond:gamma2} and \eqref{gamma-II} for more details. The spatial derivative will be treated via the Fourier-based spectral method, and the velocity discretization will be done by using the mapped Chebyshev-polynomial-based pseudo-spectral method. Details on the velocity discretization can be found in \cite[Section 2]{XW21}. The scheme described in \eqref{semi_scheme} is slightly different from the one presented in \cite{XW21}, where in \eqref{semis3} we define $\eta^{n+1}$ implicitly and use $\eta^{n+1}$ as a source to derive $g^\nh$ in \eqref{semis1}. This adjustment will not result in any additional computational overhead, given that we employ a Fourier-based spectral method for the variable $x$. The well-posedness of~\eqref{semi_scheme} is stated in Lemma~\ref{lem:well-posedness-scheme}.

\begin{figure}[!h]
\includegraphics[width=0.36\textwidth]{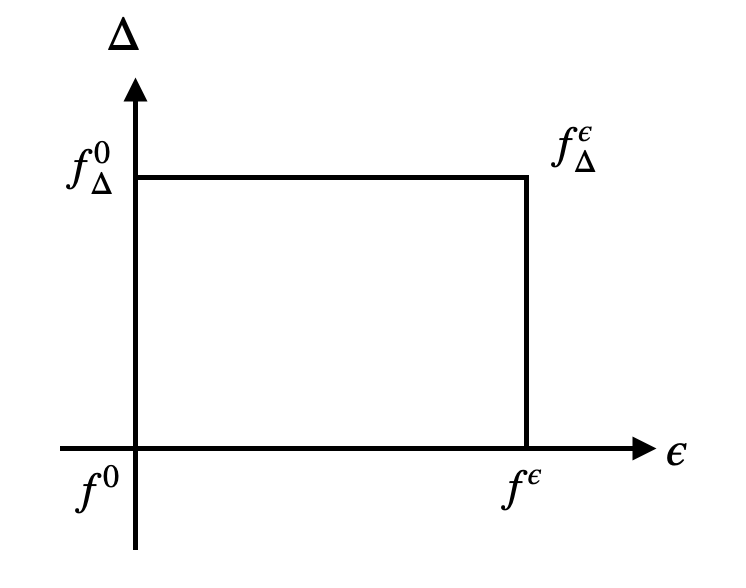}
\caption{\small Illustration of AP schemes.}
\label{fig: AP}
\end{figure}

The goal of this paper is to provide a rigorous foundation to the above scheme. In particular, we will show that the scheme is asymptotic preserving (AP). Hence the convergence of the scheme is uniform with respect to both $\eps$ and $\Delta t$. The general framework of AP schemes was laid out in \cite{golse1999convergence} and the main idea is illustrated in Fig.~\ref{fig: AP}. 
Denote by $f^\eps$ and $f^0$ the exact solutions to the kinetic equation (e.g., \eqref{eqn:LFP}) and macroscopic limit (e.g, \eqref{eqn: limit_system}), respectively. The same notation with subscript $\Delta $ represents the associated numerical solution. The key idea is to establish a uniform bound for $\norm{f^\eps_\Delta  -f^\eps}$ by optimizing over two routes: one is the direct error estimate in the kinetic regimes which gives
\begin{equation} \label{0824}
\norm{f^\eps_\Delta  -f^\eps} \leq C \frac{\Delta^p}{\eps},
\end{equation}
where $p$ is the order of the adopted numerical discretization. The other is through the asymptotic relation
\begin{equation} \label{0823}
\norm{f^\eps_\Delta  -f^\eps}  \leq \norm{f^\eps_\Delta  -f^0_\Delta} + \norm{f^0_\Delta - f^0} + \norm{f^0 - f^\eps},
\end{equation}
where $\norm{f^0_\Delta - f^0} \leq C \Delta^p$ as there is no $\eps$ dependence in the equation, while $\norm{f^0 - f^\eps} \leq C \eps$ is guaranteed by the asymptotic relation at the continuum level. {\it If} the scheme is asymptotic  preserving, then we have 
\begin{equation} \label{0826}
\norm{f^\eps_\Delta  -f^0_\Delta} \leq C \eps. 
\end{equation}
Consequently, \eqref{0823} becomes
\begin{equation} \label{0825}
\norm{f^\eps_\Delta  -f^\eps} \leq C( \Delta^p + \eps)\,.
\end{equation}
Finally, an optimization between \eqref{0824} and \eqref{0825} gives us  the uniform estimate 
\begin{equation} \label{0827}
\norm{f^\eps_\Delta  -f^\eps} \leq C \Delta^{p/2}\,.
\end{equation}
We emphasize that throughout the procedures above, $C$ is a constant {\it independent} of both $\eps$ and $\Delta$, a property that we view as a certain uniformity of the traditional AP estimates. 

This framework, albeit universal, is not easy to carry out for specific problems. Oftentimes the bottleneck is to obtain the estimate \eqref{0826}, which can be as difficult as obtaining the uniform error estimate directly. For this reason, despite the proliferation of AP schemes, there is a notable scarcity in rigorously justifying their uniform accuracy. To the best of our knowledge, only a limited number of works \cite{golse1999convergence, gosse2004asymptotic, hu2021uniform, jin2009uniformly, klar2002uniform, liu2010analysis, filbet2013analysis, jang2014analysis, li2017uniform, bessemoulin2020hypocoercivity, filbet2021convergence, peng2021stability} provide a rigorous stability or error estimate. Notably, except for \cite{filbet2021convergence}, which focuses on plasma models, these works are predominantly confined to linear transport or relaxation-type equations.

\begin{figure}[!h]
\includegraphics[width=0.35\textwidth]{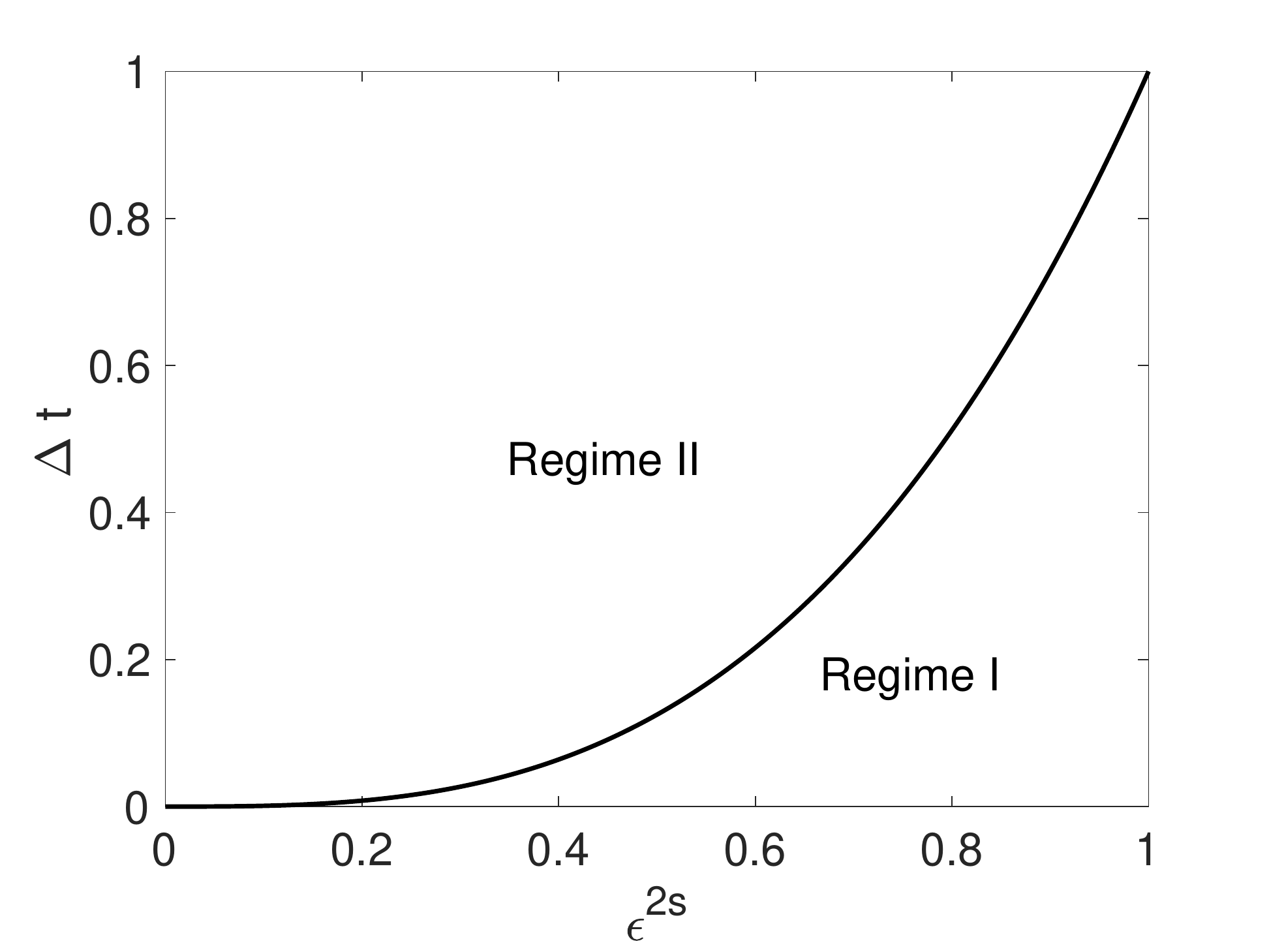}
\includegraphics[width=0.35\textwidth]{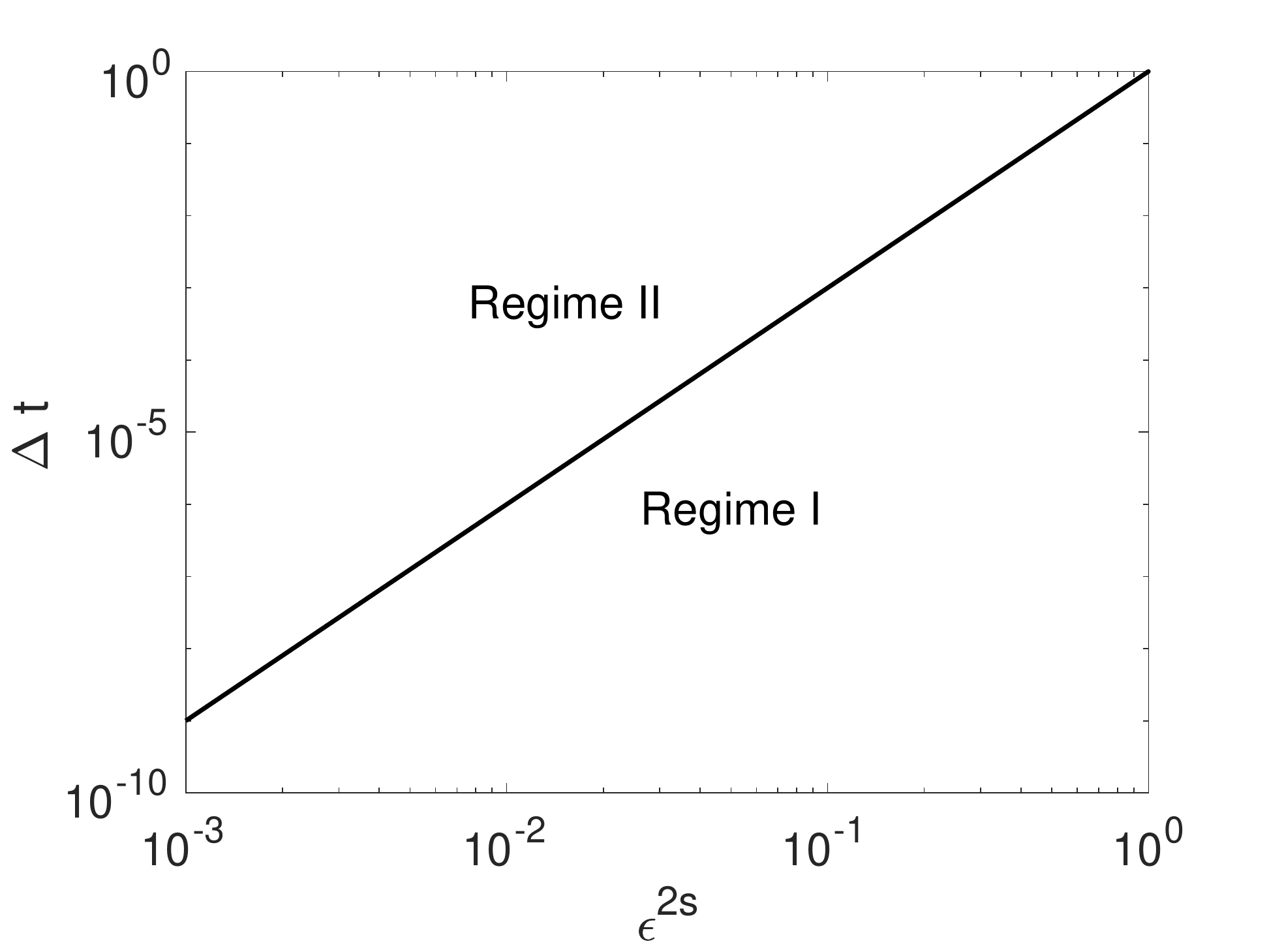}
\caption{\small Two regimes separated by the relation between $\eps^{2s}$ and $\Delta t $. Left: linear scale. Right: log scale.}
\label{fig: regimes}
\end{figure}

In this paper, we prove a uniform error estimate for the first time for the {\it fractional} kinetic equation. In our case, the strong non-locality, both in the L\'evy-Fokker-Planck operator itself and in the slow algebraic decay of the equilibrium,  poses tremendous difficulties in the proof as the basic energy estimate fails immediately. To solve this problem, we propose a domain decomposition as depicted in Fig.~\ref{fig: regimes}, which can be viewed as a generalization of the classical AP framework described above. In particular, we divide the domain into two regimes---Regime I with $\eps^{2s} \geq \Delta t^{2s\beta}$ and Regime II with $\eps^{2s} < \Delta t^{2s\beta}$, where the constant $\beta \in (0, 1/6s)$ is specified in Theorem~\ref{thm:2}. 

Denote $W^{k, p}$ and $H^k$ as the usual Sobolev spaces. Introduce the notation $L^2_{\equilibrium^{-1}}$ such that 
\begin{align}
&f \in L^2_{\equilibrium^{-1}}
\Longleftrightarrow
\normM f^2 = \int \int \frac{f^2}{\equilibrium} \rd x \rd v < \infty, \label{normM}
\\
&f \in L^\infty_{\equilibrium^{-1}} 
\Longleftrightarrow
 \norm{\frac{f}{\CalM}}_{L^\infty_{x, v}} < \infty. \label{normMinfty}
\end{align}
Leaving out some details, our main result is to establish the following error estimate 
\begin{Thm*}
Suppose the initial data are smooth enough such that 
\begin{align} 
v f_{in} \in W^{1,1}_{x, v},  
\quad 
  |v|^2 D^2_x f _{in}\in L^2_{\CalM^{-1}}, 
\quad 
  |v|^2 D^2_v f _{in}\in L^2_{\CalM^{-1}}, \label{cond:reg-initial-1}
\\
  \rho_{in} \in H^2_x \cap W_x^{4s, \infty}\cap W_x^{1, \infty}, 
\quad
  g_{in}, \, \nabla_x g_{in}, \, D^2_v g_{in} \in L^2_{\equilibrium^{-1}},  \label{cond:reg-initial-2}
\end{align}
where $D^2_x, D^2_v$ denote the collection of second-order derivatives in $x, v$ respectively.
Then outside of an initial layer, there exist constants $C, b, \zeta > 0$ independent of $\eps$ and $\Delta t$ such that 
\begin{align*}
    \normM{\average{v}^{-b} (f^\eps_\Delta  - f^\eps)} 
\leq  
   C \Delta t^{\zeta}. 
\end{align*}
The parameters $b, \zeta$ are specified in Theorem~\ref{thm:2} and Theorem~\ref{thm:3}. 
\end{Thm*}

\begin{remark} \label{rmk:1}
The well-posedness of~\eqref{eqn:111}, thus~\eqref{split0}, is given in \cite{cesbron2012anomalous}. Given the initial data in~\eqref{cond:reg-initial-1}-\eqref{cond:reg-initial-2}, the PDE solution in~\eqref{split0} satisfies at least the same regularity as the initial data for $\rho, g$ respectively.
\end{remark}

The main theorem is proved in different ways for Regimes I and II: 
\begin{itemize}
\item in Regime I when $\eps$ is relatively large compared with $\Delta t$ and the kinetic equation may not be well-approximated by the macroscopic equation, we bound the error by performing a weighted energy estimate for the semi-discrete kinetic error equation \eqref{erroreqn2} directly. In particular, Lemmas~\ref{lemma2}-\ref{lemma6} are to bound the source term generated by the operator splitting. In the energy estimate, a weight $(1 + \delta \average{v})^{-b}$ is introduced to compensate the slow decay tail of the solution. Since this weight is not commutable with the equation, Lemmas~\ref{lemt1}-\ref{lemma:commutator} provide bounds for the commutator terms; 
\item In Regime II where $\eps$ is small compared with $\Delta t$, the analysis in done in a more straightforward way by comparing the kinetic equation with the macroscopic equation, where the continuous asymptotic error $g(t^n, x, v)$ and the discrete error $g^n(x, v)$ are explicitly bounded in Theorem~\ref{APlimit} and Lemma~\ref{lemma4} respectively.
\end{itemize}

Compared with the general framework outlined in \eqref{0824}--\eqref{0827}, the main novelty in our proof lies in the relaxation of the two key inequalities \eqref{0824} and \eqref{0826}. In fact, owing to the slow decay of the equilibrium $\mathcal{M}$ which may not even have a finite first moment, we cannot obtain estimates like \eqref{0824} or \eqref{0826} with a uniform constant $C$. Instead, we relax both upper bounds to a form of $\Delta t^a \eps^b$, wherein one of $a$ and $b$ can be negative while  the product as a whole is a quantity with positive power in $\Delta t$. This underpins the relation between $a$ and $b$ and leads to the division of the parameter domain. We expect that this generalization provides a new route for proving the uniform accuracy for asymptotic preserving schemes beyond our current work. 

The rest of the paper is organized as follows. In Section ~\ref{sec:APlimit}, we prove the strong asymptotic limit from the kinetic to the macroscopic systems in the $L^2_{x,v}$-norm with an explicit convergence rate. This result will be used in obtaining the error estimate in Regime II.  
Section 3 consists of technical lemmas on commutator estimates, which will appear in the error estimate in Regime I. Finally, Section 4 is devoted to the proof of the main theorem.

In this paper, we follow the convention of allowing constants $C, C_{s, t_0}$ to change from line to line.

\section{Strong asymptotic limit} \label{sec:APlimit}
This section is devoted to establishing a strong convergence from \eqref{eqn:111} to \eqref{eqn: limit_system} with an explicit convergence rate. Our main tool is the Fourier transform, through which solutions are derived explicitly. In particular, denote $\hat f(t,\xi, k)$ as the Fourier transform of $f(t,x,v)$, i.e., 
\begin{align*}
\hat f(t,\xi,k) = \int_{\mathbb R^d} \int_{\mathbb R^d} f(t,x,v) e^{-i (x \cdot \xi + v \cdot k)} \rd x \rd v.
\end{align*}
Then \eqref{eqn:111} rewrites as 
\begin{equation*}
\begin{cases}{}
\eps^{2s} \partial_t \hat f - \eps \xi \cdot \nabla_k \hat f = - k \cdot \nabla_k \hat f - |k|^{2s} \hat f,
\\ \hat f(0,\xi, k) = \hat f_{in}(\xi, k).
\end{cases}
\end{equation*}
By the method of characteristics, its solution is 
\begin{equation} \label{fhat}
\hat f(t,\xi, k) 
= \hat f_{in} (\xi, \eps \xi + \expo (k - \eps \xi))  \,
   e^{-\int_0^{\epst} |e^{-w}(k-\eps \xi) + \eps \xi|^{2s} \rd w}.
\end{equation}
Similarly, if we denote $\hat \rho(t,\xi)$ as the Fourier transform of $\rho(t,x)$, then from \eqref{eqn: limit_system} we have 
\begin{equation*} 
\hat \rho(t,\xi) = \hat f_{in} (\xi,0) e^{-|\xi|^{2s} t }.
\end{equation*}
Thus the equilibrium state in the Fourier space is
\begin{equation} \label{rhoMhat}
\hat \rho(t,\xi) \hat \equilibrium(k) = \hat f_{in} (\xi,0) e^{-|\xi|^{2s} t -\frac{1}{2s}|k|^{2s}}  .
\end{equation}
The difference of the solutions in \eqref{fhat} and \eqref{rhoMhat} is as follows. 
\begin{theorem} \label{APlimit}
Let $f$ and $\rho$ be the solutions to \eqref{eqn:111} and \eqref{eqn: limit_system} respectively. Let $t_0 > 0$ be a fixed constant. Then for any $t \geq t_0$, we have
\begin{align*}
& \norm{f - \rho \equilibrium}_{L_{x,v}^2}^2 \leq 
C_{s,t_0}  \eps^{2s} \norm{\rho_{in}}_{L^1}^2
+
C \left[ \eps^2 \norm{v \cdot \nabla_x f_{in}}_{L^1_{x,v}}^2  + C_{t_0} \eps^{4s} \norm{  v f_{in}}_{L^1_{x,v}}^2 \right],
\end{align*}
where $C_{s,t_0}$, $C$ and $C_{t_0}$ are three constants that do not depend on $\eps$. 
\end{theorem}

Some preparations for proving Theorem~\ref{APlimit} are in place. Frist we prove an elementary lemma that will be used repeatedly. 
\begin{lemma} For all  $s \in (0,1)$ and $a, b >0$, we have the following inequalities: 
\begin{align}
& (a+b)^s \leq a^s + b^s; \label{0221}
\\ & |(a+b)^{2s} - (a^{2s} + b^{2s}) |\leq 2 a^s b^s. \label{0525}
\end{align}
\end{lemma}
\begin{proof}
The first inequality follows from the concavity of the function $x \mapsto x^s$ with $x \geq 0$:
\begin{align*}
   \vpran{\frac{a}{a+b}}^s + \vpran{\frac{b}{a+b}}^s 
\geq 
  \vpran{\frac{a}{a+b} + \frac{b}{a+b}}^s = 1.
\end{align*}
To prove~\eqref{0525}, by symmetry we can assume that $a \geq b$. Denote $z = \frac{b}{a} \in [0, 1]$. Then~\eqref{0525} is equivalent to
\begin{align*}
   \abs{(1+z)^{2s} - (1 + z^{2s})} \leq 2 z^s,
\end{align*}
or equivalently, 
\begin{align} \label{ineq:equiv-2}
   (1+z)^{2s} - (1 + z^{2s}) \leq 2 z^s
\quad \text{and} \quad
   (1 + z^{2s}) - (1+z)^{2s} \leq 2 z^s. 
\end{align}
The first inequality in~\eqref{ineq:equiv-2} holds since by~\eqref{0221},
\begin{align*}
    (1+z)^{2s} \leq (1+z^s)^2 = (1 + z^{2s}) + 2 z^s. 
\end{align*}
The second inequality in~\eqref{ineq:equiv-2} holds since
\begin{equation*}
   1 + z^{2s} - 2 z^s
= (1 - z^s)^2
\leq (1+z)^{2s}, 
\qquad z \in [0, 1].  \qedhere
\end{equation*}   
\end{proof}

By Parseval's identity, we will prove Theorem~\ref{APlimit} in the Fourier space. Subtracting \eqref{rhoMhat} from \eqref{fhat}, we have
\begin{align}
& \hat f(t,\xi, k) - \hat\rho(t,\xi) \hat \equilibrium (k). \nonumber
\\ = &  ~\hat f_{in} (\xi, 0) [  e^{-\int_0^{\epst} |e^{-w}(k-\eps \xi) + \eps \xi|^{2s} \rd w} - e^{-|\xi|^{2s} t - \frac{1}{2s}|k|^{2s}} ]  \nonumber
\\  & \qquad + [\hat f_{in} (\xi, \eps \xi + \expo (k - \eps \xi))- \hat f_{in} (\xi, 0)] e^{-\int_0^{\epst} |e^{-w}(k-\eps \xi) + \eps \xi|^{2s} \rd w}    \nonumber  
\\ =: & ~\hat f_{in} (\xi, 0) A_1 + A_2.  \label{diff0} 
\end{align}
The following lemma quantifies the difference in the exponents in $A_1$.
\begin{lemma} \label{lem1}
Suppose $t_0 > 0$ is a fixed constant. Then for any fixed $k$ and $\xi$, there exists a constant $C_{s,t_0}$ depending only on $s$ and $t_0$, such that 
\begin{align*}
 \abs{\int_0^{\eps^{-2s}t} |e^{-w}(k-\eps \xi) + \eps \xi|^{2s} \rd w   
         - \vpran{|\xi|^{2s} t +  \frac{|k|^{2s}}{2s}}}
\leq C_{s,t_0} \eps^{s}  \vpran{|\xi|^{2s} + |k|^{2s}} , \quad \forall~ t  \geq t_0.
\end{align*}
\end{lemma}
\begin{proof}
First we reformulate the integral by decomposing $\xi$ into the component in the direction of $k$ and the one perpendicular to it: $\xi_k := \xi \cdot \frac{k}{|k|}$ and $\xi^\perp = \xi - \xi_k$. Then
\begin{align*}
 J:=& \int_0^{\eps^{-2s}t} |e^{-w}(k-\eps \xi) + \eps \xi|^{2s} \rd w
 \\ 
 = &  \int_0^{\eps^{-2s}t} \left[ |e^{-w} |k| + (1-e^{-w}) \eps \xi_k |^2 +  (1-e^{-w})^2 \eps^2 |\xi^\perp|^2 \right]^s \rd w 
 \\ 
 = & \int_0^{\eps^{-2s}t} \big[ \underbrace{e^{-2w}|k|^2 + 2e^{-w} (1-e^{-w}) \eps \, \xi_k \, |k|}_a + \underbrace{(1-e^{-w})^2 \eps^2|\xi|^2}_{b} \big]^s  \rd w 
 \\ = & \int_0^{\eps^{-2s}t} \vpran{(a+b)^s - a^s - b^s + a^s + b^s}  \rd w. 
\end{align*} 
Hence,
\begin{align*}
  \abs{J - \vpran{|\xi|^{2s} t +  \frac{|k|^{2s}}{2s}}}
&\leq
  \int_0^{\eps^{-2s}t} \abs{(a+b)^s - a^s - b^s} \rd w
\\
& \quad \, 
  +  \abs{\int_0^{\eps^{-2s}t} a^s \rd w - \frac{|k|^{2s}}{2s}}
  +  \abs{\int_0^{\eps^{-2s}t} b^s \rd w - |\xi|^{2s}t}. 
\end{align*}
We show that each term is of order $\eps^s$. By~\eqref{0525}, 
\begin{align} \label{0621}
   \int_0^{\eps^{-2s}t} \abs{(a+b)^s - a^s - b^s} \rd w
& \leq
   2 \int_0^{\eps^{-2s}t} a^{\frac{s}{2}} b^{\frac{s}{2}} \rd w  \nn
\\
& \hspace{-1.2cm}
= 2\int_0^{\eps^{-2s}t} [(e^{-w}|k|)^2 + 2e^{-w} (1-e^{-w}) \eps \xi_k |k|]^{\frac{s}{2}} 
           (1-e^{-w})^s\eps^s|\xi|^s \rd w \nonumber
\\ 
& \hspace{-1.2cm}
\leq 
  2\int_0^{\eps^{-2s}t} [(e^{-w}|k|)^s + (2e^{-w} (1-e^{-w}) \eps \xi_k |k|)^{\frac{s}{2}}] 
           (1-e^{-w})^s\eps^s|\xi|^s \rd w \nonumber
\\ 
& \hspace{-1.2cm}
\leq 
  C \eps^s[|k|^{2s} + |\xi|^{2s}],
\end{align}
where the first inequality follows from \eqref{0221} and $C$ is a generic constant. Next,
\begin{align} \label{I3-may}
 \left|  \int_0^{\eps^{-2s}t}  b^s \rd w -  |\xi|^{2s}t  \right| 
 & = \eps^{2s} |\xi|^{2s} \int_0^{\eps^{-2s}t} \abs{(1-e^{-w})^{2s}  -1} \rd w \leq 2 \eps^{2s} |\xi|^{2s} 
\end{align}
where we have used the fact that for $s \in (0,1)$,
\begin{align*}
 1 - (1-e^{-w})^{2s}  
&= [1 - (1-e^{-w})^{s}][1 + (1-e^{-w})^{s}] 
\\
&\leq 2 [1 - (1-e^{-w})^{s}] 
\leq 2 [1 - (1-e^{-w})] 
= 2e^{-w}\,.
\end{align*}
Using~\eqref{0221} again, we have
\begin{align*}
 \int_0^{\eps^{-2s}t} (e^{-w}|k|)^{2s} \rd w
 \leq 
   \int_0^{\eps^{-2s}t} a^s \rd w 
& \leq  
  \int_0^{\eps^{-2s}t} (e^{-w}|k|)^{2s} + (2e^{-w}(1-e^{-w}) \eps \xi_k |k|)^s \rd w \nonumber
 \\ 
 & \leq  \int_0^{\eps^{-2s}t} (e^{-w}|k|)^{2s} \rd w 
 +  C \eps^s \vpran{|\xi_k|^{2s} + |k|^{2s}} \,,
\end{align*}
where $C$ is a generic constant. Therefore, 
\begin{align*}
   \abs{\int_0^{\eps^{-2s}t} a^s \rd w - \int_0^{\eps^{-2s}t} (e^{-w}|k|)^{2s} \rd w}
\leq
   C \eps^s \vpran{|\xi_k|^{2s} + |k|^{2s}}. 
\end{align*}
Since 
\begin{align*}
\left|  \int_0^{\eps^{-2s}t} (e^{-w}|k|)^{2s}  \rd w - \frac{|k|^{2s}}{2s} \right| = \int_{\eps^{-2s}t}^\infty  (e^{-|w|}|k|)^{2s} \rd w 
\leq |k|^{2s} e^{-\frac{2st}{\eps^{2s}}} ,  
\end{align*}
we have, for $t \geq t_0$, 
\begin{align} \label{I2-may}
  \abs{\int_0^{\eps^{-2s}t} a^s \rd w - \frac{|k|^{2s}}{2s}}
\leq
  C_{s, t_0} \eps^s \vpran{|\xi|^{2s} + |k|^{2s}}
  + |k|^{2s} e^{-\frac{2st_0}{\eps^{2s}}}. 
\end{align}
To estimate the second term on the right-hand side of~\eqref{I2-may}, we consider two cases. First, if $t_0 \geq \sqrt{\eps^{2s}}$, then 
\begin{align*}
   |k|^{2s} e^{-\frac{2st_0}{\eps^{2s}}}
\leq
  |k|^{2s} e^{-\frac{2s}{\eps^{s}}}
\leq 
  C |k|^{2s} \eps^{2s}.
\end{align*}
where $C$ only depends on $s$. If $t_0 \leq \sqrt{\eps^{2s}}$, then $\eps \geq (t_0)^{1/s}$ and
\begin{align*}
  |k|^{2s} e^{-\frac{2st_0}{\eps^{2s}}}
\leq
  |k|^{2s}
= \eps^{s} |k|^{2s} \frac{1}{\eps^{s}}
\leq 
 C \eps^{s} |k|^{2s},
\end{align*}
where $C = 1/t_0$ only depends on $t_0$. In both cases, we have $|k|^{2s} e^{-\frac{2st_0}{\eps^{2s}}} \leq C_{s, t_0} \eps^s |k|^{2s}$ where $C_{s, t_0}$ only depends on $s, t_0$. Combining \eqref{0621}, \eqref{I3-may} and \eqref{I2-may} gives the desired bound for this lemma.
\end{proof}

A corollary follows immediately from Lemma~\ref{lem1}.
\begin{corollary} \label{col1}
Consider $t > t_0 > 0$ where $t_0$ is a fixed constant. Then 
for any fixed $k$ and $\xi$, there exists a constant $C_{s,t_0}$ that only depends on $s$ and $t_0$ such that 
\begin{align*} 
(1-C_{s,t_0} \eps^{s}) \left(|\xi|^{2s} t + \frac{1}{2s}|k|^{2s}  \right)   
\leq 
  \int_0^{\epst} |e^{-w}(k-\eps \xi) + \eps \xi|^{2s} \rd w  \nn
\leq (1+C_{s,t_0}\eps^{s}) 
\left(|\xi|^{2s} t + \frac{1}{2s}|k|^{2s}  \right) .
\end{align*}
Then for $\eps$ sufficiently small, we have
\begin{align} \label{0527}
\int_0^{\epst} |e^{-w}(k-\eps \xi) + \eps \xi|^{2s} \rd w \geq \half \left(|\xi|^{2s} t + \frac{1}{2s}|k|^{2s}  \right)  .
\end{align}
In addition, 
\begin{align*}
   \abs{e^{-\int_0^{\epst} |e^{-w}(k-\eps \xi) + \eps \xi|^{2s} \rd w} - e^{-|\xi|^{2s} t - \frac{1}{2s}|k|^{2s}}} 
 \leq 
   e^{-\half \left(|\xi|^{2s} t + \frac{1}{2s}|k|^{2s}  \right) } C_{s,t_0} \eps^{s}  [|\xi|^{2s} + |k|^{2s} ].
\end{align*}
\end{corollary}

We are now ready to prove Theorem~\ref{APlimit}. 
\begin{proof}[Proof of Theorem~\ref{APlimit}]
By Lemma~\ref{lem1} and Corollary~\ref{col1}, the $A_1$-term in~\eqref{diff0} satisfies
\begin{align}
 & \quad \,
 \int_{\RR^d} \int_{\RR^d} |\hat f_{in} (\xi, 0)A_1 |^2 \rd k \rd \xi  \nonumber
 \\
 &\leq  
   C_{s,t_0}  \eps^{2s} \int_{\RR^d} \int_{\RR^d}  |\hat f_{in}(\xi,0)|^2 e^{-(|\xi|^{2s}t + \frac{|k|^{2s}}{2s} )}  [|\xi|^{2s} + |k|^{2s}]^2 \rd k \rd \xi
\leq 
  C_{s,t_0}  \eps^{2s} \norm{\rho_{in}}_{L^1}^2.  \label{0211}
\end{align}
We bound $A_2$ as follows. 
\begin{align} \label{A222}
  \norm{A_2}_{L_{k,\xi}^2}^2 
&= \int\int \left|\frac{1}{\average{k}}[\hat f_{in} (\xi, \eps \xi + \expo (k - \eps \xi))- \hat f_{in} (\xi, 0)] \average{k} e^{-\int_0^{\epst} |e^{-w}(k-\eps \xi) + \eps \xi|^{2s} \rd w} \right|^2 \rd k \rd \xi \nonumber
\\ 
& \leq  
  \norm{ \frac{1}{\average{k}}[\hat f_{in} (\xi, \eps \xi + \expo (k - \eps \xi))- \hat f_{in} (\xi, 0)]}_{L^\infty_{k,\xi}}^2
\norm{\average{k} e^{-\int_0^{\epst} |e^{-w}(k-\eps \xi) + \eps \xi|^{2s} \rd w}}_{L_{k,\xi}^2}^2 \nonumber
\\  
& \leq 
  \norm{\eps \abs{\xi \cdot \nabla_k \hat f_{in}}
    + C_{t_0} \eps^{2s} \abs{\nabla_k \hat f_{in}}}_{L^\infty_{k,\xi}}^2
\norm{\average{k} e^{-\int_0^{\epst} |e^{-w}(k-\eps \xi) + \eps \xi|^{2s} \rd w}}_{L_{k,\xi}^2}^2 \nonumber
\\ 
& \leq 
   \left[ \eps \norm{v \cdot \nabla_x f_{in}}_{L^1_{x,v}}  + C_{t_0} \eps^{2s} \norm{  v f_{in}}_{L^1_{x,v}} \right]^2
\norm{\average{k} e^{-\half(|\xi|^{2s}t + \frac{|k|^{2s}}{2s} )}}_{L_{k,\xi}^2}^2  \nonumber
\\ 
& \leq 
   C \left[ \eps^2 \norm{v \cdot \nabla_x f_{in}}_{L^1_{x,v}}^2  + C_{t_0} \eps^{4s} \norm{  v f_{in}}_{L^1_{x,v}}^2 \right],
\end{align}
where we have used that for $t>t_0$, it holds that $e^{-\eps^{-2s} t} \leq C_{t_0} \eps^{2s}$ for the same reason stated in the proof of Lemma~\ref{lem1} following~\eqref{I2-may}. The second last inequality follows from \eqref{0527}.  Combining \eqref{0211} and \eqref{A222} leads to the desired result. 
\end{proof}

At the end of this section, we state a useful property regarding the decay of the derivatives of the equilibrium state $\equilibrium(v)$. 
There properties have been previously shown (see for example \cite{doetsch2012introduction, hawkes1971lower, uchaikin2011chance}). We provide a proof in the appendix for the convenience of the reader. 

\begin{proposition} \label{M-decay}
Derivatives of the equilibrium state $\equilibrium(v)$ satisfies
\begin{equation*}
|\nabla_v^m \equilibrium (v)| \leq C_{d, m} |\equilibrium (v)|, 
\end{equation*}
where $\nabla_v^m \equilibrium = \partial_{v_1}^{\alpha_1} \cdots \partial_{v_d}^{\alpha_d} \equilibrium (v)$, for $\alpha_i \geq 0$ and $\sum_{i=1}^d \alpha_i = m$ and the constant $C_{d,m}$ only depends on $m$ and $d$. 
\end{proposition}

\section{Two technical lemmas}
In this section we group some technical commutator estimates that will appear in the error analysis in the next section.  
\begin{lemma} \label{lemt1}
For any function $\bar g(v)$ smooth enough, the commutator
\begin{align} \label{c1228}
 [\average{v}^{p}, (-\Delta_v)^s] \bar g = \int \frac{\average{w}^{p} - \average{v}^{p}}{|v-w|^{d+2s}} \bar g (w) \rd w 
\end{align}
 admits the following estimate for $-2s< p< \frac{d}{2} + 2s$: 
\begin{align*}
\norm{ [\average{v}^{p}, (-\Delta_v)^s] \bar g}_{L_v^2} \leq C  \normv{\average{v}^{p} \bar g} + 
\normv{\average{v}^{p-1} |\nabla_v \bar g|}.
\end{align*}
\end{lemma}

\begin{proof}
We divide the integral into four parts and estimate them separately: 
\begin{align*}
  [\average{v}^{p}, (-\Delta_v)^s] \bar g  
  = \int_{|v-w|>1 , |w|> \half |v|}  \cdot \rd w  
       + \int_{|v-w|>1 , |w| < \half |v|}  \cdot \rd w  
       + \int_{|v-w|\leq1}  \cdot \rd w.
\end{align*}

When $|v-w| > 1$ and $|w|> \half |v|$,  we have  
\begin{align} \label{part1-0}
& \norm{\int_{|v-w| > 1, |w|> \half |v|}  \frac{|\average{w}^{p} - \average{v}^{p}|}{|v-w|^{d+2s}} |\bar g (w)|  \rd w}_{L_{v}^2} \nonumber 
\\
& \hspace{2cm}
\leq \left\{  \begin{array}{cc} C\norm{\int_{|v-w|>1,|w|> \half |v|} \frac{\average{v}^p |\bar g(w)|}{|v-w|^{d+2s}} \rd w }_{L_{v}^2}  &  \text{ if }~ p< 0 ,
\\ C\norm{\int_{|v-w|>1,|w|> \half |v|} \frac{\average{w}^p |\bar g(w)|}{|v-w|^{d+2s}} \rd w }_{L_{v}^2}  &  \text{ if }~ p \geq 0 .
\end{array} \right.
\end{align}
If $p < 0$,  then \eqref{part1-0} can be estimated as follows. If $\half |v| \leq  |w| \leq  2|v|$, then 
\begin{align} \label{part1-1}
 & \norm{\int_{|v-w| > 1,\half |v| \leq |w| \leq 2|v| }  \frac{|\average{w}^{p} - \average{v}^{p}|}{|v-w|^{d+2s}} |\bar g (w)|  \rd w}_{L_{v}^2}  \nonumber
 \\  & \qquad \leq C \norm{ (\average{v}^p |\bar g(v)|) \ast \left( {\bf 1}_{|v|>1} \frac{1}{|v|^{d+2s}} \right)}_{L_{v}^2}  \leq C \norm{ \average{v}^{p} \bar g}_{L_{v}^2}, \quad p<0,
\end{align}
where we have used Young's inequality and the integrability of $\frac{1}{|v|^{d+2s}}$ for $|v|>1$ in the last inequality.
When $|w| \geq 2|v|$, we have
\begin{align*} 
  & \norm{\int_{|v-w| > 1,  |w| \geq 2|v| }  \frac{|\average{w}^{p} - \average{v}^{p}|}{|v-w|^{d+2s}} |\bar g (w)|  \rd w}_{L_{v}^2} 
\leq 
  C \norm{\int_{|v-w| > 1,  |w| \geq 2|v| }  \frac{\average{v}^p}{|v-w|^{d+2s}} |\bar g (w)|  \rd w}_{L_{v}^2}   \nonumber
\\ 
&\leq  
  C \norm{\int_{|v-w| > 1,  |w| \geq 2|v| }  \frac{1}{|v-w|^{d+2s+p}} \average{w}^p |\bar g (w)|  \rd w}_{L_{v}^2}  \nonumber 
\\
&\leq 
   C \norm{ (\average{v}^p |\bar g(v)|) \ast \left( {\bf 1}_{|v|>1} \frac{1}{|v|^{d+2s+p}}  \right) }_{L_{v}^2}  \nonumber
\\ &  \leq C  \norm{ \average{v}^{p} \bar g(v) }_{L_{v}^2} \norm{\left( {\bf 1}_{|v|>1} \frac{1}{|v|^{d+2s+p}}  \right) }_{L_v^1} \nonumber
\\ 
&\leq 
  C  \norm{ \average{v}^{p} \bar g(v) }_{L_{v}^2} , \qquad \text{for}\quad d+2s+p>d \quad \text{or equivalently} \quad p > -2s,
\end{align*}
where we have used Young's inequality again. The lower bound $p>-2s$ is to guarantee the integrability of $|v|^{-d-2s-p}$. 
If $p \geq 0$, then the estimate \eqref{part1-0} takes the following form:
\begin{align*} 
  & \norm{\int_{|v-w| > 1, |w|> \half |v|}  \frac{|\average{w}^{p} - \average{v}^{p}|}{|v-w|^{d+2s}} |\bar g (w)|  \rd w}_{L_{v}^2} \nonumber
  \\ & \qquad 
  \leq C \norm{ (\average{v}^p \bar g) \ast \left( {\bf 1}_{|v|>1} \frac{1}{|v|^{d+2s}} \right) ]}_{L_{v}^2} \leq   C \norm{\average{v}^p \bar g}_{L_{v}^2}, \quad p \geq 0.
\end{align*}
Similarly, if $|v-w| > 1$ and $|w|< \half |v|$, then 
\begin{align*} 
& \norm{\int_{|v-w| > 1, |w|<\half |v|}  \frac{|\average{v}^{p} - \average{w}^{p}|}{|v-w|^{d+2s}} |\bar g (w)|  \rd w}_{L_{v}^2}
\leq \left\{  \begin{array}{cc} C\norm{\int_{|v-w|>1} \frac{\average{w}^p |\bar g(w)|}{|v-w|^{d+2s}} \rd w }_{L_{v}^2} &  \text{ if }~ p\leq 0 ,
\\ C\norm{\int_{|v-w|>1} \frac{\average{v}^p |\bar g(w)|}{|v-w|^{d+2s}} \rd w }_{L_{v}^2} &  \text{ if }~ p >0 .
\end{array} \right.
\end{align*}
The case when $p\leq 0$ can be estimated exactly in the same way as \eqref{part1-1}, which gives
\begin{align*} 
 \norm{\int_{|v-w| > 1, |w| <\half |v|}  \frac{|\average{w}^{p} - \average{v}^{p}|}{|v-w|^{d+2s}} |\bar g (w)|  \rd w}_{L_{v}^2}
  \leq C \norm{ \average{v}^p \bar g}_{L_{v}^2} , \quad p\leq 0.
\end{align*}
A more involved estimate is needed for $p>0$. Note that in this region, $|v| \geq \frac{2}{3}$ and $|v-w|\geq \half |v|$. Thus we have
\begin{align} \label{part2-2}
& \quad\,  \norm{\int_{|v-w|>1, |w| <\half|v|} \frac{\average{v}^p |\bar g(w)|}{|v-w|^{d+2s}} \rd w }_{L_{v}^2}
\leq  C \norm{{\mathbf{1}}_{|v|>\frac{2}{3}} \int_{|w| <\half|v|} \frac{\average{v}^p |\bar g(w)|}{|v|^{d+2s}} \rd w  }_{L_{v}^2}  \nonumber 
\\ 
& =  C \left\{\int_{|v|>\frac{2}{3}} \average{v}^{2p-2d-4s}  \left[\int_{|w|<\half |v|} |\bar g(w)| \rd w \right]^2  \rd v \right\}^\half  \nonumber
\\ 
&\leq 
  C \left\{\int_{|v|>\frac{2}{3}} \average{v}^{2p-2d-4s}  \left[ \int_{|w| < \half |v|} \average{w}^{-2p} \rd w \right] \norm{\average{w}^p \bar{g}(w)}_{L_{v}^2}^2  \rd v\right\}^\half \nonumber
\\ 
&= C \left\{  \int_{|v|>\frac{2}{3}} \average{v}^{2p-2d-4s}  \left[ \int_{|w| < \half |v|} \average{w}^{-2p} \rd w \right]   \rd v \right\}^\half 
\norm{\average{w}^p \bar{g}(w)}_{L_{v}^2}   \nonumber
\\
&\leq 
   \left\{ \begin{array}{cc}  \text{if}~ p<  \frac{d}{2} & C \left\{   \int_{|v|>\frac{2}{3}} \average{v}^{2p-2d-4s} \average{v}^{d-2p}   \rd v  \right\}^\half \norm{\average{w}^p \bar{g}(w)}_{L_{v}^2} ,
\\  \text{if}~ p >  \frac{d}{2}  & C \left\{   \int_{|v|>\frac{2}{3}} \average{v}^{2p-2d-4s}    \rd v  \right\}^\half 
\norm{\average{w}^p \bar{g}(w)}_{L_{v}^2} ,
\\ \text{if}~ p = \frac{d}{2} & C \left\{   \int_{|v|>\frac{2}{3}} \average{v}^{2p-2d-4s} \log\average{v}   \rd v  \right\}^\half 
\norm{\average{w}^p \bar{g}(w)}_{L_{v}^2}.
\end{array} \right.
\end{align}
The integrals in \eqref{part2-2} converge since $p < \frac{d}{2} + 2s$. Combining \eqref{part1-1}--\eqref{part2-2}, 
we have
\begin{align} \label{part12}
 \norm{\int_{|v-w| > 1}  \frac{|\average{w}^{p} - \average{v}^{p}|}{|v-w|^{d+2s}} |\bar g (w)|  \rd w}_{L_{v}^2}
 \leq C \normv{\average{v}^{p} \bar g}, \qquad -2s <p < \frac{d}{2} + 2s.
\end{align}

When  $|v-w| \leq 1$, we consider $s<\half$ and $s\geq \half$ separately. 
Since  $|v-w| \leq 1$, we have 
\begin{equation} \label{0618}
|\average{w}^p - \average{v}^p| \leq C \average{\xi}^{p-1}|v-w|\,,
\end{equation}
for any $\xi$ between $v$ and $w$, where $C$ is independent of $\xi, v, w$.
Thus for $s < 1/2$,
\begin{align}  \label{part3}
& \quad \, 
\norm{\int_{|v-w| \leq 1}  \frac{|\average{w}^{p} - \average{v}^{p}|}{|v-w|^{d+2s}} |\bar g (w)|  \rd w}_{L_{v}^2} \nonumber
 \leq  
   C \norm{\int_{|v-w| \leq 1} \frac{\average{w}^{p-1}}{|v-w|^{2s+d-1}} |\bar g(w)| \rd w}_{L_{v}^2}  
 \nonumber
\\ 
&\leq  
  C \norm{ (\average{v}^{p-1} |\bar g| )\ast \frac{{\bf 1}_{|v|<1}}{|v|^{2s+d-1}}}_{L_{v}^2} \leq C \norm{\average{v}^{p-1}\bar g}_{L_{v}^2}, 
\qquad s< \half.
\end{align}

For $s\geq \half$, we rewrite \eqref{c1228} in the region $|v-w| \leq 1$ as 
\begin{align*}
 \int_{|v-w| \leq 1} \frac{\average{w}^{p} - \average{v}^{p}}{|v-w|^{d+2s}} \bar g (w) \rd w &= 
 \int_{|v-w| \leq 1} \frac{\average{w}^{p} - \average{v}^{p}}{|v-w|^{d+2s}} [ \bar g (w) -\bar g(v)]\rd w
+  \int_{|v-w| \leq 1} \frac{\average{w}^{p} - \average{v}^{p}}{|v-w|^{d+2s}}  \bar g(v) \rd w \nonumber 
\\ 
& =: I_1 + I_2.
\end{align*}
We start with the estimate of $I_2$. Write 
\begin{equation*}
\average{w}^p = \average{v}^p + p \average{v}^{p-2} v\cdot (w-v) + \half D : (w-v)\otimes (w-v) , 
\end{equation*}
where $D$ is the Hessian matrix
\begin{align*}
D = p(p-2) \average{\zeta}^{p-4} \zeta \otimes \zeta, \quad \text{where~} |\zeta| \in \lfloor  |v|, |w|\rceil.
\end{align*}
If $\half |v|< |w|<2 |v|$, then $|D| \leq C \average{v}^{p-2}$. If $|w|\leq \half |v|$ or $|w|>2|v|$, then both $|v|$ and $|w|$ are bounded and it holds that $|D| \leq C \average{v}^{p-2}$. Altogether, we have
\begin{align}\label{2022-I2}
\norm{I_2}_{L_{v}^2} \leq C \norm{\average{v}^{p-2} \bar g}_{L_{v}^2} .
\end{align}

Now we estimate $I_1$. Expand $\bar g(w)$ as
\begin{equation*}
\bar g(w) = \bar g (v) + \int_0^1 \nabla_v \bar g(tw  + (1-t) v) \cdot (w-v) \rd t. 
\end{equation*}
Then 
\begin{align*} 
\norm{I_1}^2_{L_{v}^2}
&\leq  
  \norm{ \int_{|v-w| \leq 1} \int_0^1 |  \nabla_v \bar g(tw  + (1-t) v)|  \frac{\average{\zeta}^{p-1}}{|v-w|^{d+2s-2}} \rd t \rd w}_{L_{v}^2}^2 \qquad   \nonumber 
\\ 
&\leq 
  C \int  \int_0^1 \left[ \int_{|v-w| \leq 1}  \frac{\average{\zeta}^{p-1}}{|v-w|^{d+2s-2}}  |  \nabla_v \bar g(tw  + (1-t) v)|  \rd w\right]^2 \rd t \rd v,  
\end{align*}
where $|\zeta|$ is between $ |v|$ and $|w|$ and the inequalities follow from Cauchy-Schwarz. Let
\begin{equation*}
  z = t w + (1-t) v 
\quad \Rightarrow \quad 
  t^d \rd w = \rd z, \quad |v-z| = t |v-w|.
\end{equation*}
Note that  \eqref{0618} holds when $|v-w|\leq1$, thus 
\begin{align} \label{2203}
\norm{I_1}^2_{L_{v}^2}  
 & \leq C \int\!\! \int  \int_0^1 \left[ \int_{|v-z| \leq t}  \frac{\average{z}^{p-1}}{|v-z|^{d+2s-2}} |\nabla_z \bar g(z)| t^{d+2s-2} t^{-d}  \rd z \right]^2  \rd t \rd v \nonumber
 \\ & = C \int_0^1 \norm{ \int_{|v-z| \leq t}  \frac{\average{z}^{p-1}}{|v-z|^{d+2s-2}} |\nabla_z \bar g(z)| t^{d+2s-2} t^{-d}  \rd z }_{L_{v}^2}^2  \rd t  \nonumber
 \\ & \leq C \int_0^1 \norm{{\mathbf 1}_{|z|\leq t} \frac{1}{|z|^{d+2s-2}}}_{L^1_{x,v}}^2 \norm{\average{z}^{p-1} |\nabla_z \bar g(z)|}_{L_{v}^2}^2  t^{4s-4}\rd t \nonumber
 \\ & = C \norm{\average{z}^{p-1} |\nabla_z \bar g(z)|}_{L_{v}^2}^2 .
\end{align}
Combining \eqref{2022-I2} and \eqref{2203}, we have
\begin{align} \label{part34-2}
\norm{ \int_{|v-w| \leq 1}  \frac{|\average{w}^{p} - \average{v}^{p}|}{|v-w|^{d+2s}} |\bar g (w)|  \rd w}_{L_{v}^2}
\leq 
  C \norm{ \average{v}^{p-2} \bar g}_{L^2_v} +  \norm{\average{v}^{p-1} |\nabla_v \bar g|}_{L^2_v} \, \quad s\geq \half .
\end{align}
The final result follows from \eqref{part12}, \eqref{part3} and \eqref{part34-2}. 
\end{proof}

The second lemma in this section is a commutator estimate for weights and the operator $\Lop^s$. 
\begin{lemma} \label{lemma:commutator}
For $0 \leq b \leq d+ 2s$, if $u$ satisfies 
\begin{equation} \label{hlessM}
|u| \leq \frac{C}{\eps^{2s}} \equilibrium ,
\end{equation}
then for $0<\delta < 1$ the commutator
\begin{align*}
[(1+\delta \average{v})^{-b}, \Lop^s] u & := - \left[ \nabla_v \cdot (v(1+\delta \average{v})^{-b} u) - (1+\delta \average{v})^{-b} \nabla_v \cdot (v u) \right] 
\\ &  \qquad + \left[ (-\Delta_v)^s ((1+\delta \average{v})^{-b}u) - (1+\delta \average{v})^{-b} (-\Delta_v)^s u \right]
\end{align*}
satisfies the bound
\begin{align*}
 \left|\int \!\! \int (1+\delta \average{v})^{-b} u [(1+\delta \average{v})^{-b}, \Lop^s] u \frac{1}{\equilibrium} \rd v \rd x \right|
\leq C\delta^{\frac{2s}{d+4}}  \normM{(1+\delta \average{v})^{-b} u}^2 + C  {\color{black} \eps^{-4s}} \delta^{\min\{1,2s\}} .
\end{align*}
\end{lemma}
\begin{proof}
Denote 
\begin{align*}
& T_1 u (v)=   \nabla_v \cdot (v(1+\delta \average{v})^{-b} u) - (1+\delta \average{v})^{-b} \nabla_v \cdot (v u), 
\\ &T_2 u (v) =  (-\Delta_v)^s ((1+\delta \average{v})^{-b} u) - (1+\delta \average{v})^{-b} (-\Delta_v)^s u. 
\end{align*}
Then $[(1+\delta \average{v})^{-b}, \Lop^s] u = -T_1 u + T_2 u .$
Note that 
\begin{align*}
T_1 u = v u \cdot  \nabla_v[(1+\delta \average{v})^{-b} ]  = -b u (1+\delta \average{v})^{-b-1}  \delta \frac{|v|^2}{\average{v}}.
\end{align*}
Then from the assumption \eqref{hlessM}, we have 
\begin{align*}
& \quad \, 
\left| \int (1+\delta \average{v})^{-b}  u T_1 u \frac{1}{\equilibrium} \rd v\right| 
\leq Cb  \eps^{-4s} \int \frac{\delta \average{v}}{1+ \delta \average{v}} \equilibrium \rd v  \nonumber
\\ &= Cb   \eps^{-4s} \int_{|v|<1} \frac{\delta \average{v}}{1+ \delta \average{v}} \equilibrium \rd v
+ Cb   \eps^{-4s} \int_{|v| \geq 1} \frac{\delta \average{v}}{1+ \delta \average{v}} \equilibrium \rd v \,. 
\end{align*}
It is immediate that 
\begin{align*}
Cb  \eps^{-4s} \int_{|v|<1} \frac{\delta \average{v}}{1+ \delta \average{v}} \equilibrium \rd v 
\leq C b\eps^{-4s} \delta \int_{|v|<1} \average{v} \equilibrium \rd v \leq C b\eps^{-4s} \delta.
\end{align*}
Moreover, 
\begin{align*}
& \quad \, 
  Cb   \eps^{-4s} \int_{|v| \geq 1} \frac{\delta \average{v}}{1+ \delta \average{v}} \equilibrium \rd v
\\ 
&= Cb  {\color{black} \eps^{-4s}} \int_{1\leq |v| \leq \frac{1}{\delta}} \frac{\delta \average{v}}{1+ \delta \average{v}} \equilibrium \rd v + 
Cb  \eps^{-4s} \int_{|v| \geq \frac{1}{\delta}} \frac{\delta \average{v}}{1+ \delta \average{v}} \equilibrium \rd v 
\\ 
&\leq 
  Cb  \eps^{-4s} \delta  \int_{1\leq |v| \leq \frac{1}{\delta}} \average{v} \equilibrium \rd v + 
Cb \eps^{-4s} \int_{|v| \geq \frac{1}{\delta}}  \equilibrium \rd v  \leq C b  \eps^{-4s}  \delta^{2s}.
\end{align*}
Therefore, 
\begin{align} \label{T1}
\left| \int (1+\delta \average{v})^{-b}  u T_1 u \frac{1}{\equilibrium} \rd v\right| 
\leq  Cb  {\color{black} \eps^{-4s}} (\delta^{2s} + \delta) .
\end{align}

To estimate $T_2u$, we first rewrite it as
\begin{align*}
T_2 u (v) = C_{s,d} \int u(w) \frac{(1+\delta \average{v})^{-b}  - (1+\delta \average{w})^{-b}  }{|v-w|^{d+2s} } \rd w.
\end{align*}
Hence, 
\begin{align*} 
& \int (1+\delta \average{v})^{-b}  u T_2 u \frac{1}{\equilibrium (v)} \rd v  \nonumber 
\\ & =  C_{s,d} \int \int \frac{(1+\delta \average{v})^{-b}  - (1+\delta \average{w})^{-b}  }{|v-w|^{d+2s} } u(w) \frac{u(v)}{(1+\delta \average{v})^b}  \frac{1}{\equilibrium (v)} \rd v \rd w  \nonumber
\\ 
& =: C_{s,d}  \int \int J \rd v \rd w.
\end{align*}
We proceed by separating the integration domain into $|v-w|\leq C_0$ and $|v-w|>C_0$ for some constant $C_0>1$ to be determined and estimate the two integrals individually. Over the domain $|v-w|> C_0$, we have 
\begin{align*}
 & \quad \, 
 \left| {\int \int}_{|v-w|> C_0} J \rd v \rd w \right| 
 \\  
 &\leq {\int \int}_{|v-w|> C_0} \frac{| (1+\delta \average{v})^{b}  - (1+\delta \average{w})^{b} |}{|v-w|^{d+2s}} \frac{1}{(1+\delta \average{v})^{b} } \frac{|u(v)|}{(1+\delta \average{v})^{b} } \frac{|u(w)|}{(1+\delta \average{w})^{b} } \frac{1}{\equilibrium (v)} \rd w \rd v
\\  
&=  {\int \int}_{|v-w|> C_0, |w|< \half |v|} +  {\int \int}_{|v-w|> C_0, |w|>2|v|} +  {\int \int}_{|v-w|> C_0, \half |v| \leq |w| \leq 2 |v|}
\\  
&=: T_{23} + T_{24} + T_{25}.
\end{align*}
If $|v-w|> C_0$ and $|w| < \half |v|$, then $|v|>\frac{2}{3} C_0$ and 
\begin{align*}
 \frac{| (1+\delta \average{v})^{b}  - (1+\delta \average{w})^{b} |}{|v-w|^{d+2s}} \frac{1}{(1+\delta \average{v})^{b} }  \frac{1}{\equilibrium(v)}
 \leq 
   C \frac{(1+\delta \average{v})^b}{|v|^{d+2s}}   \frac{1}{(1+\delta \average{v})^{b} }  \frac{1}{\equilibrium(v)} 
 \leq C.
\end{align*}
Thus, 
\begin{align}
  T_{23} 
&\leq  
  C {\int \int}_{|v| > \frac{2}{3} C_0}  \frac{|u(w)|}{(1+\delta \average{w})^{b} } \frac{|u(v)|}{(1+\delta \average{v})^{b} }  \rd w \rd v \nonumber 
\\ 
&\leq 
  C \int \frac{|u(w)|}{(1+\delta \average{w})^{b}} \frac{1}{\sqrt{\equilibrium(w)}} \sqrt{\equilibrium(w)}\rd w 
                 \int_{|v| > \frac{2}{3} C_0} \frac{|u(v)|}{(1+\delta \average{w})^{b}} \frac{1}{\sqrt{\equilibrium(v)}} \sqrt{\equilibrium(v)}\rd v \nonumber
\\ 
&\leq 
  C \normM{(1+\delta \average{v})^{-b}  u}^2  \left( \int_{|v| > \frac{2}{3} C_0} \equilibrium(v) \rd v \right)^{\half}
\nonumber 
\\ 
&\leq 
  C  C_0^{-s} \normM{(1+\delta \average{v})^{-b}  u}^2.  \label{T23}
\end{align}

If $|v-w|> C_0$ and $|w| > 2 |v|$, then $|w| > \frac{2}{3} C_0$ and $|v-w| \geq |w|-|v| \geq \frac{1}{2}|w|$, thus 
\begin{align} \label{1208}
& \quad\,  
\frac{| (1+\delta \average{v})^{b}  - (1+\delta \average{w})^{b} |}{|v-w|^{d+2s}} \frac{1}{(1+\delta \average{v})^{b}}  \frac{1}{\equilibrium(v)} 
 \nonumber
 \\ 
 &\leq 
   C \frac{(1+\delta \average{w})^{b}}{|w|^{d+2s} }\frac{1}{(1+\delta \average{v})^{b}}  \frac{1}{\equilibrium(v)} 
  \leq 
    C \frac{ (1+\delta \average{w})^{b} \average{w}^{-d-2s}}{  (1+\delta \average{v})^{b} \average{v}^{-d-2s} },
\end{align}
where in the last inequality we have used $C_1 \average{v}^{-d-2s} \leq \equilibrium(v) \leq C_2 \average{v}^{-d-2s}$ for $0 < C_1 \leq C_2$.
If we choose $b \leq d+ 2s$,
then $(1+\delta \average{w})^{b} \average{w}^{-d-2s}$ is a decreasing function in $|w|$ and \eqref{1208} implies 
\begin{align*}
 \frac{| (1+\delta \average{v})^{b}  - (1+\delta \average{w})^{b} |}{|v-w|^{d+2s}} \frac{1}{(1+\delta \average{v})^{b}}  \frac{1}{\equilibrium(v)}  \leq C.
\end{align*}
The rest of  the estimate is the same as  that  of $T_{23}$ and we have 
\begin{align} \label{T24}
T_{24}  \leq C  C_0^{-s} \normM{(1+\delta \average{v})^{-b}  u}^2.
\end{align}

If $|v-w|> C_0$ and $\half |v| < |w| < 2 |v|$, then 
\begin{align*}
 \frac{| (1+\delta \average{v})^{b}  - (1+\delta \average{w})^{b} |}{(1+\delta \average{v})^b}  \leq C.
\end{align*}
Hence, 
\begin{align}
  T_{25} 
& \leq  
  C \int\int_{|v-w|>C_0, \half |v|<|w|<2|v|} \frac{1}{|v-w|^{2s+d}}  \frac{|u(w)|}{(1+\delta \average{w})^b} \frac{|u(v)|}{(1+\delta \average{v})^b}  \frac{1}{\equilibrium(v)} \rd w \rd v \nonumber
\\ 
&\leq 
  C \int\int_{|v-w|>C_0, \half |v|<|w|<2|v|} \frac{1}{|v-w|^{2s+d}}  \frac{|u(w)|}{(1+\delta \average{w})^b} \frac{|u(v)|}{(1+\delta \average{v})^b}  \frac{1}{\sqrt{\equilibrium(v)}} \frac{1}{\sqrt{\equilibrium(w)}} \rd w \rd v \nonumber
\\ 
&\leq 
  C C_0^{-2s} \normM{(1+\delta \average{v})^{-b}  u}^2.  \label{T25}
\end{align}

Combining \eqref{T23}, \eqref{T24} and \eqref{T25}, we have 
\begin{align}
\left| {\int \int}_{|v-w| > C_0} J  \rd v \rd w \right|   \leq C C_0^{-s} \normM{(1+\delta \average{v})^{-b}  u}^2 .   \label{Cg}
\end{align}

For the domain $|v-w| \leq C_0$, we symmetrize the integral as 
\begin{align*}
&  {\int \int}_{|v-w|\leq C_0} J \rd v \rd w 
\\ =&~ \half {\int \int}_{|v-w|\leq C_0} \frac{(1+\delta \average{v})^{b}  - (1+\delta \average{w})^{b}  }{|v-w|^{d+2s} } \frac{u(w)}{(1+\delta \average{w})^b}  \frac{u(v)}{(1+\delta \average{v})^b} \times
\\ 
& \hspace{6cm} 
  \times \left[ \frac{1}{\equilibrium(v) (1+\delta \average{v})^b} - \frac{1}{\equilibrium(w) (1+\delta \average{w})^b}  \right] \rd v \rd w.
\end{align*} 
By the mean value theorem, for $\average{\eta_1}, \average{\eta_2}$ between  $ \average{v}$ and $ \average{w}$, we have
\begin{align} \label{2207}
& \quad \, \left| {\int \int}_{|v-w|\leq C_0} J \rd v \rd w  \right|  \nonumber
\\ 
&\leq 
  b \delta \int\!\!\int \frac{(1+\delta \average{\eta_1})^{b-1} \frac{|\eta_1|}{\average{\eta_1}}}{|v-w|^{d+2s-2}} 
\frac{|u(w)|}{(1+\delta \average{w})^b}  \frac{|u(v)|}{(1+\delta \average{v})^b} \times  \nn
\\
& \hspace{4cm}
\times
\frac{1}{\equilibrium(\eta_2) (1+\delta \average{\eta_2})^b} 
\left| \frac{\nabla \equilibrium(\eta_2)}{\equilibrium(\eta_2)} + \frac{b\delta  \frac{\eta_2}{\average{\eta_2} }}{(1+\delta \average{\eta_2})}  \right| \rd v \rd w  \nonumber
\\ 
& \leq 
  b \delta \int\!\!\int \frac{1}{|v-w|^{d+2s-2}} 
\frac{|u(w)|}{(1+\delta \average{w})^b}  \frac{|u(v)|}{(1+\delta \average{v})^b} \times  \nn
\\
& \hspace{4cm}
\times\frac{(1+\delta \average{\eta_1})^{b-1} }{\equilibrium(\eta_2) (1+\delta \average{\eta_2})^b} \left| \frac{\nabla \equilibrium(\eta_2)}{\equilibrium(\eta_2)} + \frac{b\delta  \frac{\eta_2}{\average{\eta_2} }}{1+\delta \average{\eta_2}}  \right| \rd v \rd w  \nonumber
\\ 
&\leq b \delta(1+C_0\delta)^{\max\{b-1,0\}} \int\!\!\int \frac{1}{|v-w|^{d+2s-2}} 
   \frac{|u(w)|}{(1+\delta \average{w})^b}  \frac{|u(v)|}{(1+\delta \average{v})^b} \frac{1}{\equilibrium(\eta_2)} \rd v \rd w, 
\end{align}
where we have used $\abs{\frac{\nabla \equilibrium(\eta_2)}{\equilibrium(\eta_2)}} \leq C$ thanks for Proposition~\ref{M-decay}. In addition, if $\half |v| < |w| \leq 2|v|$, then 
\begin{align*}
   \average{\eta}\ \leq C \average{v} 
\quad \text{and} \quad
   \frac{(1+\delta \average{\eta})^{b-1}}{(1+\delta \average{v})^{b} } 
\leq C. 
\end{align*} 
If $|w| \geq 2|v|$, then we have $2|v| \leq |w| \leq 2C_0 $ and
\begin{equation*} 
(1+\delta \average{v})^{b-1} \leq \left\{ \begin{array}{cc} C (1+\delta C_0)^{b-1} & b>1, \\ C & b\leq 1. \end{array} \right. 
\end{equation*}
Likewise, when $|w| \leq \frac{1}{2} |v|$, we have $2|w| \leq |v| \leq 2C_0 $ and the similar bound holds.
Therefore, we always have 
\begin{align*}
  \frac{(1+\delta \average{\eta_1})^{b-1} }{ (1+\delta \average{\eta_2})^b} \leq 
  C (1+C_0\delta)^{\max\{b-1,0\}}.
\end{align*} 
Since $|v-w| \leq C_0$, if $|w| > 2 |v|$, then we have $|v|\leq C_0$ and $|w| \leq 2 C_0$. If $|w| < \half |v|$, then $|v| \leq 2C_0$ and $|w| \leq C_0$. In both cases,
\begin{equation} \label{2208}
  \frac{1}{\equilibrium(\eta_2)} 
\leq 
   C \max \left\{\frac{1}{\equilibrium(w)}, \frac{1}{\equilibrium(v)} \right\} \leq 
   C \frac{1}{\sqrt{\equilibrium(w)}} \frac{1}{\sqrt{\equilibrium(v)}} C_0^{\frac{d}{2}+s}.
\end{equation}
When $\half|v| \leq |w|\leq 2|v|$, $\average{\eta_2} $ is comparable to both $ \average{v}$ and $ \average{w}$. Therefore, we have
\begin{equation} \label{2209}
\frac{1}{\equilibrium(\eta_2)} \leq  C \frac{1}{\sqrt{\equilibrium(w)}} \frac{1}{\sqrt{\equilibrium(v)}}.
\end{equation}
Substituting \eqref{2208} and \eqref{2209} into \eqref{2207}, we have 
\begin{align} \label{CI-2}
 & \quad \, 
 \left| {\int \int}_{|v-w|\leq C_0} J \rd v \rd w  \right|  \nonumber
 \\  
 &\leq 
   b \delta(1+C_0\delta)^{\max\{b-1,0\}} C_0^{\frac{d}{2}+s} \int\!\!\int_{|v-w|< C_0} \frac{1}{|v-w|^{d+2s-2}} 
   \frac{|u(w)|}{(1+\delta \average{w})^b}  \times  \nn
\\
& \hspace{7.5cm} \times   
   \frac{|u(v)|}{(1+\delta \average{v})^b}  \frac{1}{\sqrt{\equilibrium(w)}} \frac{1}{\sqrt{\equilibrium(v)}} \rd v \rd w \nonumber
\\ 
&\leq Cb \delta (1+\delta C_0)^{\max\{b-1,0\}} C_0^{\frac{d}{2}+s} \left(  \int \!\! \int_{|v-w|<C_0}  \frac{1}{|v-w|^{2s+d-2}} \frac{|u(v)|^2}{(1+\delta \average{v})^{2b} } \frac{1}{\equilibrium(v)} \rd v \rd w \right)^\half \times \nonumber
\\  
& \hspace{5cm} 
  \times \left(  \int \!\! \int_{|v-w|<C_0}  \frac{1}{|v-w|^{2s+d-2}} \frac{|u(w)|^2}{(1+\delta \average{w})^{2b} }\frac{1}{\equilibrium(w)}  \rd v \rd w \right)^\half  \nonumber
\\ 
&\leq  Cb \delta (1+\delta C_0)^{\max\{b-1,0\}} C_0^{\frac{d}{2}+s} C_0^{2-2s} \normM{{(1+\delta\average{v})^{-b}} u }^2 .
\end{align}

Summation of \eqref{Cg} and  \eqref{CI-2} leads to
\begin{align} 
 \int\!\! \int (1+\delta \average{v})^{-b}  u T_2 u \frac{1}{\equilibrium} \rd x \rd v \nonumber
 \leq 
   C \left[ b\delta (1+\delta C_0)^{\max\{b-1,0\}} C_0^{\frac{d}{2}+2-s}  + C_0^{-s} \right] \normM{(1+\delta \average{v})^{-b}  u}^2.
\end{align}

Choose $C_0 = \delta^{-\frac{2}{d+4}}$. 
Then the inequality above becomes
\begin{align*}
 \int \!\! \int  (1+\delta \average{v})^{-b}  u T_2 u \frac{1}{\equilibrium} \rd x \rd v  \leq C \delta^{\frac{2s}{d+4}} \normM{(1+\delta \average{v})^{-b}  u}^2.
\end{align*}
Together with \eqref{T1}, this implies
\begin{align*}
\left|\int \!\! \int (1+\delta \average{v})^{-b} u [(1+\delta \average{v})^{-b}, \Lop^s] u \frac{1}{\equilibrium} \rd x  \rd v \right|
\leq  C\delta^{\frac{2s}{d+4}}  \normM{(1+\delta \average{v})^{-b} u}^2 + C  {\color{black} \eps^{-4s}} \delta^{\min\{1,2s\}} .
\end{align*}
We thereby finish the proof of the commutator estimate. 
\end{proof}

\section{Well-posedness}
In this short section, we present the basic well-posedness results for system~\eqref{semi_scheme}. We start with the properties of operators involved in these systems. 

Let the operator $\Gamma_a$ be defined as a Fourier multiplier:
\begin{align*}
   \Gamma_a = \vpran{I + a (-\Delta_x)^s}^{-1}, 
\qquad
   a > 0, 
\quad x \in \RR^d,
\end{align*}
In particular, for any $h \in H^k(\RR^d)$ with $k \in \RR$, 
\begin{align*}
   \widehat {\Gamma_a h} =  \frac{1}{1 + a |\xi|^{2s}} \hat h. 
\end{align*}

Then we have
\begin{proposition} \label{prop:basic-operators}
(a) The operator $\Gamma_a$ is bounded on $L^\infty(\RR^d)$ and $H^k(\RR^d)$ for any $k \in \RR$. Moreover, 
\begin{align*}
   \norm{\Gamma_a}_{H^k \to H^k} \leq 1, 
\qquad
   \norm{\Gamma_a}_{L^\infty \to L^\infty} \leq 1, 
\end{align*}
Similar bound holds on $W^{m, \infty}(\RR^d)$ for any $m \geq 0$.

\Ni (b) Denote $L^2({\mathcal{M}^{-1} \rd v})$ as the space where
\begin{align*}
   \norm{g}^2_{L^2({\mathcal{M}^{-1} \rd v})}
= \int_{\RR^d} \frac{g^2}{\mathcal{M}} \rd v < \infty.   
\end{align*}
Then any $\lambda > 0$ is in the resolvent set of $\mathcal{L}^s$ on $L^2({\mathcal{M}^{-1} \rd v})$ and 
\begin{align*}
   \norm{(\lambda - \mathcal{L}^s)^{-1}}_{L^2({\mathcal{M}^{-1} \rd v}) \to L^2({\mathcal{M}^{-1} \rd v})} \leq 1/\lambda.
\end{align*}

\end{proposition}
\begin{proof}
(a) First, the $H^k$-bound can be seen from the Fourier transform and the Parseval's identity. Second, suppose $g \in C^2(\RR^d) \cap L^\infty(\RR^d)$ is the solution to 
\begin{align*}
  g = - a (-\Delta_x)^{s} g + h. 
\end{align*}
and $g$ attains its maximum at $x_0$. Then by \eqref{def_fl}
\begin{align*}
  \norm{g}_{L^\infty} 
= g(x_0) 
= - a \, C_{s, d} \text{P.V.} \int_{\mathbb{R}^d} \frac{g(x_0)-\tilde g(w)}{|x_0-w|^{d + 2s}} \rd w + h(x_0)
\leq  h(x_0) \leq \norm{h}_{L^\infty}.
\end{align*}
General $L^\infty$-bound can be deduced by a regular approximation via convolution. The $W^{m \infty}$-bound holds since the operator $\nabla_x^k$ commutes with $\Gamma_a$.

(b) The invertibility of $\lambda - \mathcal{L}^s$ for $\lambda > 0$ is shown in~\cite{cesbron2012anomalous}. In particular it follows from a theorem of Lions~\cite{JLL1961}, which is recalled in Appendix A of~\cite{Degond1986}. It can also been seen by the Fourier transform as $\lambda g = \mathcal{L}^s g + h$ is transformed into
\begin{align} \label{eq:FT-transport}
   (\lambda + |\xi|^{2s}) \, \hat g - \xi \cdot \nabla_\xi \hat g = \hat h, 
\end{align}
where $\hat g, \hat h$ are the Fourier transform of $g, h$. For $\lambda > 0$, this is a linear transport equation with a Lipschitz velocity field. Hence its solution is well-defined via characteristics. Moreover, if $h$ is smoother, then $g$ is smoother. For example, if $D_v^2 h \in L^2({\mathcal{M}^{-1} \rd v})$, then $|\xi|^2 \hat g$ satisfies
\begin{align*}
   (\lambda + |\xi|^{2s} + 2) \, (|\xi|^2 \hat g) - \xi \cdot \nabla_\xi (|\xi|^2 \hat g) = |\xi|^2 \hat h, 
\end{align*}
which shows that $D^2_v g \in L^2({\mathcal{M}^{-1} \rd v})$. The operator bound can be derived from the energy estimate: suppose $g, h$ are smoother enough such that 
\begin{align*}
   \lambda g = \mathcal{L}^s g + h
\end{align*}
holds in a strong sense. Multiply the equation by $g/M$ and integrate in $v$. Then
\begin{align*}
   \lambda \norm{g}^2_{L^2({\mathcal{M}^{-1} \rd v})}
= \int_{\RR^d} \frac{g}{\mathcal{M}} \mathcal{L}^s g \, \rd v
   + \int_{\RR^d} f g \frac{1}{\mathcal{M}} \, \rd v
\leq 
   \int_{\RR^d} f g \frac{1}{\mathcal{M}} \, \rd v
\leq
  \norm{g}_{L^2({\mathcal{M}^{-1} \rd v})}
  \norm{h}_{L^2({\mathcal{M}^{-1} \rd v})}. 
\end{align*}
following the non-negativity of the operator $\mathcal{L}^s$ in $L^2({\mathcal{M}^{-1} \rd v})$ (see \cite{cesbron2012anomalous} for a proof). Therefore, $\norm{g}_{L^2({\mathcal{M}^{-1} \rd v})} \leq (1/\lambda) \norm{h}^2_{L^2({\mathcal{M}^{-1} \rd v})}$ for $h, g$ smooth enough. Similar bound for general $L^2({\mathcal{M}^{-1} \rd v})$ solutions follows from the linearity of the equation and a density argument by smoothing $h$. 
\end{proof}

The well-posedness of~\eqref{semi_scheme} follows immediately from Proposition~\ref{prop:basic-operators}. In particular, we have 
\begin{lemma} \label{lem:well-posedness-scheme}
(a) The solution $\eta^n$ to \eqref{semis3} with initial condition $\eta^0(x,v) = \rho_{in}(x+\eps v)$ has the form
\begin{equation} \label{heqn0}
\eta^n(x,v) = h^n(x+\eps v)\,,
\end{equation}
where  $h^n(y)$ solves
\begin{equation} \label{heqn}
 \frac{1}{\Delta t} (h^{n+1}-h^n)  =  - (-\Delta_y)^s h^{n+1}, \quad h^0(y) = \rho_{\text{in}}(y)\,.
\end{equation} 
Moreover, we have
\begin{equation*}
\norm{h^n}_{H_x^1} \leq \norm{h^0}_{H_x^1}=\norm{\rho_{in}}_{H_x^1},  \qquad
\norm{h^n}_{W_x^{1,\infty}} \leq \norm{h^0}_{W_x^{1,\infty}} = \norm{\rho_{in}}_{W_x^{1,\infty}}.
\end{equation*}

\Ni (b) Given $\eta^{n+1}(x, v) = h^{n+1}(x + \eps v)$ with $h \in H^1_x$, solutions $g^{n+\half}$ and $g^{n+1}$ to~\eqref{semis1} and~\eqref{semis2} are both well-defined in the space $L^2_{\mathcal{M}^{-1}}$ given in~\eqref{normM}. 
\end{lemma}
\begin{proof}
(a) The structure in~\eqref{heqn0} can be seen from induction. First, it holds for $n=0$. Suppose $\eta^n$ has the desired form. Then by~\eqref{scheme22}, 
\begin{align*}
   \eta^{n+1} (x, v) 
= (I + \Delta t (-\Delta_x)^s)^{-1} \eta^n (x, v)
&= (I + \Delta t (-\Delta_x)^s)^{-1} \vpran{h^n(x + \eps v)}
\\
& = \vpran{(I + \Delta t (-\Delta_x)^s)^{-1} h^n} (x + \eps v),
\end{align*}
by the translation invariance of $(I + \Delta t (-\Delta_x)^s)^{-1}$. The bounds follow from Proposition~\ref{prop:basic-operators}. 

(b) By Lemma~\ref{lemma2} part 1) (which is independent of the current lemma), 
we have that the source term in equation~\eqref{semis1} for $g^{n+\half}$ is in $L^2_{\mathcal{M}^{-1}}$. Note that equation~\eqref{semis1} only acts on $v$-variable and $x$ is a parameter. Therefore, for $a.e.$ $x$ where $I(\eta^{n+1}(x, \cdot), \CalM) \in L^2(\CalM^{-1} \dv)$, we have
\begin{align*}
   g^{n+\half}
= \vpran{\tfrac{\eps^{2s}}{\Delta t} + \gamma - \CalL^s}^{-1}
   \vpran{\tfrac{\eps^{2s}}{\Delta t} g^n - I(\eta^{n+1}, \CalM)},
\qquad
   \gamma > 0. 
\end{align*}
By induction, if $g^n \in L^2_{\CalM^{-1}}$, then $g^{n+ \half} \in L^2_{\CalM^{-1}}$.  Using $g^{n+ \half}$ as a source term in equation~\eqref{semis2}, we have
\begin{align*}
  g^{n+1}
= \tfrac{\eps^{2s}}{\Delta t} \vpran{\tfrac{\eps^{2s}}{\Delta t} - \gamma+ \eps v \partial_x}^{-1} g^{n+\half},
\qquad
  \gamma \neq \tfrac{\eps^{2s}}{\Delta t}. 
\end{align*}
where $\tfrac{\eps^{2s}}{\Delta t} - \gamma+ \eps v \partial_x$ is a linear transport operator, so its inverse is well-defined.  Applying the Fourier transform in $x$, we have
\begin{align*}
   \norm{g^{n+1}}_{L^2_{\CalM^{-1}}}
\leq 
  \frac{\eps^{2s}/\Delta t}{\abs{\frac{\eps^{2s}}{\Delta t} - \gamma}} \norm{g^{n+\half}}_{L^2_{\CalM^{-1}}}. 
\end{align*}
Note that the choice of $\gamma$ in \eqref{cond:gamma1}, \eqref{cond:gamma2} and \eqref{gamma-II} always guarantee $\gamma \neq \alpha := \frac{\eps^{2s}}{\Delta t}$. Therefore, $g^n \in L^2_{\CalM^{-1}}$ implies that $g^{n+1} \in L^2_{\CalM^{-1}}$, which closes the induction proof. 
\end{proof}

\section{Error estimate}
In this section we prove the uniform accuracy in time for the semi-discrete scheme \eqref{semi_scheme}. When $\eps$ is large compared with $\Delta t$, it corresponds to the kinetic regime, in which we show that the numerical error can be controlled by some positive power of $\Delta t$. This shows the consistency of the method. When $\eps$ is small compared with $\Delta t$, we will leverage the asymptotic preserving property of our scheme to control the numerical error by a certain positive power of $\eps$. The specific regimes will be made clear later. 

Given the initial data $f_{in}(x,v)$, we have the initial data for $\rho$ and $g$ as
\[
\rho_{in} =  \average{f_{in}}_v, \quad 
g_{in }(x,v) =f_{in}(x,v)-\rho_{in}(x+\eps v) \equilibrium(v)\,,
\]
where the initial data satisfy the regularity condition in~\eqref{cond:reg-initial-1}-\eqref{cond:reg-initial-2}.
Add \eqref{semis1} and \eqref{semis2} to get 
\begin{subequations}\label{scheme2}
\begin{numcases}{}
 \frac{1}{\Delta t} (\eta^{n+1}-\eta^n)  =  - (-\Delta_x)^s \eta^{n+1} , \label{scheme22}
 \\ \frac{\eps^{2s}}{\Delta t}(g^{n+1}-g^n)  + \eps v  \partial_x g^{n+1} = \Lop^s(g^\nh) - \gamma g^\nh + \gamma  g^{n+1} - I(\eta^{n+1}, \equilibrium). \label{scheme21} 
\end{numcases}
\end{subequations}

Let $\egn$ and $\een$ be the local truncation errors defined by 
\begin{align}
   e_1^n &= \frac{1}{\Delta t^2} (g(t^\np) - g(t^n) - \Delta t \, \partial_t g(t^\np)),  \label{e1}
\\
  e_2^n &= \frac{1}{\Delta t^2} (\eta(t^\np) - \eta(t^n) - \Delta t \, \partial_t \eta(t^\np)).  \label{e2}
\end{align}
Then~\eqref{split0} becomes
\begin{subequations}  \label{0327}
\begin{numcases}{}
 \frac{1}{\Delta t} \left[ \eta(\tnp) - \eta(t^n)\right] = -(-\Delta_x)^s \eta (\tnp) + \Delta t e_2^n,
\\ 
\epsdt \left[ g(\tnp) - g(t^n)\right] + \eps v \partial_x g(\tnp) = \Lop^s  g(\tnp) - I(\eta(\tnp), \equilibrium) + \eps^{2s} \Delta t e_1^n.
\end{numcases}
\end{subequations}
Denote 
\begin{equation*} 
\ten(x,v) = \eta(t^n,x,v) - \eta^n(x,v), \qquad \tgn(x,v) = g(t^n,x,v)-g^n(x,v)
\end{equation*}
as the numerical error. By subtracting \eqref{scheme2} from \eqref{0327}, we get 
\begin{subequations}\label{erroreqn}
\begin{numcases}{}
 \frac{1}{\Delta t} (\tenp - \ten) = -(-\Delta_x)^s \tenp + \Delta t e_2^n,  \hspace{10cm}  \label{err_eta}
\\ \epsdt (\tgnp - \tgn) + \eps v \partial_x \tgnp = \Lop^s \tgnp - I(\tenp, \equilibrium) + (\Lop^s - \gamma \id) (g^{n+1}-g^\nh) + \Delta t \eps^{2s} e_1^n  . \label{err_g}
\end{numcases}
\end{subequations}
Denote the error from the operator splitting as 
\begin{equation*} 
e_3^n = (\Lop^s - \gamma \id) (g^{n+1}-g^\nh)\,,
\end{equation*}
and denote
\begin{align} \label{def:alpha-tilde-f}
  \alpha = \frac{\eps^{2s}}{\Delta t}, 
\qquad 
  \tilde f^n = \tilde \eta^n \equilibrium + \tilde g^n.
\end{align}
Multiply \eqref{err_eta} by $\eps^{2s}\equilibrium$ and add to \eqref{err_g}. Then the error equation \eqref{erroreqn} has the form 
\begin{equation} \label{erroreqn2}
\alpha (\ftilde^\np - \ftilde^n) + \eps v \partial_x \ftilde^\np = \Lop^s \ftilde^\np + \Delta t \eps^{2s} e_1^n + e_3^n + \eps^{2s} \Delta t e_2^n \equilibrium.
\end{equation}
In view of \eqref{scheme22}, $\tilde \eta^n$ evolves independently from $g^n$.   We outline the properties of $\tilde \eta^n$ and the estimates $e_1^n$ and $e_2^n$  as follows.

\begin{lemma} \label{lemma1}
If $(\eta, g)$ satisfies \eqref{split0}, then $e_1^n$, $e_2^n$ defined in \eqref{e1} and \eqref{e2} satisfy
\begin{align*}
& \normM{e_1^n} \leq C \eps^{-4s} \sum_{i,j = 1}^2 \normM{|v|^2 (\partial^2_{v_i v_j}  + \partial_{x_i x_j}  )f_{in} } + \half \norm{(-\Delta_x)^{2s}\rho_{in}}_{L_x^2},
\\ &\normM{e_2^n \equilibrium} \leq \half \norm{(-\Delta_x)^{2s}\rho_{in}}_{L_x^2}.
\end{align*}
Consequently,  $\tilde \eta^{n+1}$ obtained in \eqref{0327} has the following estimate
\begin{align*}
  \norm{\tilde\eta^{n+1} \equilibrium}_{L^{\infty}_{\equilibrium^{-1}}} 
\leq 
 \Delta t \norm{(-\Delta_x)^{2s}\rho_{in}}_{L_x^\infty}.
\end{align*}
\end{lemma}

\begin{proof}
Rewrite $e_1^n$ and $e_2^n$ in their integral forms:
\begin{align*} 
 e_1^n 
=  \frac{1}{\Delta t^2}\int_{t^n}^{t^\np}  (\partial_s g(s, x,v) - \partial_s g(t^\np, x,v)) \rd s 
= \frac{1}{\Delta t^2}\int_{t^n}^{t^\np} \int_s^{t^\np} \partial_{\tau \tau } g(\tau, x, v) \rd \tau \rd s. 
\end{align*}
Similarly,
\begin{align} \label{1210}
\qquad e_2^n = \frac{1}{\Delta t^2} \int_{t^n}^{t^\np} \int_s^{t^\np} \partial_{\tau \tau} \eta (\tau, x, v) \rd \tau \rd s.
\end{align}
Recall that 
$\eta(t,x,v) = h(t,x+\eps v)$ where $h(t,y)$ solves 
\begin{equation*}
\partial_t h = - (-\Delta_y)^s h, \qquad h(0,y) = \rho_{in}(y).
\end{equation*}   
 Then $\partial_{tt}h(t,y)$ satisfies 
\begin{equation} \label{1211}
\partial_t (\partial_{tt}h) = - (-\Delta_y)^s (\partial_{tt}h), \qquad \partial_{tt}h(0,y) = (-\Delta_y)^{2s}\rho_{in}(y).
\end{equation} 
Note that $\partial_{tt} \eta(t, x, v) =  \partial_{tt} h(t, x+\eps v)$. Thus by the Minkowski inequality 
we have
\begin{align} \label{0619}
\normM{e_2^{n} \equilibrium} \leq \frac{1}{\Delta t^2} \int_{t^n}^{t^\np} \int_s^{t^\np} \normM{\partial_{\tau \tau} \eta (\tau, \cdot, \cdot) \mathcal M} \rd \tau \rd s
\leq 
  \norm{(-\Delta_x)^{2s} \rho_{in}}_{L_x^2} .
\end{align}

To bound $e_1^n$, we make use of equation \eqref{eqn:111} for $f$. Note that $\partial_{tt}f$ solves the same equation as $f$ with initial condition 
\begin{align*}
\partial_{tt}f(0,x,v) & = \eps^{2-4s}(v\cdot \nabla_x)(v\cdot \nabla_x) f_{in}(x,v) - \eps^{1-4s} v \cdot \nabla_x (\Lop^s f_{in})
\\ 
& \quad \,
   - \eps^{1-4s} \Lop^s(v \cdot \nabla_x f_{in}) + \eps^{-4s} \Lop^s \Lop^s f_{in}.
\end{align*}
Due to the non-negativity of the operator $\mathcal{L}^s$ in $L^2({\mathcal{M}^{-1} \rd v})$ (see \cite{cesbron2012anomalous} for a proof) , we have
\begin{align} \label{0620}
\normM{\partial_{tt}f} \leq \normM{\partial_{tt}f_{in}} \leq C \eps^{-4s} \sum_{i,j = 1}^2 \normM{|v|^2 (\partial^2_{v_i v_j}  + \partial_{x_i x_j}  )f_{in} }.
\end{align}
Then 
\begin{align*}
\normM{e_1^n} & = \frac{1}{\Delta t^2}  \normM{\int_{t^n}^{t^\np} \int_s^{t^\np} \partial_{\tau \tau} g(\tau, x, v) \rd \tau \rd s}
\\& \leq  \frac{1}{\Delta t^2} \int_{t^n}^{t^\np} \int_s^{t^\np}  \normM{\partial_{\tau \tau} (f - \eta \equilibrium)(\tau, \cdot, \cdot)} \rd \tau \rd s
\\& \leq  \frac{1}{\Delta t^2} \int_{t^n}^{t^\np} \int_s^{t^\np}  \normM{\partial_{\tau \tau} f} + \normM{\partial_{\tau \tau } \eta \equilibrium} \rd \tau \rd s.
\end{align*}
The desired bound for $e_1^n$ follows from \eqref{0619} and \eqref{0620}.

For $\tilde \eta^{n}$, note from \eqref{err_eta} that 
\[
\norm{\tilde \eta^{n+1}\equilibrium}_{L_{\equilibrium^{-1}}^\infty} =  \norm{\tilde \eta^{n+1}}_{L_{x,v}^\infty} =  \norm{\tilde \eta^{n+1}}_{L_x^\infty} \leq \norm{\tilde \eta^{n}}_{L_x^\infty} +  \Delta t^2\norm{e_2^n}_{L_x^\infty},
\]
thanks to the Maximum principle that the operator $(I + \Delta t (-\Delta_x)^s)^{-1}$ satisfies (see Lemma~\ref{lem:well-posedness-scheme} a)). Then summing over $n$ on both sides, we have
\begin{align*}
\norm{\tilde \eta^{n+1}}_{L_x^\infty} \leq \sum_n \Delta t^2 \norm{e_2^n}_{L_x^\infty}\,.
\end{align*}
By \eqref{1210} and \eqref{1211}, we have 
\begin{equation*}
  \norm{e_2^n}_{L_x^\infty} 
\leq 
  \norm{\partial_{tt} \eta(t, \cdot)}_{L_{x}^\infty}  
\leq 
  \norm{(-\Delta_x)^{2s}\rho_{in}}_{L_x^\infty}\,. \qedhere
\end{equation*}
\end{proof}

Next, we bound $g^{n+1}-\gstar$ and $g^\np$. These are for estimating the splitting and asymptotic errors in the fractional diffusion regime. In both cases, we start with the estimate for $I(\eta^\np, \equilibrium)$.

\begin{lemma}[Estimate of $I(\eta, \equilibrium)$] \label{lemma2}
Denote $\eta (x,v) = h(x+\eps v)$. Then for  $s\in (0, 1)$, we have
\begin{itemize}
\item[1)] $\normM{D^kI(\eta, \equilibrium)} \leq C \|h\|_{H_x^{k+1}} \eps^{s}$,  $0\leq k \leq 2$;
\item[2)] $|I(\eta,\equilibrium)| \leq C \norm{h}_{W^{1,\infty}}\equilibrium$.
\end{itemize}
Here $\normM{\cdot}$ is defined in \eqref{normM}, and $D^1 f = \partial_{v_i} f$, $D^2 f = \partial_{v_i v_j} f$.  
\end{lemma}

\begin{proof}
1) Recall the definition of $I$ in \eqref{I} and apply the Minkowski inequality for $L^2_x$. Then we have
\begin{align*}
\normM{I(\eta, \equilibrium)}
&\leq C \norm{\int_{\RR^d} \frac{\norm{(h(x+\eps v) - h(x+\eps w))}_{L^2_x} |\equilibrium(v) - \equilibrium(w)|}{|v-w|^{d+2s} \sqrt{\equilibrium(v)}} \rd w}_{L^2_v} 
\\ &= C \norm{\int_{| v - w| >\frac{1}{\eps}} \cdot  \rd w}_{L^2_v} + C \norm{ \int_{|v - w| <\frac{1}{\eps}}  \cdot  \rd w}_{L^2_v}
 =: C   \norm{I_1}_{L^2_v} 
     +  C \norm{I_2}_{L^2_v}.
\end{align*}
For $I_1$, we  have
\begin{align*}
   \norm{(h(x+\eps v) - h(x+\eps w))}_{L^2_x}
\leq
  \norm{h(x+\eps v)}_{L^2_x}
  + \norm{h(x+\eps w)}_{L^2_x}
\leq
  2 \norm{h}_{L^2_x}.
\end{align*}
Thus, 
\begin{align*}
 \|I_1\|_{L_v^2}& \leq C \norm{h}_{L_x^2}
 \norm{\int_{|v - w| > \frac{1}{\eps}} \frac{ |\equilibrium(v) - \equilibrium(w)| }{| v - w |^{d+2s} } \frac{1}{\sqrt{\equilibrium(v)}}\rd w}_{L^2_v} .
 \end{align*}
 If $|v-w|>\frac{1}{\eps}$ and $|w| > \half |v|$, then $\equilibrium(w) \leq C \equilibrium(v)$. In this case we have
\begin{align} 
& \quad \,
\norm{\int_{|v - w| > \frac{1}{\eps}, |w| > \half |v|} \frac{ |\equilibrium(v) - \equilibrium(w)| }{| v - w |^{d+2s} } \frac{1}{\sqrt{\equilibrium(v)}}\rd w}_{L^2_v} \nonumber
\\ 
&\leq 
  C \norm{\int_{|v - w| > \frac{1}{\eps}, |w| > \half |v|}  \frac{\sqrt{\equilibrium(v)}}{|v-w|^{d+2s}} \rd w}_{L_v^2}  \nonumber 
\\ 
&\leq 
 C \int_{|w| > \frac{1}{\eps}}  \frac{1}{|w|^{d+2s}} \rd w \norm{\equilibrium}_{L_v^1}^\half = C \eps^{2s}\,, \label{1226}
\end{align}
where the last inequality is a consequence of Young's convolution inequality.
If $|v-w|>\frac{1}{\eps}$ and $|w| \leq \half |v|$, then $\equilibrium(v) \leq C \equilibrium(w)$, $|v| > \frac{2}{3\eps}$ and $|v-w| \geq \half |v|$, and we have
\begin{align}
& \quad \, 
\norm{\int_{|v - w| > \frac{1}{\eps}, |w| \leq \half |v|} \frac{ |\equilibrium(v) - \equilibrium(w)| }{| v - w |^{d+2s} } \frac{1}{\sqrt{\equilibrium(v)}}\rd w}_{L^2_v}\nonumber 
\\ 
&\leq 
  C \norm{\int_{|v - w| > \frac{1}{\eps}, |w| \leq \half |v|}  \frac{\equilibrium(w)}{|v-w|^{d+2s}} \frac{1}{\sqrt{\equilibrium(v)}} \rd w}_{L_v^2}  \nonumber 
\\ 
&\leq 
  C \norm{\int_{\RR^d} \frac{\equilibrium(w)}{|v|^{d+2s}} \frac{1}{\sqrt{\equilibrium(v)}} \mathbf{1}_{|v|> \frac{2}{3\eps}}  \rd w  }_{L_v^2} \nonumber
\\ 
&\leq 
  C (\int_{|v|> \frac{2}{3\eps}} \equilibrium(v) \rd v)^\half \int \equilibrium(w) \rd w = C \eps^s .   \label{1227}
\end{align}
Combining \eqref{1226} and \eqref{1227} gives 
\begin{equation*}
\norm{I_1}_{L_v^2} \leq C \norm{h}_{L_x^2} \eps^{s}.
\end{equation*}

For $I_2$, we rewrite it as
\begin{align*}
  \normv{I_2} 
& = \norm{\int_{|v-w| \leq 1} \frac{\norm{(h(x+\eps v) - h(x+\eps w))}_{L^2_x} |\equilibrium(v) - \equilibrium(w)|}{|v-w|^{d+2s} \sqrt{\equilibrium(v)}} \rd w}_{L^2_v} 
\\ 
& \quad \,  
  + \norm{\int_{1\leq |v-w| \leq \frac{1}{\eps}} \frac{\norm{(h(x+\eps v) - h(x+\eps w))}_{L^2_x} |\equilibrium(v) - \equilibrium(w)|}{|v-w|^{d+2s} \sqrt{\equilibrium(v)}} \rd w}_{L^2_v} 
\\ & =: \normv{I_{21}} + \normv{I_{22}}.
\end{align*}
Note that 
\begin{align*}
& \norm{h(x+\eps v) - h(x+\eps w)}_{L_x^2} \leq \eps \norm{\nabla h}_{L_x^2} |v-w|,
\\ & \equilibrium(v) - \equilibrium(w) = \nabla \equilibrium(\xi)\cdot (v-w), \quad  \text{where} \quad \xi = (1-t) v+ tw, ~t \in (0,1).
\end{align*}
Then for $|v-w| \leq 1$, we have 
\begin{equation} \label{0812}
  |\nabla \equilibrium(\xi)| 
\leq 
  C |\nabla \equilibrium(v)| \leq C \equilibrium (v),
\end{equation}
where for the first inequality we have used the fact that $|\xi- v| \leq 1$ and the second comes as a result of Proposition~\ref{M-decay}.
Therefore $I_{21}$ satisfies
\begin{align} \label{I21-new}
\normv{I_{21}} \leq \eps \norm{\nabla h}_{L_x^2} 
\norm{\int_{|v-w|<1} \frac{1}{|v-w|^{d+2s-2}} \sqrt{\equilibrium(v)} \rd w}_{L_v^2} \leq C \eps \norm{\nabla h}_{L_x^2}.
\end{align}
Similarly, for $I_{22}$, we have 
\begin{align} \label{I22-8}
  \normv{I_{22}} 
&\leq 
  C \eps  \norm{\nabla h}_{L_x^2} \norm{ \int_{1\leq |v-w| \leq \frac{1}{\eps}, |w| \leq \half |v|} \frac{1}{|v-w|^{d+2s-1}} \frac{|\equilibrium(v)-\equilibrium(w)|}{\sqrt{\equilibrium(v)}} \rd w}_{L_v^2} \nonumber
\\ 
& \quad \, 
  +  C \eps \norm{\nabla h}_{L_x^2} \norm{ \int_{1\leq |v-w| \leq \frac{1}{\eps}, |w| > \half |v|} \frac{1}{|v-w|^{d+2s-1}} \frac{|\equilibrium(v)-\equilibrium(w)|}{\sqrt{\equilibrium(v)}} \rd w}_{L_v^2} \nonumber
\\ 
&=: C \eps  \norm{\nabla h}_{L_x^2} (I_{22,1} + I_{22,2}).
\end{align}
For $I_{22,1}$,  using $|v-w|>\half |v|$, we get
\begin{align}
I_{22,1} & \leq C \int \equilibrium(w) \rd w \left[\int_{\frac{2}{3} < |v|< \frac{2}{\eps}}   \left(\frac{1}{|v|^{d+2s-1}} \frac{1}{\sqrt{\equilibrium(v)}} \right)^2 \rd v  \right]^\half \leq C \eps^{s-1}.  \label{I22-new}
\end{align}
For $I_{22,2}$, we have $\equilibrium(w) \leq C \equilibrium(v)$. Thus, by Young's convolution inequality, we have
\begin{align} \label{I222}
I_{22,2} \leq C\norm{\equilibrium}_{L_v^1} \int_{1\leq|w|\leq \frac{1}{\eps}} \frac{1}{|w|^{d+2s-1}} \rd w \leq C \eps^{2s-1}.
\end{align}
Applying \eqref{I22-new} and \eqref{I222} to \eqref{I22-8}, we get
\begin{equation} \label{I22-new}
\normv{I_{22}} \leq  C   \norm{\nabla h}_{L_x^2} \eps^{s}.
\end{equation}
The combination of \eqref{I21-new} and \eqref{I22-new} leads to 
\begin{equation*}
\normM{I(\eta, \equilibrium)} \leq C   \norm{\nabla h}_{L_x^2} \eps^{s}.
\end{equation*}

The first-order derivatives of $\equilibrium$ are
\begin{align*}
\partial_{v_i} I(\eta, \equilibrium) &= \partial_{v_i} \int_{\RR^d} \frac{( h(x+\eps v) -h(x+\eps w)) (\equilibrium(v) - \equilibrium(w))}{|v-w|^{1+2s}} \rd w
\\ & = \partial_{v_i} \int_{\RR^d} \frac{( h(x+\eps v) -h(x+\eps v- \eps w)) (\equilibrium(v) - \equilibrium(v-w))}{|w|^{1+2s}} \rd w
\\ & = \eps I(\partial_{x_i} h, \equilibrium) + I (\eta, \partial_{v_i} \equilibrium).
\end{align*}
These can be estimated similarly as $I(\eta, \equilibrium)$. In particular, we have $\eps I(\partial_x h, \equilibrium) \leq  C \norm{h}_{H_x^2} \eps^{s+1}$ directly following the first estimate in part 2). For  $I (h, \partial_{v_i} \equilibrium)$, the same estimate holds since 
\begin{align*}
|\partial_v^k \equilibrium| \leq C_k \equilibrium,
\end{align*}
thanks to Proposition~\ref{M-decay}.
Likewise, 
\begin{align*}
\partial_{v_iv_j} I(\eta, \equilibrium) &= \partial_{v_j} (\eps I(\partial_{x_i} \eta, \equilibrium) + I(\eta, \partial_{v_i} \equilibrium))
\\ &= \eps^2 I(\partial_{x_i x_j} \eta, \equilibrium) +  \eps I(\partial_{x_j}\eta, \partial_{v_i} \equilibrium) 
+ \eps I(\partial_{x_i}\eta, \partial_{v_j} \equilibrium) + I(h, \partial_{v_iv_j} \equilibrium),
\end{align*}
and the rest of the estimate proceeds exactly the same as before.

3) We divide the integral in $I(h, \equilibrium)$ into several subdomains and estimate them individually:
\begin{align*}
|I(\eta, \equilibrium) |& \leq \int \frac{|h(x+\eps v) - h(x+ \eps w)|| \equilibrium(v) - \equilibrium(w)|}{|v-w|^{d+2s}} \rd w
 \\ & = \int_{|v-w|>1} \cdot ~ \rd w  + \int_{|v-w|\leq 1} \cdot ~ \rd w =: I_3 + I_4.
\end{align*}

For $I_3$, we separate the two cases where $|w|\geq\half |v|$ and $|w|<\half |v|$. When $|w|\geq\half |v|$, we have $\equilibrium(w) \leq C\equilibrium(v)$. Thus,
\begin{align*}
& \quad \,
\int_{|v-w|>1,|w|\geq\half|v|} \frac{|h(x+\eps v) - h(x+ \eps w)||\equilibrium(v) - \equilibrium(w)|}{|v-w|^{d+2s}} \rd w 
\\ 
&\leq 
  C\norm{h}_{\infty} \equilibrium(v) \int_{|v-w|>1} \frac{1}{|v-w|^{d+2s}} \rd w = C  \norm{h}_{\infty}  \equilibrium(v) .
\end{align*}
When $|w|<\half |v|$, we have $|v| > \frac{2}{3}$, $|v-w|>\half |v|$ and $\equilibrium(v) \leq C\equilibrium(w)$. Therefore,
\begin{align*}
& \quad \, 
\int_{|v-w|>1,|w|<\half|v|} \frac{|h(x+\eps v) - h(x+ \eps w)||\equilibrium(v) - \equilibrium(w)|}{|v-w|^{d+2s}} \rd w 
\\ 
&\leq 
  C\norm{h}_{\infty} \int_{\mathbb{R}^d} {\bf 1}_{|v|>\frac{2}{3}}\frac{\equilibrium(w)}{|v|^{d+2s}} \rd w  \leq  C \norm{h}_{\infty}  \left( \frac{1}{|v|^{d+2s}} {\bf 1}_{|v|>\frac{2}{3}}\right) \int \equilibrium(w) \rd w
\\ &  \leq C \norm{h}_{\infty}  \equilibrium(v).
\end{align*}

For $I_4$,  we have  
\begin{align*}
I_4 \leq C \eps \norm{\nabla h}_\infty \int_{|v-w|<1} \frac{|\nabla \equilibrium(\xi)|}{|v-w|^{d+2s-2}} \rd w, \quad  \xi = (1-t)v + tw, ~ t\in (0,1),
\end{align*}
where $|\nabla \equilibrium(\xi)|$ has the same estimate as in \eqref{0812}. Therefore,
\begin{align*}
I_4 \leq C \eps \norm{\nabla h}_\infty \equilibrium(v).
\end{align*}
Combining the estimates above gives the pointwise bound of $I$.
\end{proof}

To proceed with the error estimate, we will consider two separate regimes sketched in Fig.~\ref{fig: regimes}.

\subsection{Regime I:  kinetic regime with $\eps^{2s} \geq \Delta t^{2s\beta}$}
First we show that when $\eps$ is large compared to $\Delta t$ and $\gamma $ is chosen to satisfy \eqref{cond:gamma1} and \eqref{cond:gamma2}, the accuracy is controlled by some positive power of $\Delta t$. In the following estimates, we do not keep track of the error dependence on $\gamma$, but it is expected that the constant increases with larger values of $\gamma$. So in practice, we always choose $\gamma$ to be a constant slightly larger than 2 that satisfies the aforementioned conditions.

\begin{lemma}[Estimate of $\norm{g^n}$ in Regime I] \label{lemma3}
Let $(\eta^n, g^n)$ be the solution to \eqref{scheme22}--\eqref{scheme21} and $h^n(x+\eps v) = \eta^n(x,v)$. 
Suppose $\eps^{2s} \geq {\Delta t^{2s\beta}}$ with $\beta < \frac{1}{4s}$. Then for any $\gamma$ such that 
\begin{equation} \label{cond:gamma1}
0< \gamma \leq \sqrt{\frac{\lambda_0 - 1}{\lambda_0} }\alpha ,  \qquad \text{for some fixed } \quad \lambda_0 >1,
\end{equation}
where $\alpha = \frac{\eps^{2s}}{\Delta t}$, 
we have the following estimates for $g^n$ and $g^\nh$:
\begin{itemize}
\item[1)] $\normM{g^n} \leq C \eps^{-s} \left( \eps^s \normM{g_{in}} + \norm{\rho_{in}}_{H_x^1} \right)  $, ~
$\norm{g^n}_{L_{x,v}^2} \leq C \eps^{-s} \left(\eps^s \normM{g_{in}} + \norm{\rho_{in}}_{H_x^1} \right)  $;  
\item[2)] $\normM{g^\nh} \leq  C\eps^{-s} \left( \eps^s\normM{g_{in}} + \norm{\rho_{in}}_{H_x^1} \right)  $, ~
$\norm{g^\nh}_{L_{x,v}^2} \leq  C \eps^{-s} \left( \eps^s \normM{g_{in}} + \norm{\rho_{in}}_{H_x^1} \right) $;  
\item[3)] $\normMinf{g^n} := \norm{\frac{g^n}{\equilibrium}}_\infty  \leq  C \eps^{-2s} ( \eps^s \normMinf{g_{in}}  +  \norm{\rho_{in}}_{W_x^{1,\infty}})$.
\end{itemize}
\end{lemma}

\begin{proof}
1)  Without loss of generality, we assume that the final time is $t=1$ such that 
\begin{align*}
   N\Delta t = 1, \qquad 0\leq n \leq N. 
\end{align*}
By \eqref{semis1} and \eqref{semis2}, we have
\begin{equation} \label{g*}
\gstar = (\alpha + \gamma - \Lop^s)^{-1} (\alpha g^n - I^{n+1}), 
\qquad
g^\np = (\alpha - \gamma + \eps v \cdot  \nabla_x )^{-1} \alpha \gstar\,,
\end{equation}
where $I^{n+1} =  I(\eta^{n+1}, \mathcal{M})$.  Then we have 
\begin{align} \label{eqn0807}
\frac{g^\np}{\sqrt{\equilibrium}} &= (\alpha - \gamma + \eps v \cdot \nabla_x)^{-1} \frac{\alpha}{\sqrt{\equilibrium}} (\alpha + \gamma- \Lop^s)^{-1} \alpha \sqrt{\equilibrium} \left( \frac{g^n}{\sqrt{\equilibrium}} - \frac{1}{\alpha}\frac{I^\np}{\sqrt{\equilibrium}} \right)  \nonumber
\\&= \Hop \left( \frac{g^n}{\sqrt{\equilibrium}}  \right) - \Hop \left( \frac{1}{\alpha}\frac{I^\np}{\sqrt{\equilibrium}} \right)  \nonumber
\\&= \Hop^{n+1} \left( \frac{g_{in}}{\sqrt{\equilibrium}}  \right) - \sum_{p=1}^{n+1} \Hop^p \left( \frac{1}{\alpha} \frac{I^{n+1-p}}{\sqrt{\equilibrium}}\right) ,
\end{align}
where
\begin{align*} 
\Hop = (\alpha - \gamma + \eps v \partial_x)^{-1} \frac{\alpha}{\sqrt{\equilibrium}} (\alpha + \gamma- \Lop^s)^{-1} \alpha \sqrt{\equilibrium} 
 =: \Hop_1 \Hop_2,   
\end{align*}
with 
\begin{align} \label{H12}
\Hop_1 = (\alpha - \gamma + \eps v \cdot \nabla_x)^{-1} \alpha, \qquad 
\Hop_2 = \frac{1}{\sqrt{\equilibrium}} (\alpha + \gamma- \Lop^s)^{-1} \alpha \sqrt{\equilibrium} .
\end{align}
It is easy to check that $\Hop_1$ is bounded operator from $L_{x,v}^2 $ to itself with bound 
\begin{equation*} 
\norm{\Hop_1}_{L_{x,v}^2 \rightarrow {L_{x,v}^2}} \leq \frac{\alpha}{\alpha - \gamma}.
\end{equation*}
For $\Hop_2$, we conduct the energy estimate. Suppose $F$ and $G$ satisfy $\Hop_2 F = G$, or equivalently,
\begin{align*}
F = \frac{\alpha  + \gamma}{\alpha} G - \frac{1}{\alpha} \frac{1}{\sqrt{\equilibrium}} \Lop^s (\sqrt{\equilibrium} G).
\end{align*}
Multiply the equation by $G$ and integrate with respect to $x$ and $v$.  Using the fact that
\begin{align*}
\int_{\RR^d} \int_{\RR^d} \frac{G}{\sqrt{\equilibrium}} \Lop^s (\sqrt{\equilibrium} G) \rd x \rd v 
= \int_{\RR^d} \int_{\RR^d} \frac{1}{\equilibrium}  \sqrt{\equilibrium} G \Lop^s (\sqrt{\equilibrium} G) \rd x \rd v  \leq 0,
\end{align*}
one has 
\begin{align} \label{H2estimate}
\norm{G}_{L_{x,v}^2} \leq \frac{\alpha}{\alpha + \gamma} \norm{F}_{L_{x,v}^2}
\Longrightarrow
\norm{\Hop_2}_{L_{x,v}^2 \rightarrow {L_{x,v}^2}} \leq \frac{\alpha}{\alpha + \gamma}.
\end{align}
Therefore $\norm{\Hop}_{L_{x,v}^2 \rightarrow {L_{x,v}^2}} \leq \frac{\alpha^2}{\alpha^2 - \gamma^2}$. Plugging it into \eqref{eqn0807}, one has
\begin{align} \label{1114}
\normM{g^\np}  
\leq \left(  \frac{\alpha^2}{\alpha^2- \gamma^2} \right)^{n+1} \normM{g_{in}} +  \frac{1}{\alpha}\sum_{p=1}^{n+1} \left( \frac{\alpha^2}{\alpha^2-\gamma^2}\right)^p\normM{I^{n+2-p}}.
\end{align}
Note that since $\eps^{2s} \geq \Delta t^{2s \beta}$ with $\beta  < \frac{1}{4s}$ , we have
\begin{align*}
\alpha^2- \gamma^2 \geq \frac{1}{\Delta t^{2-4s \beta}}-\gamma^2 \geq N- \gamma^2 \geq n- \gamma^2 .
 \end{align*}
Therefore,
 \begin{equation} \label{eqn:0806}
 \left(  \frac{\alpha^2}{\alpha^2- \gamma^2} \right)^{n+1}  \leq \left(1+\frac{\gamma^2}{n-\gamma^2} \right)^{n+1} \leq Ce^{\gamma^2}, \quad \textrm{for all}~ n \geq \gamma^2+1.
 \end{equation}
 In the case when $n < \gamma^2+1$,  as long as $\gamma$ is chosen such that $\gamma^2 \leq \frac{\lambda_0 - 1}{\lambda_0} \alpha^2$ for some fixed $\lambda_0 >1$, we have
 \begin{equation}\label{eqn:0807}
  \left(  \frac{\alpha^2}{\alpha^2- \gamma^2} \right)^{n+1}  \leq  \lambda_0^{n+1}.
 \end{equation}
 In sum, \eqref{1114} continues as 
 \begin{align*}
 \normM{g^\np} \leq C\normM{g_{in}} + \frac{C}{\alpha} \sum_{p=1}^{n+1} \normM{I^p}, \quad C = \max \left\{ e^\gamma, \lambda_0^{(\gamma^2+1)} \right\}.
 \end{align*}
 Using Lemma~\ref{lemma2} part 1) and Lemma~\ref{lem:well-posedness-scheme}, one has
 \begin{align*}
 \normM{g^\np} \leq C \eps^{-s} \left( \eps^s \normM{g_{in}} + \norm{\rho_{in}}_{H_x^1}\right) .
 \end{align*}
 The bound for $\norm{g^\np}_{L_{x,v}^2}$ is a direct consequence of the boundedness of $\normM{g^\np}$. 

2) By~\eqref{g*} we have
\begin{align*}
\frac{g^\nh}{\sqrt{\equilibrium}} 
&= \frac{1}{\sqrt{\equilibrium}} (\alpha + \gamma -\Lop^s)^{-1} \sqrt{\equilibrium} (\alpha \frac{g^n}{\sqrt{\equilibrium} } -  \frac{ I^\np}{\sqrt{\equilibrium}})
\\ & = \Hop_2 \left(  \frac{g^n}{\sqrt{\equilibrium}} - \frac{1}{\alpha}  \frac{ I^\np}{\sqrt{\equilibrium}} \right)
\\ & = \Hop_2^{n+1} \left( \frac{g_{in}}{\sqrt{\equilibrium}}  \right) - \sum_{p=1}^{n+1} \Hop_2^p \left( \frac{1}{\alpha} \frac{I^{n+1-p}}{\sqrt{\equilibrium}}\right). 
\end{align*}
Using the bound \eqref{H2estimate}, we have for any $\gamma >0$, 
\begin{align*}
 \normM{g^\nh} 
\leq 
  \normM{g_{in}} + \frac{1}{\alpha} \sum_{p=1}^{n+1} \normM{I^p} .
\end{align*}
The estimate for $\normM{g^\nh}$ again follows from Lemma~\ref{lemma2} part 1) and Lemma~\ref{lem:well-posedness-scheme}.

3)  Denote
\begin{align*}  
\Hop_0 =  (\alpha - \gamma + \eps v \cdot \nabla_x )^{-1} \alpha (\alpha + \gamma - \Lop^s)^{-1}  \alpha, 
\qquad 
 \Hop_3 = (\alpha + \gamma - \Lop^s)^{-1}  \alpha.
\end{align*}
Then by using \eqref{g*}, 
one can write $g^\np$ as 
\begin{align*}
g^\np = \Hop_0 g^n - \frac{1}{\alpha} \Hop_0 I.
\end{align*}
Then our goal is to find the bounds of $\Hop_3$ and $\Hop_1$ defined in \eqref{H12} with respect to the norm $L^\infty_{\equilibrium^{-1}}$. 
To this end, we will use the maximum principle argument. Let $q_1$ satisfies 
\begin{align*}
q_1 = \Hop_3 g = (\alpha + \gamma - \Lop^s)^{-1} \alpha g.
\end{align*}
Then $g = \left( \frac{\alpha + \gamma}{\alpha} - \frac{1}{\alpha} \Lop^s \right) q_1 $. Since $g \leq \norm{\frac{g}{\equilibrium}}_\infty \equilibrium$,
we have
\begin{align*}
g -  \norm{\frac{g}{\equilibrium}}_\infty  \equilibrium = \left( \frac{\alpha + \gamma}{\alpha} - \frac{1}{\alpha} \Lop^s \right) \left( q_1 - \frac{\alpha}{\alpha + \gamma}  \norm{\frac{g}{\equilibrium}}_\infty  \equilibrium  \right) \leq 0.
\end{align*}
Multiply the above inequality by $\sgn \left( q_1 - \frac{\alpha}{\alpha + \gamma}  \norm{\frac{g}{\equilibrium}}_\infty \equilibrium  \right)$, where $\sgn(x) =1$ if $x>0$ and $0$ otherwise. Then it becomes 
\begin{align}
& \frac{\alpha+\gamma}{\alpha}  \left( q_1 - \frac{\alpha}{\alpha + \gamma}  \norm{\frac{g}{\equilibrium}}_\infty  \equilibrium  \right)^+ - \frac{1}{\alpha} \nabla_v  \cdot \left[v  \left( q_1 - \frac{\alpha}{\alpha + \gamma}  \norm{\frac{g}{\equilibrium}}_\infty  \equilibrium  \right)^+ \right]  \nonumber 
\\ & \qquad + \frac{1}{\alpha}  \sgn \left( q_1 - \frac{\alpha}{\alpha + \gamma} \norm{\frac{g}{\equilibrium}}_\infty \equilibrium  \right)^+ (-\Delta_v)^s  \left( q_1 - \frac{\alpha}{\alpha + \gamma}  \norm{\frac{g}{\equilibrium}}_\infty  \equilibrium  \right) \leq 0.  \label{1215}
\end{align}
Note that for any $b(v)$, we have 
\begin{align*}
(\sgn b(v) ) (-\Delta_v)^s b(v) &= C_{s,d} (\sgn b)(v) \int \frac{b(v)-b(w)}{|v-w|^{1+2s}} \rd w 
\geq C_{s,d} (\sgn b)(v) \int \frac{b(v)-b^+(w)}{|v-w|^{1+2s}} \rd w
\\ 
& \hspace{-1cm} = C_{s,d}  \int \frac{b^+(v)-(\sgn b)(v) b^+(w)}{|v-w|^{1+2s}} \rd w 
\geq C_{s,d}  \int \frac{b^+(v)-b^+(w)}{|v-w|^{1+2s}} \rd w  = (-\Delta_v)^s b^+(v).
\end{align*}
Therefore, \eqref{1215} leads to
\begin{align*}
\frac{\alpha+\gamma}{\alpha}  \left( q_1 - \frac{\alpha}{\alpha + \gamma} \norm{\frac{g}{\equilibrium}}_\infty  \equilibrium  \right)^+ - \Lop^s \left( q_1 - \frac{\alpha}{\alpha + \gamma}  \norm{\frac{g}{\equilibrium}}_\infty  \equilibrium  \right)^+ \leq 0.
\end{align*}
Multiply the above equation by $\frac{1}{\equilibrium}\left( q_1 - \frac{\alpha}{\alpha + \gamma}  \norm{\frac{g}{\equilibrium}}_\infty  \equilibrium  \right)^+$ and integrate in $v$. We get 
\begin{align*}
\norm{\vpran{q_1 - \frac{\alpha}{\alpha + \gamma}  \norm{\frac{g}{\equilibrium}}_\infty  \equilibrium}^+}_{\equilibrium^{-1}} \leq 0 ,
\end{align*}
which implies $q_1 \leq \frac{\alpha}{\alpha + \gamma} \norm{\frac{g}{\equilibrium}}_\infty   \equilibrium$. Similarly, we can show that $-q_1 \leq  \frac{\alpha}{\alpha + \gamma}  \norm{\frac{g}{\equilibrium}}_\infty  \equilibrium$. Therefore,
\begin{align*} 
\norm{\frac{q_1}{\equilibrium}}_\infty \leq \frac{\alpha}{\alpha + \gamma} \norm{\frac{g}{\equilibrium}}_\infty \qquad \textrm{or} \qquad 
\normMinf{\Hop_3} \leq \frac{\alpha}{\alpha + \gamma}.
\end{align*}
Likewise, to bound $\Hop_1$, let $q_2 = \Hop_1 g = (\alpha -\gamma + \eps v \cdot  \nabla_x)^{-1} \alpha g$. Then
\begin{align*}
\left( \frac{\alpha -\gamma}{\alpha} + \frac{\eps}{\alpha} v \cdot  \nabla_x \right) q_2 = g.
\end{align*}
Since $g- \norm{\frac{g}{\equilibrium}}_\infty \equilibrium \leq 0$, we have 
\begin{align*}
\left( \frac{\alpha -\gamma}{\alpha} + \frac{\eps}{\alpha} v \cdot  \nabla \right) \left( q_2 - \frac{\alpha}{\alpha-\gamma}  \norm{\frac{g}{\equilibrium}}_\infty  \equilibrium \right) = g-  \norm{\frac{g}{\equilibrium}}_\infty  \equilibrium \leq 0.
\end{align*}
Multiplying the above inequality by $\left( q_2 - \frac{\alpha}{\alpha-\gamma}  \norm{\frac{g}{\equilibrium}}_\infty  \equilibrium \right)^+$ and conducting an energy estimate leads to 
$
q_2 \leq \frac{\alpha}{\alpha - \gamma}  \norm{\frac{g}{\equilibrium}}_\infty  \equilibrium.
$
Consequently, we have
\begin{equation*} 
\normMinf{\Hop_1} \leq \frac{\alpha}{\alpha-\gamma}
\qquad \text{and} \qquad
\normMinf{\Hop_0} \leq \frac{\alpha^2}{\alpha^2-\gamma^2}.
\end{equation*}
Then following the same proof as for part 1), we have 
\begin{align*} 
\normMinf{g^\np}  
& \leq \left(  \frac{\alpha^2}{\alpha^2- \gamma^2} \right)^{n+1} \normMinf{g_{in}} +   \frac{1}{\alpha}\sum_{p=1}^{n+1} \left( \frac{\alpha^2}{\alpha^2-\gamma^2}\right)^p\normMinf{I^{n+1-p}}.
\end{align*}
This leads to 
\begin{align*}
\normMinf{g^\np}  \leq  C\normMinf{g_{in}} + \frac{C}{\alpha} \sum_{p=1}^{n+1} \normMinf{I^p}, 
 \leq  C\normMinf{g_{in}}  + \frac{C}{\eps^{2s}}  \norm{\rho_{in}}_{W^{1,\infty}},  
\end{align*}
where $C = \max\{ e^\gamma, \lambda_0^{[\gamma^2+1]}\}$ and
 we have used part 2) of Lemma~\ref{lemma2} and Lemma~\ref{lem:well-posedness-scheme}. 
\end{proof}

\vspace{-6pt}
Next we bound the derivatives of $g^\np$ and $g^\nh$.  
\vspace{-1pt}
\begin{lemma}[Estimate of $\normM{\partial_{v_i} g^\np}$ and $\normM{\partial_{v_i} g^\nh}$ in Regime I]  \label{lemma7}
Suppose $\eps^{2s} \geq \Delta t^{2s\beta}$ and 
\begin{align} \label{cond:gamma2}
2 < \gamma \leq  \min\left\{\sqrt{\frac{\lambda_2-4}{\lambda_2}} \frac{\eps^{2s}}{\Delta t} = \sqrt{\frac{\lambda_2-4}{\lambda_2}} \alpha\,,~ 4 \right\} \quad \text{for some fixed} ~ \lambda_2 >4.
\end{align}
Then the following bounds hold:
\begin{itemize}
\item[1)] for any $i= 1,2, \cdots, d$,
\begin{align*}
 \normM{\partial_{v_i} g^\np} \leq C (\eps^{-s} + \eps^{1-3s}) L_1^0, 
\quad
 \normM{\partial_{v_i} g^\nh} \leq C (\eps^{-s} + \eps^{1-3s}) L_1^0,
\end{align*}
where
\begin{equation} \label{eqn:L01}
L_1^0  =  \eps^s\normM{\partial_{x_i}g_{in}} + \eps^s \normM{\partial_{v_i} g_{in}}  + \norm{\rho_{in}}_{H_x^2};
\end{equation}

\item[2)] for any $i, j= 1,2, \cdots, d$,
\begin{align*}
\normM{\partial_{v_iv_j} g^\np}
  + \normM{\partial_{v_iv_j} g^\nh} 
\leq 
  C(\eps^{-s} + \eps^{1-3s} + \eps^{2-5s})L_2^0,
\end{align*}
where
\begin{align}\label{eqn:L02}
 L_2^0 
&= \eps^s \normM{\partial_{x_i}g_{in}} 
     + \eps^s \normM{\partial_{x_j}g_{in}} 
     + \eps^s \normM{\partial_{v_i}g_{in}}  \nn
\\
& \quad \,
     + \eps^s \normM{\partial_{v_i v_j} g_{in}}
     + \norm{\rho_{in}}_{H_x^2}.
\end{align}
\end{itemize}
\end{lemma}
\Ni The proof of Lemma~\ref{lemma7} is technical and resembles the proof of Lemma~\ref{lemma3}. We have left it to the appendix.

Now we estimate the difference $g^\np - g^\nh$. 
\begin{lemma}[Estimate of $g^\np - g^\nh$ in Regime I] \label{lemma6} 
For $\eps^{2s} \geq \Delta t^{2s\beta}$ with $\beta < \frac{1}{4s}$, we have the following estimates for $g^\np - g^\nh$: 
\begin{itemize}
\item[1)] $\norm{\average{v}^{m} (g^\np - g^\nh) }_{L_{x,v}^2}  \leq C\left( \frac{\Delta t}{\eps^{2s}}\right)^{\min\{ s+\halfd - m , 1\}} \eps^{-s} L_0^0$,
\\
where $m< s+ \halfd$ and 
$L_0^0 = \eps^s \normM{g_{in}} + \norm{\rho_{in}}_{H_x^2}$;

\item[2)] $\norm{\average{v}^{m} \partial_{v_i} (g^\np - g^\nh) }_{L_{x,v}^2}  \leq C   \left( \frac{\Delta t}{\eps^{2s}}\right)^{\min\{ s+\halfd - m , 1\}}(\eps^{-s} + \eps^{1-3s}) L_1^0$, \quad  
\\
where $m< s+ \halfd$, $i = 1, 2, \cdots, d$, and $L_1^0$ is defined in \eqref{eqn:L01};
\item[3)] $\norm{\average{v}^{m} \partial_{v_i v_j} (g^\np - g^\nh) }_{L_{x,v}^2}  \leq C   \left( \frac{\Delta t}{\eps^{2s}}\right)^{\min\{ s+\halfd - m , 1\}} (\eps^{-s} + \eps^{2-5s})L_2^0$, 
\\
where $m< s+ \halfd$, $i = 1, 2, \cdots, d$ and $L_2^0$ is defined in \eqref{eqn:L02}; 
\item[4)] $\norm{\average{v}^{p}(-\Delta_v)^s (g^\np-g^\nh)}_{L_{x,v}^2} \leq C  \left( \frac{\Delta t}{\eps^{2s}}\right)^{\min\{s+\frac{d}{2} -p, 1\}}  (\eps^{-s} + \eps^{2-5s})(L_0^0 + L_2^0)$,
\\
where $-2s <p< s + \frac{d}{2}$.  
\end{itemize}
\end{lemma}

\begin{proof}
1) From \eqref{semis2}, one sees that
\begin{equation} \label{1003}
g^\np - g^\nh = \frac{1}{\alpha} (\gamma g^\np - \eps v \cdot  \nabla_x g^\np ), 
\qquad \alpha = \frac{\eps^{2s}}{\Delta t}.
\end{equation}
Then for a fixed $m$, 
\begin{align*}
& \quad \, 
\norm{\average{v}^m (g^\np - g^\nh)}_{L_{x,v}^2}^2 = \int\!\! \int  \average{v}^{2m} | g^\np - g^\nh|^{2\theta} | g^\np - g^\nh|^{2(1-\theta)} \rd x \rd v
\\& = \frac{C}{\alpha^{2(1-\theta)}} \int  \!\! \int \average{v}^{2m} | g^\np - g^\nh|^{2\theta} |\gamma  g^\np - \eps v \cdot  \nabla_x g^\np |^{2(1-\theta)} \rd x \rd v 
\\& \leq \frac{C}{\alpha^{2(1-\theta)}} \int \!\! \int \average{v}^{2m} (| g^\np|^{2\theta} + |g^\nh|^{2\theta})
(|\gamma  g^\np|^{2(1-\theta)} + \eps^{2(1-\theta)}|v \cdot \nabla_x g^\np|^{2(1-\theta)} ) \rd x \rd v .
\end{align*}
Expand the product in the above inequality. The most difficult term is 
\begin{align*}
\int \!\! \int  \average{v}^{2m} |g^\nh|^{2\theta} \eps^{2(1-\theta)}|v \cdot \nabla_x g^\np|^{2(1-\theta)} \rd x  \rd v,
\end{align*}
which can be estimated as follows:
\begin{align*}
& \quad \, 
\int\!\! \int  \average{v}^{2m} |g^\nh|^{2\theta} \eps^{2(1-\theta)}|v \cdot \nabla_x g^\np|^{2(1-\theta)}  \rd x \rd v 
\\ & = \eps^{2(1-\theta)} \int \!\! \int  \average{v}^{2m+2(1-\theta)}  | g^\nh|^{2\theta} |\nabla_x g^\np|^{2(1-\theta)}  \rd x \rd v 
\\ & \leq  \eps^{2(1-\theta)}  \left( \int \!\! \int \average{v}^{2m+2(1-\theta)}  | g^\nh|^2 \rd x  \rd v \right)^\theta   \left( \int \!\! \int \average{v}^{2m+2(1-\theta)}|\nabla_x g^\np|^2 \rd x \rd v \right)^{1-\theta}.
\end{align*}
To use Lemma~\ref{lemma3} to bound the right hand side of the inequality above, we need 
\begin{align*}
2m+2(1-\theta) \leq d+2s.
\end{align*}
Since $\theta \in (0,1)$, this requires $m < \halfd + s$. In this case, we can choose
\begin{align*}
\theta = \max \left\{ 1-({d}/{2} + s -m), 0 \right\}.
\end{align*}
and apply Lemma~\ref{lemma3}. The other terms are bounded similarly.  Note that while Lemma~\ref{lemma3} is for the estimates of $\| g^{n+1}\|_{M^{-1}}$ rather than $\|\nabla_x g^{n+1}\|_{M^{-1}}$, $(\nabla_x g^{n+\frac{1}{2}}, \nabla_x g^{n+1})$ satisfy similar equations as $(g^{n+\frac{1}{2}}, g^{n+1})$ in \eqref{semis1}-\eqref{semis2}, with $\eta$ replaced by $\nabla_x \eta$. Consequently, a similar estimate for $\nabla_x g^{n+1}$ applies with the bound now depends on $\|\rho_{in}\|_{H_x^2}$ instead of $\|\rho_{in}\|_{H_x^1}$.

2) Equation \eqref{1003} implies that
\begin{equation} \label{dgg}
\partial_{v_i}(g^\np - g^\nh) = \frac{1}{\alpha}\left\{  \gamma \partial_{v_i} g^\np - \eps [ v\cdot \nabla_x (\partial_{v_i} g^\np) + \partial_{x_i} g^\np]  \right\}.
\end{equation}
Then for a fixed $m$, 
\begin{align*}
& \quad \, 
\norm{\average{v}^m \partial_{v_i}(g^\np - g^\nh)}_{L_{x,v}^2}^2 = \int \!\! \int \average{v}^{2m} |\partial_{v_i} (g^\np - g^\nh)|^{2\theta} |\partial_{v_i} (g^\np - g^\nh)|^{2(1-\theta)} \rd x \rd v 
\\
&= \frac{C}{\alpha^{2(1-\theta)}} \int \!\! \int  \average{v}^{2m} |\partial_{v_i} (g^\np - g^\nh)|^{2\theta} |\gamma \partial_{v_i} g^\np - \eps [v\cdot \nabla_x (\partial_{v_i} g^\np) + \partial_{x_i} g^\np ]|^{2(1-\theta)} \rd x \rd v 
\\
&\leq 
  \frac{C}{\alpha^{2(1-\theta)}} \int \!\! \int  \average{v}^{2m} [|\partial_{v_i} g^\np|^{2\theta} + |\partial_{v_i} g^\nh|^{2\theta}]   \times
\\ 
& \hspace{3cm} 
  \times [|\gamma \partial_{v_i} g^\np|^{2(1-\theta)}+ \eps^{2(1-\theta)} |v \cdot \nabla_x (\partial_{v_i} g^\np)|^{2(1-\theta)} + |\partial_{x_i} g^\np|^{2(1-\theta)}] \rd x \rd v .
\end{align*}
The rest follows similarly as for part 1) by replacing Lemma~\ref{lemma3} with Lemma~\ref{lemma7} in the proof.

3) This part is similar to Part 2). From \eqref{dgg}, we have 
\begin{equation*} 
\partial_{v_i v_j}(g^\np - g^\nh) = \frac{1}{\alpha}\left\{  \gamma \partial_{v_i} g^\np - \eps [ v\cdot \nabla_x (\partial_{v_i v_j} g^\np) + \partial_{x_i}\partial_{v_j} g^\np + \partial_{x_j}\partial_{v_i} g^\np]  \right\}.
\end{equation*}
Then for a fixed $m$, 
\begin{align*}
& \quad \, \norm{\average{v}^m \partial_{v_i v_j}(g^\np - g^\nh)}_{L_{x,v}^2}^2 
\\
&= \int \!\! \int  \average{v}^{2m} |\partial_{v_iv_j} (g^\np - g^\nh)|^{2\theta} |\partial_{v_i v_j} (g^\np - g^\nh)|^{2(1-\theta)}  \rd x \rd v 
\\
& = \frac{C}{\alpha^{2(1-\theta)}} \int \!\! \int  \average{v}^{2m} |\partial_{v_i v_j} (g^\np - g^\nh)|^{2\theta} |\gamma \partial_{v_i v_j} g^\np 
\\ 
& \hspace{3cm}
  - \eps [v\cdot \nabla_x (\partial_{v_i v_j} g^\np) + \partial_{x_i} \partial_{v_j}g^\np + \partial_{x_j} \partial_{v_i} g^\np ]|^{2(1-\theta)} \rd x  \rd v 
\\
& \leq \frac{C}{\alpha^{2(1-\theta)}} \int \!\! \int  \average{v}^{2m} [|\partial_{v_i v_j} g^\np|^{2\theta} + |\partial_{v_i v_j} g^\nh|^{2\theta}]  
[|\gamma \partial_{v_i v_j} g^\np|^{2(1-\theta)} + |\partial_{x_i} \partial_{v_j} g^\np|^{2(1-\theta)} 
\\ 
& \hspace{3cm}  
  + \eps^{2(1-\theta)} |v \cdot \nabla_x (\partial_{v_i v_j} g^\np)|^{2(1-\theta)} + |\partial_{x_j} \partial_{v_i} g^\np|^{2(1-\theta)}]  \rd x \rd v .
\end{align*}
The rest of the estimate follows similarly as that in part 2).

4) For simplicity, we denote $\bar g = g^\np-g^\nh$ and rewrite  
\begin{align*}
\average{v}^{p}(-\Delta_v)^s \bar g = (-\Delta_v)^s (\average{v}^{p} \bar g) - [\average{v}^{p}, (-\Delta_v)^s] \bar g,
\end{align*}
where $[\cdot, \cdot]$ is the commutator. For the first term, by interpolation we have for any $p < s + \halfd$, 
\begin{align}
  \norm{(-\Delta_v)^s (\average{v}^{p} \bar g) }_{L_{x,v}^2} ^2  \nonumber
&\leq C\int\!\! \int  |\average{v}^{p} \bar g|^2 \rd x \rd v  +  \int \!\! \int  |\nabla_v^2 (\average{v}^{p} \bar g)|^2 \rd x \rd v 
\nonumber
\\ 
&\leq 
  C (L_2^0)^2 \left( \frac{\Delta t}{\eps^{2s}}\right)^{\min\{ 2s+d-2p, 2\}}  (\eps^{-s} + \eps^{1-3s}+ \eps^{2-5s})^2.\label{est1}
\end{align}
For the second term,  from Lemma~\ref{lemt1}, we have, for $-2s< p< \frac{d}{2} + 2s$,
\begin{align*}
\norm{ [\average{v}^{p}, (-\Delta_v)^s] \bar g}_{L_{x,v}^2} \leq C \normxv{\average{v}^{p} \bar g} + 
\normxv{\average{v}^{p-1} |\nabla_v \bar g|}.
\end{align*}
Substitute part 1) and 2) into the right-hand side of the inequality. Then it becomes 
\begin{align*}
 \norm{ [\average{v}^{p}, (-\Delta_v)^s] \bar g  }_{L_{x,v}^2} 
& \leq C  \left( \frac{\Delta t}{\eps^{2s}}\right)^{\min\{s+\frac{d}{2} -p, 1\}} \eps^{-s} L_0^0 
\\ & \qquad  + C \left( \frac{\Delta t}{\eps^{2s}}\right)^{\min\{s+\frac{d}{2}+1-p , 1\}} 
 (\eps^{-s} +\eps^{1-3s})L_1^0, \quad \quad  -2s <p< s + \frac{d}{2},
\end{align*}
which together with \eqref{est1} implies
\begin{gather*}
\norm{\average{v}^{p}(-\Delta_v)^s \bar g}_{L_{x,v}^2} 
\leq  C  \left( \frac{\Delta t}{\eps^{2s}}\right)^{\min\{s+\frac{d}{2} -p, 1\}}  (\eps^{-s} + \eps^{1-3s}+ \eps^{2-5s})(L_0^0 + L_2^0),  
\end{gather*}
with $-2s<p < s + \frac{d}{2}$.
\end{proof}

\smallskip
We are now ready to state our main theorem in Regime I. The main idea is to construct a new weighted norm that can compensate the slow decay of the equilibrium at the tail.

\begin{theorem}[Error estimate in Regime I] \label{thm:2}
Assume \eqref{cond1} and suppose 
\begin{align*}
   \eps^{2s} \geq \Delta t^{ 2s\beta }, 
\qquad
  \beta \leq \beta_0 <  \frac{1}{6s}, 
\end{align*} 
where $\beta_0$ is the biggest constant that satisfies \eqref{beta1}--\eqref{beta5}.
Then the solution $f^n = \eta^n \equilibrium + g^n$ obtained from solving \eqref{semi_scheme} with $\gamma$ chosen to satisfy \eqref{cond:gamma1} and \eqref{cond:gamma2} has the following estimate
\begin{align*}
\normM{\average{v}^{-b} \ftilde^n} \leq C \Delta t^\zeta  
\end{align*}
for $1<b\leq \min\{d+2s, \frac{d}{2}+3s^-\}$ and  $\zeta(\beta_0, b) >0$ explicitly computable. 
Here $\tilde f^n$ is defined in \eqref{def:alpha-tilde-f}, $n \Delta t = T$, and $C$ is a constant depending on $T$ and initial data.
\end{theorem}
\begin{proof}
Multiply \eqref{erroreqn2} by the weight $(1+\delta \average{v})^{-b}$, where $b$ and $\delta$ are two constants that will be made precise later. We have
\begin{align*}
& \quad \, 
\alpha \left[ (1+\delta \average{v})^{-b} \ftilde^\np - (1+\delta \average{v})^{-b} \ftilde^n \right]  + \eps v \partial_x ((1+\delta \average{v})^{-b} \ftilde^\np)   
\\ 
& = \Lop^s ((1+\delta \average{v})^{-b} \ftilde^\np) + [(1+\delta \average{v})^{-b}, \Lop^s] \ftilde^\np 
\\ & \hspace{2cm} 
+ (1+\delta \average{v})^{-b} (\Lop^s- \gamma \id)(g^\np - g^\nh)
+ (1+\delta \average{v})^{-b} e^n, 
\end{align*}
where
\begin{equation} \label{en}
e^n =\Delta t( e_1^n + \eps^{2s} e_2^n \equilibrium ).
\end{equation}
Conduct an energy estimate for the equation by multiplying it by $(1+\delta \average{v})^{-b} \frac{\ftilde^\np}{\equilibrium} $ and integrating with respect to $x$ and $v$. We obtain 
\begin{align} \label{energy1}
& \quad \, \frac{\alpha}{2} \left( \normM{(1+\delta \average{v})^{-b} \ftilde^\np}^2 - \normM{(1+\delta \average{v})^{-b} \ftilde^n}^2\right) \nonumber
\\  
&\leq \int \!\! \int (1+\delta \average{v})^{-b}\ftilde^\np [(1+\delta \average{v})^{-b}, \Lop^s] \ftilde^\np  \frac{1}{\equilibrium} \rd v \rd x \nonumber
\\ 
& \quad \,
+ \int \!\! \int (1+\delta \average{v})^{-b}  \Lop^s(g^\np - g^\nh) (1+\delta \average{v})^{-b}  \ftilde^\np  \frac{1}{\equilibrium} \rd v\rd x \nonumber
\\ & \qquad 
+ \int\!\! \int (1+\delta \average{v})^{-2b} e^n \frac{\tilde f^\np}{\equilibrium} \rd v \rd x,
\end{align}
where we have used the fact that $\int \frac{f \Lop^s f}{\equilibrium} \rd v \leq 0$ . 
By the definition of $\Lop^s$ in \eqref{eqn:LFP},
\begin{align*}
& \quad \, 
\frac{1}{\sqrt{\equilibrium}}(1+\delta \average{v})^{-b} \Lop^s (g^\np - g^\nh) 
\\ 
& = \frac{1}{\sqrt{\equilibrium}} (1+\delta \average{v})^{-b} \left(  \nabla_v \cdot (v (g^\np - g^\nh) )- (-\Delta_v)^{s} (g^\np - g^\nh) \right).
\end{align*}
To proceed, first note that the lower bound in Part 4) of Lemma~\ref{lemma6} can be removed. This is because for $p < -2s$, we have
\begin{align*}
 \norm{\average{v}^{p}(-\Delta_v)^s (g^\np-g^\nh)}_{L_{x,v}^2} 
&\leq 
 C \norm{\average{v}^{-2s^-}(-\Delta_v)^s (g^\np-g^\nh)}_{L_{x,v}^2} 
\\ 
&\leq 
  C  \left( \frac{\Delta t}{\eps^{2s}}\right)^{\min\{\frac{d}{2}+3s^-, 1\}}  (\eps^{-s} + \eps^{2-5s})(L_0^0 + L_2^0).
\end{align*}
where $s^-$ denotes any number smaller than $s$. Therefore, if $b$ satisfies 
\begin{align} \label{b-cond2}
   -b+\half(d+2s)+1<s+\halfd, 
\quad \text{or equivalently,}\quad
   b > 1,  
\end{align}
then by Lemma~\ref{lemma6},
we have the estimate
\begin{align*} 
& \quad \,\normM{(1+\delta \average{v})^{-b} \Lop^s (g^\np - g^\nh) }  \nonumber
\\ 
& \leq 
  C \delta^{-b} \left\{  \left( \frac{\Delta t}{\eps^{2s}}\right)^{\min\{b-1,1\}} (\eps^{-s} + \eps^{1-2s}) L_1^0 
+ \left( \frac{\Delta t}{\eps^{2s}} \right)^{\min\{\frac{d}{2}+3s^-,1\}} (\eps^{-s} + \eps^{2-5s})(L_0^0 + L_2^0)  \right\} \,.
\end{align*}
If we further assume that 
\begin{align} \label{b-cond3}
b< 3s^- + \frac{d}{2} + 1\,,
\end{align}
the above inequality becomes 
\begin{align}\label{energy2}
\normM{(1+\delta \average{v})^{-b} \Lop^s (g^\np - g^\nh) }
\leq C \delta^{-b}  \left(  \frac{\Delta t}{\eps^{2s}}\right)^{\min\{b-1,1\}} (\eps^{-s} + \eps^{2-5s}) L^0,
\end{align}
where $ L^0 = L^0_0+L_2^0$. This is because in Region I we have $\Delta t/\eps^{2s} \leq 1$.

Next, we control the commutator term $ [(1+\delta \average{v})^{-b}, \Lop^s] \ftilde^\np$. By Lemma~\ref{lemma:commutator}, if
\begin{equation} \label{finfM}
\normMinf{\ftilde^\np} \leq C \eps^{-2s},
\end{equation} 
then we have
\begin{align} \label{energy3}
& \left|\int \!\! \int (1+\delta \average{v})^{-b}  \ftilde^\np [(1+\delta \average{v})^{-b}, \Lop^s]  \ftilde^\np \frac{1}{\equilibrium} \rd x \rd v \right|  \nonumber
\\ & \qquad  \leq  C\delta^{\frac{2s}{d+4}}  \normM{(1+\delta \average{v})^{-b}  \ftilde^\np}^2 + C  {\color{black} \eps^{-4s}} \delta^{\min\{1, 2s\}} .
\end{align}
Here \eqref{finfM} is fulfilled by combining part 2) of Lemma~\ref{lemma3}, Lemma~\ref{lemma1},  the bound $\normMinf{g(t^n)}  \leq C \eps^{-2s}$ and the fact that $\tilde f^n = \tilde\eta^n \equilibrium  + \tilde g^n$. Indeed, from the argument of maximum principle that similar to the proof of part 2) in Lemma~\ref{lemma3}, the boundedness of $\normMinf{g(t^n)} $ is guaranteed. 
Denote 
\begin{align*}
   F^{n}= \normM{(1+\delta \average{v})^{-b} \ftilde^n}.
\end{align*} 
Plug \eqref{energy2} and \eqref{energy3} into \eqref{energy1}.
Then we get
\begin{align} \label{Fneqn}
(F^\np)^2 - (F^n)^2  &\leq a \Delta t(F^\np)^2 + r,
\end{align}
where
\begin{align*} 
& a = 1 + C \delta^{\frac{2s}{d+4}} \eps^{-2s} ,
\\ &r = C \frac{\Delta t }{\eps^{6s}} \delta^{\min\{1, 2s\}} 
+ C \delta^{-2b} \frac{\Delta t}{\eps^{4s}}  \left( \frac{\Delta t}{\eps^{2s}}\right)^{\min\{ 2, 2b-2\}} \hspace{-1.5cm}(\eps^{-s} + \eps^{2-5s})^2(L^0)^2  
+ \frac{\Delta t}{\eps^{4s}} \normM{e^n }^2.
\end{align*}
Note that here we have used the fact that for $b>0$, $\normM{(1+\delta \average{v})^{-b}e^n } \leq \normM{e^n}$.
Optimize in $r$ by choosing 
\begin{equation*} 
\delta = \left \{ \begin{array}{cc} \left( \frac{\Delta t}{\eps^{2s}}\right)^{\frac{\min\{ 2, 2b-2\}}{1+2b}} \eps^{\frac{4-8s}{2b+1}}, \quad & s>\half,
\\ \left( \frac{\Delta t}{\eps^{2s}}\right)^{\frac{\min\{ 1, b-1\}}{b+s}} , \quad & s< \half.
\end{array} \right.
\end{equation*}
This correspondingly gives 
\begin{align*}
a = \left\{ \begin{array}{cc} 
1+ C  \left( \frac{\Delta t}{\eps^{2s}}\right)^{\frac{\min\{ 2, 2b-2\}}{1+2b}  \frac{2s}{d+4} }  \eps^{\frac{4-8s}{2b+1}\frac{2s}{d+4}-2s}, & s>\half
\\ 1 + C\left( \frac{\Delta t}{\eps^{2s}}\right)^{\frac{\min\{ 1, b-1\}}{b+s} {\frac{2s}{d+4}}}  \eps^{-2s}  , & s\leq \half
\end{array} \right. 
\end{align*}
and
\begin{align*}
\frac{r}{\Delta t} = \left\{ \begin{array}{cc}
C \left[  \left( \frac{\Delta t}{\eps^{2s}}\right)^{\frac{\min\{ 2, 2b-2\}}{1+2b}} \eps^{\frac{4-8s}{2b+1}-6s}
+ \left( \frac{\Delta t}{\eps^{2s}}\right)^{\frac{\min\{ 2, 2b-2\}}{1+2b}} \eps^{4-14s-2b\frac{4-8s}{2b+1}}(L^0)^2  
+ \frac{1}{\eps^{4s}} \normM{e^n}^2 \right]  & s>\half,
\\  C \left[  \left( \frac{\Delta t}{\eps^{2s}}\right)^{2s\frac{\min\{ 1, b-1\}}{b+s}} \eps^{-6s}
+  \left( \frac{\Delta t}{\eps^{2s}}\right)^{\frac{2s\min\{1,b-1\}}{b+s}}\eps^{-6s}(L^0)^2  
+ \frac{1}{\eps^{4s}} \normM{e^n}^2 \right] & s\leq \half .
\end{array} \right. 
\end{align*}
Recall that in this regime we have the relation $\eps^{2s} > \Delta t^{2s\beta}$ with $\beta < 1/4s$. In order for the power of $\frac{\Delta t}{\eps^{2s}}$ in both $a$ and $\frac{r}{\Delta t}$  to be nonnegative, $\beta$ must satisfy the following inequalities:
\begin{align}
s> \half :\qquad & (1-2s\beta) \frac{\min\{ 2, 2b-2\}}{1+2b} \frac{2s}{d+4}-2s\beta \left( \frac{8s-4}{(2b+1)(d+4)} +1\right) \geq 0, \label{beta1}
\\ & (1-2s\beta) \frac{\min\{ 2, 2b-2\}}{1+2b} -2\beta\left(  \frac{4s-2}{2b+1} +3s\right) \geq 0,
\\ & (1-2s\beta) \frac{\min\{ 2, 2b-2\}}{1+2b}-2\beta \left( 2-7s+b \frac{4s-2}{2b+1} \right) \geq 0\,,
\\ & 1- 6s\beta >0 \,; \label{0114}
\\ s\leq \half: \qquad  & (1-2s\beta) \frac{\min\{ 1, b-1\}}{b+s} \frac{2s}{d+4}-2s\beta \geq 0,
\\ & (1-2s\beta)\frac{\min\{ 1, b-1\}}{b+s} -3\beta \geq 0,
\\ & 1- 6s\beta >0 \,.  \label{beta5}
\end{align}
Here \eqref{0114} and \eqref{beta5} arise from the definition of $e^n$ in \eqref{en}, specifically $e^n = \Delta t(e_1^n + \eps^{2s} e_2^n \equilibrium)$ and the estimate of $e_1^n$, $e_2^n$ in Lemma~\ref{lemma1}. This results in $\eps^{-4s} \|e^n\|_{\equilibrium^{-1}}^2 \lesssim \Delta t \eps^{-12s}$. 
Since $\beta \in (0,1)$, it is clear that there exists  $\beta_0 \in (0,1)$ such that if $\beta \in (0, \beta_0)$ then the inequalities above are satisfied. 
Then \eqref{Fneqn} leads to 
\begin{align} \label{energy4}
(F^\np)^2 \leq e^{aT} (F^0)^2 + C \frac{e^{aT}}{a} \frac{r}{\Delta t}.
\end{align}
Moreover, since $F^0 =0$ and
\begin{align*}
\normM{(1+\delta \average{v})^{-b} \ftilde^n} \geq \normM{((1+\delta) \average{v})^{-b} \ftilde^n} \geq 2^{-b} \normM{\average{v}^{-b} \ftilde^n},
\end{align*}
the bound \eqref{energy4} reduces to 
\begin{align*} 
\normM{\average{v}^{-b} \ftilde^n} \leq C e^{\frac{aT}{2}} \sqrt{ \frac{r}{\Delta t}}  ,
\end{align*}
where $C$ is a constant only depends on final time $T$. 
\end{proof}

\begin{remark}
We can choose $b=1+2s$, then one choice of $\beta_0$ is 
\begin{align*}
\beta_0 = \left\{ \begin{array}{cc} 
\frac{1}{2s+(d+4)(12s-2)}, & s>\half,
\\  \frac{s}{(1+3s)(d+4)+2s^2}, & s\leq \half. 
\end{array} \right. 
\end{align*}
In this case we have $\normM{\average{v}^{-b} \ftilde^n} \sim \Delta t^\zeta$, where 
\begin{align*}
\zeta = \left\{ \begin{array}{cc} 
\frac{4s}{3+4s} - 2\beta \frac{13s-2+16s^2}{3+4s},   & s>\half,
\\  \frac{2s}{1+3s} - \beta \frac{3+9s+4s^2}{1+3s}.  & s \leq \half,
\end{array} \right.
\qquad \beta \in (0, \beta_0).
\end{align*}
\end{remark}

\subsection{Regime II: diffusive and intermediate regime with $\eps^{2s} \leq \Delta t^{ 2s\beta}$}
We show that when $\eps$ is small compared with $\Delta t$, $g^n$ can be bounded in terms of $\eps$. In this regime, the choice of $\gamma$ will depend on the relative relationship between $\Delta t$ and $\eps$ specified as follows.
\begin{lemma}[Estimate of $g^n$ in Regime II] \label{lemma4}
For  $\eps^{2s} \leq \Delta t^{2s \beta}$, where 
\begin{align*}
   \beta \leq \beta_0 <  \frac{1}{6s} \,,
\end{align*} 
and $\beta_0$ is the biggest constant that satisfies \eqref{beta1}--\eqref{beta5}.
We choose $\gamma$ according to  
\begin{align} \label{gamma-II}
\gamma = \left \{ \begin{array}{cc} \sqrt{3} & \text{if}~ \eps^{2s} < \Delta t \\ \sqrt{3} \Delta t^{2s\beta -1} & \text{if}~ \Delta t \leq \eps^{2s} \leq \Delta t^{2s\beta}.
\end{array} \right.
\end{align}
Then 
\begin{align*}
\|g^n\|_2 \leq C \sqrt \Delta t( \norm{\rho_{in}}_{H_x^1} +  \normM{g_{in}} ), 
\end{align*}
for $n\geq n_0$ outside an initial layer. Here $C$ is a constant that does not depend on $\gamma$.

\end{lemma}
\begin{proof}
Combining the two equations in \eqref{g*} leads to 
\begin{align*}
g^\np = (\alpha - \gamma + \eps v \partial_x )^{-1} \alpha  (\alpha + \gamma - \Lop^s)^{-1} (\alpha g^n - I^\np ).
\end{align*}
Note that the derivation of \eqref{1114} is independent of the regime, and it remains valid within this regime. Therefore, we have,
\begin{align} \label{12162}
\normM{g^\np}  
\leq \left| \frac{\alpha^2}{\alpha^2- \gamma^2} \right|^{n+1} \normM{g_{in}} +  \frac{\alpha}{\alpha^2- \gamma^2} \left[1 + \cdots +  \left|  \frac{\alpha^2}{\alpha^2- \gamma^2} \right|^{n}  \right] \normM{I^{n+2-p}}.
\end{align}

When $\eps^{2s} \leq \Delta t$ and $\gamma = \sqrt{3}$, we have $\alpha = \frac{\eps^{2s}}{\Delta t} \leq 1$, $\left| \frac{\alpha^2}{ \alpha^2 - \gamma^2 }  \right| <\half$ and $\left| \frac{\alpha}{ \alpha^2 - \gamma^2 }  \right| <\half$. Furthermore, by Lemma~\ref{lemma2} Part 1) and Lemma~\ref{lem:well-posedness-scheme} we have $\normM{I} \leq C \eps^{s} \norm{\rho_{in}}_{H_x^1}$. Therefore,
\[
\normM{g^\np} \leq  \frac{1}{2^n} \normM{g_{in}}  + C \eps^{s} \left( 1 + \half + \cdots + \frac{1}{2^{n+1}}\right)\norm{\rho_{in}}_{H_x^1}  \leq  C \eps^{s} \norm{\rho_{in}}_{H_x^1} + \frac{1}{2^n} \normM{g_{in}} .
\]
Choose $n\geq n_0$ with $n_0 = \frac{\log \Delta t}{\log \half}$ such that $\frac{1}{2^n} \leq \Delta t$. Then the conclusion follows from the bound $\normxv{g^n} \leq C \normM{g^n}$ since $\CalM^{-1}$ has a positive lower bound. 

When $\Delta t < \eps^{2s} \leq \Delta t^{2s\beta}$, by letting $\gamma = \sqrt{3} \Delta t^{2s\beta -1} $, we again have $\left| \frac{\alpha^2}{ \alpha^2 - \gamma^2 }  \right| <\half$. Then \eqref{12162} leads to
\begin{align*}
 \normM{g^\np} & \leq  \frac{1}{2^n} \normM{g_{in}}  + C \eps^{s} \left|\frac{\alpha}{\alpha^2-\gamma^2}\right| \norm{\rho_{in}}_{H_x^1} 
 \\ & =     \frac{1}{2^n} \normM{g_{in}}  + C \frac{\eps^{s}}{ \alpha }\left|\frac{\alpha^2 }{\alpha^2-\gamma^2}\right| \norm{\rho_{in}}_{H_x^1} 
 =   \frac{1}{2^n} \normM{g_{in}}  + C \frac{\Delta t}{ \eps^{s} }\left|\frac{\alpha^2 }{\alpha^2-\gamma^2}\right| \norm{\rho_{in}}_{H_x^1} 
 \\ & \leq  \Delta t \normM{g_{in}}  + C \sqrt{\dt}  \norm{\rho_{in}}_{H_x^1}
\end{align*}
for sufficiently large $n$.
\end{proof}

\begin{remark}
The selection of $\gamma$ aims to ensure  $\left| \frac{\alpha^2}{ \alpha^2 - \gamma^2 }  \right| <\half$. It is sufficient to choose $\gamma > \sqrt{3} \alpha = \sqrt{3} \frac{\eps^{2s}}{ \Delta t}$. Although such choice works for both ranges in~\eqref{gamma-II}, in practice, a small $\gamma$ is preferred to minimize error. Therefore, when $\epsilon^{2s} < \Delta t$, we choose $\gamma = \sqrt{3}$. When $\Delta t \leq \epsilon^{2s}$, however, we have to use a large $\gamma$ by letting $\gamma = \sqrt{3} \Delta t^{2s\beta -1}$.
\end{remark}

\begin{theorem}[Error estimate in Regime II] \label{thm:3}
Suppose $\eps^{2s} < \Delta t^{2s\beta}$ and $t \geq t_0$ with $t_0$ being a fixed constant. Then the numerical error $\tilde f^n$ satisfies
\begin{align*}
\norm{\tilde f^n}_{L_{x,v}^2}^2  & \leq  
C\Delta t ^2  \norm{(-\Delta_x)^{2s}\rho_{in}}_{L_x^2}
+C_{s,t_0}  \eps^{2s} \norm{\rho_{in}}_{L^1}^2
+
C \left[ \eps \norm{v \cdot \nabla_x f_{in}}_{L^1_{x,v}}  + C_{t_0} \eps^{2s} \norm{  v f_{in}}_{L^1_{x,v}} \right]
\\ & \qquad \qquad  + C \eps^{\frac{d}{2}+2s}\norm{\rho_{in}}_{\infty}  + C \eps^{d+4s} \norm{(-\Delta_x)^s \rho_{in}}_{x}^2 + C  \sqrt{\Delta t}( \norm{\rho_{in}}_{H_x^1} +  \normM{g_{in}} ).
\end{align*}
\end{theorem}

\begin{proof}
Recall the definition of the numerical error in \eqref{def:alpha-tilde-f}. We have
\begin{align*}
\norm{\tilde f^n}_{L_{x,v}^2}^2 &= \norm{(\eta(t^n,x,v)-\eta^n(x,v))\equilibrium + (g(t^n,x,v) - g^n(x,v))}_{L_{x,v}^2}^2
\\ &\leq 2\Delta t^2 \norm{e_2^n \equilibrium}_{L_{x,v}^2}^2 + 2\norm{g(t^n,x,v)}_{L_{x,v}^2}^2 + 2\norm{g^n(x,v)}_{L_{x,v}^2}^2.
\end{align*}
Note that
\begin{align*}
\normxv{g(t^n,x,v)} &= \normxv{f(t^n,x,v) - \rho(t^n,x) \equilibrium(v) + \rho(t^n,x) \equilibrium(v) - \eta(t^n,x,v) \equilibrium(v) }
\\ & \leq \normxv{f(t^n,x,v) - \rho(t^n,x) \equilibrium(v)} + \normxv{(\rho(t^n,x) - \rho(t^n, x+\eps v)) \equilibrium(v)}.
\end{align*}
Bound of the first term follows from Theorem~\ref{APlimit}. For the second term, write 
\begin{align*}
& \quad \, 
  \normxv{(\rho(t^n,x) - \rho(t^n, x+\eps v)) \equilibrium(v)}^2
\\ 
&\leq 
  C \int\int_{|v|\leq \frac{1}{\eps}} \frac{|\rho(x+\eps v)- \rho(x)|^2}{\average{v}^{2(d+2s)}} \rd x \rd v
+  C \int\int_{|v|> \frac{1}{\eps}} \frac{|\rho(x+\eps v)- \rho(x)|^2}{\average{v}^{2(d+2s)}} \rd x \rd v
 =: I_5 + I_6.
\end{align*}
For $I_6$, we have 
\begin{align*}
I_6 \leq C \norm{\rho_{in}}_{\infty} \int\int_{|v|> \frac{1}{\eps}} \frac{1}{\average{v}^{2(d+2s)}} \rd x \rd v
\leq C \eps^{\frac{d}{2}+2s}\norm{\rho_{in}}_{\infty} .
\end{align*}
For $I_5$, use the change of variables
\begin{align*}
y = x+ \eps  v, \quad \text{then} ~ v= \frac{y-x}{\eps}, ~ \rd v = \eps^{-d} \rd y.
\end{align*}
Then 
\begin{align*}
I_5 \leq C \eps^{d + 4s} \int\int_{|x-y| 
\leq  1} \frac{|\rho(y)- \rho(x)|^2}{|y-x|^{d+4s}} \rd x \rd y 
\leq C \eps^{d+4s} \norm{(-\Delta_x)^s \rho}_{x}^2  
\leq C \eps^{d+4s} \norm{(-\Delta_x)^s \rho_{in}}_{x}^2   .
\end{align*} 
The rest follows from Theorem~\ref{APlimit},  Lemma~\ref{lemma1} and Lemma \ref{lemma4}.
\end{proof}

\section{Appendix}

\subsection{Appendix A}
This part shows the proof of Proposition~\ref{M-decay}. 
\begin{proof}[Proof of Proposition 1]
Let $\chi(k) \in C_c^\infty$ be a cutoff function with 
\begin{align*}
\chi(k) = \left\{ \begin{array}{cc} 1 & \text{~for~~} k \in B(0,1) ,
\\ 0 & \text{~for~~} |k|>2.
\end{array} \right.
\end{align*}
Recall that the equilibrium $\equilibrium(v)$ has the following Fourier transform 
\begin{equation*}
\hat \equilibrium(k) = C e^{-\frac{1}{2s}|k|^{2s}}, \quad s \in (0,1).
\end{equation*}
Thus its derivative in the Fourier space reads
\begin{align*}
  \widehat{ \nabla^m \equilibrium} (k)  
= C (ik)^m \chi(k) \, e^{-\frac{1}{2s}|k|^{2s}} 
   + C (ik)^m (1-\chi(k)) \, e^{-\frac{1}{2s}|k|^{2s}}, 
\quad 
   (ik)^m = (i k_1)^{\alpha_1} \cdots (ik_d)^{\alpha_d}.
\end{align*}
Taking the inverse Fourier transform, one arrives at
\begin{align*}
\nabla^m \equilibrium(v) &= \finv (\widehat{ \nabla^m \equilibrium} (k) ) 
\\ & =  \finv (C (ik)^m \chi(k) e^{-\frac{1}{2s}|k|^{2s}} ) + \finv (C (ik)^m (1-\chi(k)) e^{-\frac{1}{2s}|k|^{2s}} ) 
\\ & =: M_1(v) + M_2(v).
\end{align*}
Since $C (ik)^m (1-\chi(k)) e^{-\frac{1}{2s}|k|^{2s}} $ is a Schwartz class function,  $M_2(v)$ decays faster than any polynomial. Thus $|M_2(v)| \leq C_{d, m} \equilibrium(v)$. For $M_1$, let 
\begin{align*}
\hat g(k):= C (ik)^m \chi(k), \quad g(v) = \finv(\hat g). 
\end{align*}
Then 
\begin{align*}
M_1(v) & = \finv(\hat g(k) \hat \equilibrium(k)) = \int_{\RR^d} g(v-w) \equilibrium(w) \rd w
\\ & = \int_{|w|\geq \half |v|} g(v-w) \equilibrium(w) \rd w +  \int_{|w|<\half |v|} g(v-w) \equilibrium(w) \rd w.
\end{align*}
When $|w|\geq \half |v|$, we have $\equilibrium(w) \leq  C \equilibrium(v)$. Thus 
\begin{align*}
   \int_{|w|\geq \half |v|} g(v-w) \equilibrium(w) \rd w  
\leq 
  C \equilibrium(v) \int_{\RR^d} |g| \rd w  
\leq 
  C_{d, m} \equilibrium(v) .
\end{align*}
When $|w|<\half |v|$, we have $|v-w| > \half |v|$. Since $\hat g$ is a Schwartz class function, in this subdomain $g(v-w) \leq C_\ell \average{v}^{-l}$ for any $l$. Therefore, 
\begin{align*}
  \int_{|w|<\half |v|} g(v-w) \equilibrium(w) \rd w 
\leq 
   C_\ell  \average{v}^{-l} \int_{\RR^d} \equilibrium(v) \rd v, 
\qquad \forall ~l. 
\end{align*}
The final result follows from combining the estimates for $M_1$ and $M_2$. 
\end{proof}

\subsection{Appendix B}
This part is devoted to the proof of Lemma~\ref{lemma7}. 
\begin{proof}[Proof of Lemma~\ref{lemma7}.]
Taking the derivative in $v_i$ of \eqref{semis1} and \eqref{semis2}, we have
\begin{subequations}\label{eqn0929}
\begin{numcases}{}
\alpha (\partial_{v_i} g^\nh - \partial_{v_i} g^n) = \Lop^s (\partial_{v_i} g^\nh) - (\gamma-1)\partial_{v_i} g^\nh - \partial_{v_i} I(\eta^\np, \equilibrium), \label{eqn1001}
\\ \alpha (\partial_{v_i} g^\np - \partial_{v_i} g^\nh) + \eps v \cdot \nabla_x (\partial_{v_i} g^\np) = \gamma \partial_{v_i} g^\np - \eps \partial_{x_i} g^\np, \label{eqn1002}
\end{numcases}
\end{subequations}
which leads to
\begin{align*}
  \partial_{v_i} g^\np 
& = (\alpha-\gamma+ \eps v \cdot \nabla_x )^{-1} \alpha (\alpha- \Lop^s + \gamma-1)^{-1} (\alpha \partial_{v_i} g^n - \partial_{v_i} I^\np)
\\
& \quad \,
- (\alpha-\gamma+\eps v \cdot  \nabla_x)^{-1} \eps \partial_{x_i} g^\np.
\end{align*}
Then similar to the proof of part 1) of Lemma~\ref{lemma3}, we write 
\begin{align} \label{eqn0930}
   \frac{\partial_{v_i} g^\np}{\sqrt{\equilibrium}} 
&= \Hop_1 \Hop_4 \left( \frac{\partial_{v_i} g^n}{\sqrt{\equilibrium}}\right) - \Hop_1 \Hop_4\left(\frac{1}{\alpha} \frac{ \partial_{v_i} I^\np}{\sqrt{\equilibrium}}\right) - \Hop_1 \eps^{1-2s} \Delta t \partial_{x_i} g^\np  \nonumber
\\ 
&= (\Hop_1 \Hop_4)^{n+1} \left( \frac{\partial_{v_i} g_{in}}{\sqrt{\equilibrium}}\right) - \sum_{p=1}^\np (\Hop_1 \Hop_4)^p \left( \frac{1}{\alpha} \frac{\partial_{v_i}I^{n+2-p}}{\sqrt{\equilibrium}}\right)   \nn
\\
& \quad \,
- \sum_{p=1}^{n+1} (\Hop_1 \Hop_4)^{p-1} \Hop_1 \eps^{1-2s} \Delta t \partial_{x_i} g^{n+2-p},
\end{align}
where  $\Hop_1$ is defined in \eqref{H12} and 
\[
\Hop_4 = \frac{1}{\sqrt{\equilibrium}} (\alpha - \Lop^s + \gamma -1)^{-1} \alpha \sqrt{\equilibrium}.
\]
Similar to \eqref{H2estimate}, we see that 
\begin{equation}\label{H3}
\norm{\Hop_4}_{L_{x,v}^2 \rightarrow L_{x,v}^2} \leq \frac{\alpha}{\alpha + \gamma - 1}.
\end{equation}
Denote $c = \frac{\alpha^2}{(\alpha-\gamma)(\alpha+ \gamma-1)}$. 
Then via iterations, \eqref{eqn0930} can be estimated as
\begin{align*}
\normM{\partial_{v_i} g^\np} &\leq  c^{n+1} \normM{\partial_{v_i} g_{in}} +  \eps^{-2s} \Delta t \sum_{p=1}^{n+1} c^p \normM{\partial_{v_i} I^{n+2-p}} 
\\ & \quad \, 
+ \sum_{p=1}^{n+1} c^{p-1} \frac{\alpha}{\alpha-\gamma} \eps^{1-2s} \Delta t \norm{\partial_{x_i} g^{n+2-p}}_{L_{x,v}^2}.
\end{align*}
By Lemma~\ref{lemma3} part 1), Lemma~\ref{lemma2} part 1) and Lemma~\ref{lem:well-posedness-scheme}, we obtain
\begin{align*}
\normM{\partial_{v_i} g^\np} &\leq  c^{n+1} \normM{\partial_{v_i} g_{in}} +  \eps^{-s} \Delta t \norm{\rho_{in}}_{H_x^2} \sum_{p=1}^{n+1} c^p 
\\ 
& \quad \, 
+    C \eps^{1-3s}  \frac{\alpha}{\alpha-\gamma}  \Delta t \left(\eps^s \normM{\partial_{x_i}g_{in}} + \norm{\rho_{in}}_{H_x^2} \right)   \sum_{p=1}^{n+1} c^{p-1} .
\end{align*}
Following the same argument as in \eqref{eqn:0806} and \eqref{eqn:0807}, if we choose $\gamma$ such that 
\begin{align} \label{cond1}
0 < \gamma \leq \sqrt{\frac{\lambda_1-2}{\lambda_1}} \frac{\eps^{2s}}{\Delta t} = \sqrt{\frac{\lambda_1-2}{\lambda_1}} \alpha
\end{align}
for some fixed $\lambda_1 >2$, then the above inequality reduces to 
\begin{align*}
\normM{\partial_{v_i} g^\np} &\leq  C \normM{\partial_{v_i} g_{in}} +  C\eps^{-s} \norm{\rho_{in}}_{H_x^2} 
+  C \eps^{1-3s}    \left( \eps^s \normM{\partial_{x_i}g_{in}} + \norm{\rho_{in}}_{H_x^2} \right)   
\\ & \leq C \normM{\partial_{v_i} g_{in}} +  C(\eps^{-s} + \eps^{1-3s}) \left( \eps^s \normM{\partial_{x_i}g_{in}} + \norm{\rho_{in}}_{H_x^2} \right).  
\end{align*}
In addition, from \eqref{eqn1001} one has
\begin{align*}
\partial_{v_i} g^\nh = (\alpha - \Lop^s + \gamma - 1)^{-1} (\alpha \partial_{v_i} g^n - \partial_{v_i} I^\np).
\end{align*}
Thus, 
\begin{align*}
\frac{\partial_{v_i} g^\nh}{\sqrt{\equilibrium}} = \Hop_4  \left( \frac{\partial_{v_i} g^n}{\sqrt{\equilibrium}}  - \frac{1}{\alpha} \frac{1}{\sqrt{\equilibrium}} \partial_{v_i} I^\np \right),
\end{align*}
which immediately leads to 
\begin{align*}
\normM{\partial_{v_i} g^\nh} \leq C(\normM{\partial_{v_i} g^n} + \frac{1}{\alpha} \normM{\partial_{v_i} I^\np}),
\end{align*}
thanks to \eqref{H3} with $\gamma >1$. Combining the bound for $\partial_{v_i} g^n$ above with Lemma~\ref{lemma2} part 1) and Lemma~\ref{lem:well-posedness-scheme}, we further have 
\begin{align*}
\normM{\partial_{v_i} g^\nh} & \leq C (\eps^{-s} + \eps^{1-3s}) \left(\eps^s  \normM{\partial_{x_i}g_{in}} + \norm{\rho_{in}}_{H_x^2} \right)  +C \normM{\partial_{v_i} g_{in}}
+ C \Delta t \eps^{-s} \|\rho_{in}\|_{H_x^{2}} 
\\ & \leq C (\eps^{-s} + \eps^{1-3s}) \left(\eps^s  \normM{\partial_{x_i}g_{in}} + \norm{\rho_{in}}_{H_x^2} \right) + C \normM{\partial_{v_i} g_{in}}  .
\end{align*}
For the second-order derivatives, differentiating \eqref{eqn1001} and \eqref{eqn1002} in $v_j$, we have
\begin{subequations}
\begin{numcases}{}
\alpha (\partial_{v_i v_j} g^\nh \!-\! \partial_{v_i v_j} g^n) 
= \Lop^s (\partial_{v_i v_j} g^\nh) - (\gamma-2)\partial_{v_i v_j} g^\nh - \partial_{v_i v_j} I(\eta^\np, \equilibrium),  \label{g1229}
\\ 
\alpha (\partial_{v_i v_j} g^\np \!-\! \partial_{v_i v_j } g^\nh) 
\!+\! \eps v \!\cdot \!\nabla_x (\partial_{v_i v_j} g^\np) 
\!=\! \gamma \partial_{v_i v_j} g^\np \! - \! \eps \partial_{x_i v_j} 
g^\np  
\!-\! \eps \partial_{x_j v_i} g^\np \!, 
\end{numcases}
\end{subequations}
which implies 
\begin{align} \label{g1227}
\partial_{v_i v_j} g^\np &= (\alpha -\gamma + \eps v \cdot \nabla_x)^{-1} \alpha (\alpha-\Lop^s+\gamma-2)^{-1}(\alpha \partial_{v_iv_j} g^n - \partial_{v_i v_j} I)  \nonumber
\\ 
& \qquad \qquad - (\alpha - \gamma + \eps v \cdot \nabla_x)^{-1} \eps(\partial_{x_i} \partial_{v_j} g^\np + \partial_{x_j} \partial_{v_i} g^\np).
\end{align}
Recall $\Hop_1$  defined in \eqref{H12} and denote 
\begin{align*}
  \Hop_5 
=  \frac{1}{\sqrt{\equilibrium}} (\alpha - \Lop^s + \gamma -2)^{-1} \alpha \sqrt{\equilibrium}.
\end{align*} 
Then \eqref{g1227} can be written as 
\begin{align} \label{g1228}
\frac{\partial_{v_i v_j} g^\np}{\sqrt{\equilibrium}} &= \Hop_1 \Hop_5 \left( \frac{\partial_{v_i v_j} g^n}{\sqrt{\equilibrium}} \right)
- \Hop_1 \Hop_5 \left( \frac{1}{\alpha} \frac{\partial_{v_i v_j} I}{\sqrt{\equilibrium}} \right)
-\Hop_1 \eps^{1-2s} \Delta t (\partial_{x_i} \partial_{v_j} g^\np + \partial_{x_j} \partial_{v_i} g^\np)  \nonumber
\\& = (\Hop_1 \Hop_5)^{n+1} \left( \frac{\partial_{v_i v_j}g_{in}}{\sqrt{\equilibrium}}\right)  \nonumber
- \sum_{p=1}^\np (\Hop_1 \Hop_5)^p \left( \frac{1}{\alpha} \frac{\partial_{v_i v_j}I^{n+1-p}}{\sqrt{\equilibrium}}\right)  
\\ & \qquad - \sum_{p=1}^{n+1} (\Hop_1 \Hop_5)^{p-1} \Hop_1 \eps^{1-2s} \Delta t (\partial_{x_i} \partial_{v_j} g^{n+1-p} +  \partial_{x_j} \partial_{v_i} g^{n+1-p}).
\end{align}
Similar to \eqref{H2estimate} and \eqref{H3}, we have
\begin{align*}
\norm{\Hop_5}_{L_{x,v}^2 \rightarrow L_{x,v}^2} \leq \frac{\alpha}{\alpha +\gamma-2}.
\end{align*}
Denote $c_2 = \frac{\alpha^2}{(\alpha-\gamma)(\alpha+ \gamma-2)}$. 
Then via iterations, \eqref{g1228} satisfies 
\begin{align*}
\normM{\partial_{v_i v_j} g^\np} &\leq  
 c^{n+1} \normM{\partial_{v_i v_j} g_{in}} 
 +  \eps^{-2s} \Delta t \sum_{p=1}^{n+1} c^p \normM{\partial_{v_i v_j} I^{n+1-p}} 
\\ 
& \quad \,
  + \sum_{p=1}^{n+1} c^{p-1} \frac{\alpha}{\alpha-\gamma} \eps^{1-2s} \Delta t  
\left( \norm{\partial_{x_i } \partial_{v_j}g^{n+1-p}}_{L_{x,v}^2}  +  \norm{\partial_{x_j} \partial_{v_i} g^{n+1-p}}_{L_{x,v}^2}  \right).
\end{align*}
By Lemma~\ref{lemma2} part 1), Lemma~\ref{lem:well-posedness-scheme} and Lemma~\ref{lemma3}, we have 
\begin{align*}
  \normM{\partial_{v_i v_j} g^\np} 
&\leq  
  c^{n+1} \normM{\partial_{v_i v_j} g_{in}} 
  + C \eps^{-2s} \|\rho_{in}\|_{H_x^{3}} \eps^{s} \Delta t 
     \sum_{p=1}^{n+1} c^p  
\\
& \quad \,
  +  C\frac{\alpha}{\alpha-\gamma} \eps^{1-2s} \Delta t  
     (\eps^{-s} + \eps^{1-3s})  \vpran{\sum_{p=1}^{n+1} c^{p-1}} \times
\\
& \qquad \quad\,
    \times \left( \eps^s \normM{\partial_{x_i} g_{in}}+ \eps^s \normM{\partial_{x_j} g_{in}} + \eps^s \normM{\partial_{v_i} g_{in}}  + \norm{\rho_{in}}_{H_x^2} \right).
\end{align*}
Following the same argument as in \eqref{eqn:0806} and \eqref{eqn:0807}, if we chose $\gamma$ such that 
\begin{align} \label{cond2}
0 < \gamma \leq \sqrt{\frac{\lambda_2-4}{\lambda_2}} \frac{\eps^{2s}}{\Delta t} = \sqrt{\frac{\lambda_2-4}{\lambda_2}} \alpha
\end{align}
for some fixed $\lambda_2 >4$, then the above inequality reduces to 
\begin{align*}
& \quad \, 
\normM{\partial_{v_i v_j} g^\np} 
\\ 
&\leq 
  C \normM{\partial_{v_i v_j} g_{in}} + C \eps^{-s} \norm{\rho_{in}}_{H_x^3}
\\ 
& \qquad 
  + C( \eps^{1-3s} + \eps^{2-5s})\left(\eps^s \normM{\partial_{x_i}g_{in}} + \eps^s \normM{\partial_{x_j}g_{in}} + \eps^s \normM{\partial_{v_i}g_{in}} + \norm{\rho_{in}}_{H_x^2} \right) 
\\ 
& \leq 
  C(\eps^{-s} + \eps^{1-3s} + \eps^{2-5s}) \times
\\ 
& \hspace{1cm}
  \times \left(\eps^s \normM{\partial_{x_i}g_{in}} + \eps^s \normM{\partial_{x_j}g_{in}} + \eps^s \normM{\partial_{v_i}g_{in}} + \eps^s \normM{\partial_{v_i v_j} g_{in}}+ \norm{\rho_{in}}_{H_x^2} \right) .
\end{align*}
Additionally, by \eqref{g1229} one has that
\begin{align*}
\partial_{v_i v_j} g^\nh = (\alpha - \Lop^s + \gamma - 2)^{-1} (\alpha \partial_{v_i v_j} g^n - \partial_{v_i v_j} I^\np).
\end{align*}
This implies
\begin{align*}
\frac{\partial_{v_i v_j} g^\nh}{\sqrt{\equilibrium}}  = \Hop_4 \left( \frac{\partial_{v_i v_j} g^n}{\sqrt{\equilibrium}}  - \frac{1}{\alpha} \frac{1}{\sqrt{\equilibrium}} \partial_{v_i v_j} I^\np \right).
\end{align*}
Therefore, for $\gamma >2$, we have 
\begin{align*}
\normM{\partial_{v_i} g^\nh} \leq C(\normM{\partial_{v_i v_j} g^n} + \frac{1}{\alpha} \normM{\partial_{v_i v_j} I})
\end{align*}
and the desired bound follows from the estimate for $\normM{\partial_{v_iv_j} g^n}$.
\end{proof}

\begin{remark}
Note that the condition \eqref{cond2} for $\gamma$ already encompasses condition \eqref{cond1}. Consequently, \eqref{cond1} is not needed in the statement of the Lemma.
\end{remark}

\bibliography{reference}

\begin{thebibliography}{10}

\bibitem{abdallah2011fractional}
{\sc N.~B. Abdallah, A.~Mellet, and M.~Puel}, {\em Fractional diffusion limit
  for collisional kinetic equations: a {H}ilbert expansion approach}, Kinetic
  \& Related Models, 4 (2011), p.~873.

\bibitem{aceves2019fractional}
{\sc P.~Aceves-Sanchez and L.~Cesbron}, {\em {F}ractional diffusion limit for a
  fractional {V}lasov--{F}okker--{P}lanck equation}, SIAM Journal on
  Mathematical Analysis, 51 (2019), pp.~469--488.

\bibitem{bessemoulin2020hypocoercivity}
{\sc M.~Bessemoulin-Chatard, M.~Herda, and T.~Rey}, {\em {H}ypocoercivity and
  diffusion limit of a finite volume scheme for linear kinetic equations},
  Mathematics of Computation, 89 (2020), pp.~1093--1133.

\bibitem{cesbron2012anomalous}
{\sc L.~Cesbron, A.~Mellet, and K.~Trivisa}, {\em Anomalous transport of
  particles in plasma physics}, Applied Mathematics Letters, 25 (2012),
  pp.~2344--2348.

\bibitem{crouseilles2016numerical}
{\sc N.~Crouseilles, H.~Hivert, and M.~Lemou}, {\em Numerical schemes for
  kinetic equations in the anomalous diffusion limit. part {I}: The case of
  heavy-tailed equilibrium}, SIAM Journal on Scientific Computing, 38 (2016),
  pp.~A737--A764.

\bibitem{CHL16}
\leavevmode\vrule height 2pt depth -1.6pt width 23pt, {\em Numerical schemes
  for kinetic equations in the anomalous diffusion limit. part {II}: Degenerate
  collision frequency}, SIAM Journal on Scientific Computing, 38 (2016),
  pp.~A2464--A2491.

\bibitem{Degond1986}
{\sc P.~Degond}, {\em Global existence of smooth solutions for the
  {V}lasov-{F}okker-{P}lanck equation in 1 and 2 space dimensions.}, Ann. Sci.
  de l'\'{E}.N.S., 19 (1986), pp.~519--542.

\bibitem{doetsch2012introduction}
{\sc G.~Doetsch}, {\em {I}ntroduction to the Theory and Application of the
  Laplace Transformation}, Springer Science \& Business Media, 2012.

\bibitem{filbet2013analysis}
{\sc F.~Filbet and A.~Rambaud}, {\em {A}nalysis of an asymptotic preserving
  scheme for relaxation systems}, ESAIM: Mathematical Modelling and Numerical
  Analysis, 47 (2013), pp.~609--633.

\bibitem{filbet2021convergence}
{\sc F.~Filbet, L.~M. Rodrigues, and H.~Zakerzadeh}, {\em {C}onvergence
  analysis of asymptotic preserving schemes for strongly magnetized plasmas},
  Numerische Mathematik, 149 (2021), pp.~549--593.

\bibitem{gentil2008levy}
{\sc I.~Gentil and C.~Imbert}, {\em The {L}{\'e}vy--{F}okker--{P}lanck
  equation: $\phi$-entropies and convergence to equilibrium}, Asymptotic
  Analysis, 59 (2008), pp.~125--138.

\bibitem{golse1999convergence}
{\sc F.~Golse, S.~Jin, and C.~D. Levermore}, {\em The convergence of numerical
  transfer schemes in diffusive regimes {I}: {D}iscrete-ordinate method}, SIAM
  journal on numerical analysis, 36 (1999), pp.~1333--1369.

\bibitem{gosse2004asymptotic}
{\sc L.~Gosse and G.~Toscani}, {\em Asymptotic-preserving \& well-balanced
  schemes for radiative transfer and the {R}osseland approximation}, Numerische
  Mathematik, 98 (2004), pp.~223--250.

\bibitem{hawkes1971lower}
{\sc J.~Hawkes}, {\em {A} lower {L}ipschitz condition for the stable
  subordinator},  (1971).

\bibitem{hu2021uniform}
{\sc J.~Hu and R.~Shu}, {\em On the uniform accuracy of implicit-explicit
  backward differentiation formulas ({IMEX-BDF}) for stiff hyperbolic
  relaxation systems and kinetic equations}, Mathematics of Computation, 90
  (2021), pp.~641--670.

\bibitem{jang2014analysis}
{\sc J.~Jang, F.~Li, J.-M. Qiu, and T.~Xiong}, {\em {A}nalysis of asymptotic
  preserving {DG-IMEX} schemes for linear kinetic transport equations in a
  diffusive scaling}, SIAM Journal on Numerical Analysis, 52 (2014),
  pp.~2048--2072.

\bibitem{jin2022asymptotic}
{\sc S.~Jin}, {\em Asymptotic-preserving schemes for multiscale physical
  problems}, Acta Numerica, 31 (2022), pp.~415--489.

\bibitem{jin2009uniformly}
{\sc S.~Jin, M.~Tang, and H.~Han}, {\em A uniformly second order numerical
  method for the one-dimensional discrete-ordinate transport equation and its
  diffusion limit with interface}, Networks \& Heterogeneous Media, 4 (2009),
  p.~35.

\bibitem{klar2002uniform}
{\sc A.~Klar and A.~Unterreiter}, {\em Uniform stability of a finite difference
  scheme for transport equations in diffusive regimes}, SIAM journal on
  numerical analysis, 40 (2002), pp.~891--913.

\bibitem{li2017uniform}
{\sc Q.~Li and L.~Wang}, {\em {U}niform regularity for linear kinetic equations
  with random input based on hypocoercivity}, SIAM/ASA Journal on Uncertainty
  Quantification, 5 (2017), pp.~1193--1219.

\bibitem{JLL1961}
{\sc J.~Lions}, {\em Equations diff\'erentielles op\'erationnelles et
  probl\`emes aux limites}, Springer, Berlin (1961).

\bibitem{liu2010analysis}
{\sc J.-G. Liu and L.~Mieussens}, {\em Analysis of an asymptotic preserving
  scheme for linear kinetic equations in the diffusion limit}, SIAM Journal on
  Numerical Analysis, 48 (2010), pp.~1474--1491.

\bibitem{mellet2011fractional}
{\sc A.~Mellet, S.~Mischler, and C.~Mouhot}, {\em Fractional diffusion limit
  for collisional kinetic equations}, Archive for Rational Mechanics and
  Analysis, 199 (2011), pp.~493--525.

\bibitem{peng2021stability}
{\sc Z.~Peng, Y.~Cheng, J.-M. Qiu, and F.~Li}, {\em {S}tability-enhanced {AP}
  {IMEX1-LDG} method: energy-based stability and rigorous {AP} property}, SIAM
  Journal on Numerical Analysis, 59 (2021), pp.~925--954.

\bibitem{uchaikin2011chance}
{\sc V.~V. Uchaikin and V.~M. Zolotarev}, {\em {C}hance and stability: stable
  distributions and their applications}, Walter de Gruyter, 2011.

\bibitem{wang2016asymptotic}
{\sc L.~Wang and B.~Yan}, {\em An asymptotic-preserving scheme for linear
  kinetic equation with fractional diffusion limit}, Journal of Computational
  Physics, 312 (2016), pp.~157--174.

\bibitem{wang2019asymptotic}
\leavevmode\vrule height 2pt depth -1.6pt width 23pt, {\em An
  asymptotic-preserving scheme for the kinetic equation with anisotropic
  scattering: Heavy tail equilibrium and degenerate collision frequency}, SIAM
  Journal on Scientific Computing, 41 (2019), pp.~A422--A451.

\bibitem{XW21}
{\sc W.~Xu and L.~Wang}, {\em An asymptotic preserving scheme for
  {L}\'{e}vy-{F}okker-{P}lanck equation with fractional diffusion limit}, Comm.
  Math. Sci. (arXiv:2103.08848),  (accepted).

\end{thebibliography}
\bibliographystyle{siam}

\end{document}